\documentclass[12pt,english]{article}
\usepackage[T1]{fontenc}
\usepackage[utf8]{inputenc}
\usepackage{geometry}
\usepackage{babel}
\usepackage{float}
\usepackage{graphicx} 
\usepackage{color} 
\usepackage{hyperref}
\usepackage{subcaption}
\usepackage{amsmath,amssymb,amsthm,bm}
\usepackage{authblk}

\usepackage{caption}
\usepackage{placeins}

\newcommand{\R}{\mathbb{R}}

\usepackage{booktabs} 
\newcommand{\br}{\toprule[1.5pt]} 
\newcommand{\mr}{\midrule}
\usepackage{comment}

\title{Inverse Electromagnetic Scattering for Doubly-Connected Cylinders using Convolutional Neural Networks
}

\author[1]{Leonidas Mindrinos\thanks{leonidas.mindrinos@aua.gr, ORCID: \href{https://orcid.org/0000-0002-7807-8153}{0000-0002-7807-8153}}}
\affil[1]{Department of Natural Resources Development and Agricultural Engineering, Agricultural University of Athens, 11855, Greece}

\author[2]{Nikolaos Pallikarakis\thanks{npall@central.ntua.gr, ORCID: \href{https://orcid.org/0000-0002-2385-6747}{0000-0002-2385-6747}, Corresponding author}}
\affil[2]{Department of Mathematics, National Technical University of Athens, Zografou Campus, Athens, 15780, Greece}

\author[3]{Nikolaos L Tsitsas\thanks{ntsitsas@csd.auth.gr, ORCID: \href{https://orcid.org/0000-0003-1409-2631}{0000-0003-1409-2631}}}
\affil[3]{School of Informatics, Aristotle University of Thessaloniki, Thessaloniki, 54124, Greece}

\date{}

\begin{document}

\maketitle

\begin{abstract}
In this work, we consider the inverse electromagnetic scattering problem for a magneto-dielectric cylinder covering an impedance cylinder of arbitrary shape. We solve it by introducing a divide-and-conquer framework using specially designed 1D multi-channel, circular-padding Convolutional Neural Networks. The solution of the direct problem provides us with the real and imaginary components of the far-field measurements representing the input data. We first classify the shape of the impedance cylinder and then reconstruct the unknown boundary curve and the impedance function. Through extensive numerical experiments, including noisy scenarios, we demonstrate the efficiency and robustness of our approach. 
\end{abstract}

\noindent\textbf{Keywords}: 
oblique incidence; inverse problem; electromagnetic scattering; convolutional neural networks; deep learning\\
\textbf{ 2020 Mathematics Subject Classification:} 68T07, 78A46, 65N21, 35R30

\section{Introduction}

Direct and inverse electromagnetic scattering theory studies the interaction of electromagnetic waves with objects. In direct scattering problems, the scattered fields are derived from known object properties by solving  Maxwell's equations together with appropriate boundary conditions. In inverse scattering, the goal is to determine the unknown properties and/or the shape of an object from scattered field data in the near- or far-field regime \cite{CKbook, KHbook}.  Such an inverse problem is challenging due to its ill-posed nature with respect to stability, where small data errors can lead to large inaccuracies in the solution. Despite these difficulties, inverse scattering has important applications in medical imaging, geophysical exploration, and non-destructive testing.

More specifically, in this work we consider the scattering problem of time-harmonic electromagnetic waves by an infinitely long magneto-dielectric cylinder. This setup finds applications in medical imaging \cite{medim}, remote sensing \cite{remsen} and antenna design \cite{antdes}. For oblique incident waves, the problem's complexity increases, as the scattering depends on the incident angle, requiring advanced techniques; see e.g. \cite{tsitsas1,Yousif:88}. Nevertheless, due to the underlying symmetry, the three-dimensional problem can be reduced to a system of two-dimensional Helmholtz equations, simplifying the analysis \cite{NakW}.

The numerical solution of the inverse problem to determine the shape and/or the properties of an object from scattering data is a broadly studied subject. Methods like the Linear Sampling Method \cite{cakoni2011linear}, the Factorization method \cite{kir04} and the Enclosure method \cite{ike16} have been effectively used to handle such problems in electromagnetics. Another commonly used method is the Boundary Integral Equation (BIE) method \cite{KreRun05}  with recent applications in the study of the inverse electromagnetic scattering by infinite long dielectric cylinders at oblique incidence; see, e.g., \cite{GinMin19, WANG2012860}. 

In particular, the BIE method is considered as a powerful technique for reconstructing the boundary of a homogeneous single- or doubly-connected scatterer with the advantage of reducing the dimensions of the problem. However, the method has two notable limitations. First, its performance strongly depends on the initial guess for the boundary curve. Second, since the boundary’s radial function is typically approximated by a series expansion (often a trigonometric polynomial), the required series order depends on prior knowledge of the boundary’s geometric complexity. For example, a peanut-shaped boundary can be accurately reconstructed using fewer coefficients than a kite-shaped one \cite{cakoni14, cakkre10, fink, KreRun01, min23, ivan16}.

Recently, machine learning, particularly deep neural networks (DNNs), has been applied to solve inverse problems, addressing their ill-posed nature by utilizing the ability of DNNs to learn complex nonlinear relationships \cite{AO, arridge2024deep, AMOS, LSAH}. These type of models have been successfully used in various inverse problems, see e.g., \cite{GLWZ,JMFU,PN}. Additionally, the so-called ``divide-and-conquer$"$ algorithms have been explored in several studies on inverse scattering \cite{ChengGuo,DLLLS, MENG2025116525, pallMLP, YIN2020109594}. In these approaches, a classification algorithm is first employed to assign the input to a specific category, followed by solving a separate regression problem within each class.

The present work introduces a novel one-dimensional Convolutional Neural Network (CNN) methodology that serves as an alternative to the conventional BIE approach for boundary reconstruction. The network is multi-channel, i.e., it jointly processes multiple electromagnetic far-field measurements as separate input channels, enabling the model to learn cross-correlations between field components at each measurement position. We also apply circular padding within the CNN to preserve the $2\pi$-periodicity of the measurements. The algorithm  classifies first the impedance cylinders according to their boundary shapes and then reconstructs the unknown boundary, thereby solving the inverse scattering problem. Unlike previous studies that focus on simply connected domains, we address the more challenging scenario of doubly-connected geometries. Furthermore, we extend beyond pure shape reconstruction by treating the impedance boundary function as an additional unknown parameter, significantly increasing the problem complexity and demonstrating the capability of neural networks to tackle inverse problems with multiple coupled unknowns.

For the needs of this study, we solve the corresponding direct scattering problems using BIE to generate the training data. Specifically, we compute the complex-valued far-field patterns--both real and imaginary components of the electric and magnetic fields--sampled in the unit circle. Our goal is to minimize data requirements by using few measurement angles and the smallest feasible dataset, depending on the problem's complexity. We restrict ourselves to measurements from a single incident wave, introducing a second incident direction only when necessary to resolve the increased inverse problem challenges.

The structure of the paper is as follows. Section \ref{scatt} formulates the scattering problem, including the definitions of the direct and inverse problems, and outlines the numerical implementation used for solving the direct scattering problem. Section \ref{CNNmeth} presents the CNN methodology, including a rigorous mathematical formulation of the end-to-end network as a composition of operators that approximate the inverse map.  Section \ref{inverse} reports the numerical results for the inverse problems, covering both classification and regression tasks. We conclude and outline future work in Section \ref{concl}. Computational details have been placed to Appendix \ref{A}, for better readability and to streamline the presentation of the main numerical findings. To ensure reproducibility and transparency, all numerical experiments include 
complete architectural specifications, training dynamics, comprehensive error analyses, 
and detailed dataset information.

\section{The scattering problem} \label{scatt}

\subsection{Formulation of the problem}\label{sec_formulation}

We study the scattering of a time-harmonic electromagnetic wave $(\mathbf{E}_{inc}, \mathbf{H}_{inc})$ incident obliquely, with direction $\hat{\bm d}_{inc},$ on an infinitely long cylinder oriented along the \(z\)-axis. The interior of the cylinder contains an arbitrarily shaped and infinitely long impedance cylinder denoted by \(\Omega_{imp} \), coated by a homogeneous magneto-dielectric layer which is denoted by \(\Omega_{int} \subset \mathbb{R}^3\). This layer is characterized by permittivity \(\varepsilon_1\) and permeability \(\mu_1\). The exterior domain, \(\Omega_{ext} = \mathbb{R}^3 \setminus (\overline{\Omega_{int}}\cup \overline{\Omega_{imp}})\), is described by the parameters \(\varepsilon_0\) and \(\mu_0\).

\begin{figure}[t!]
\centering
\includegraphics[width=0.6\textwidth]{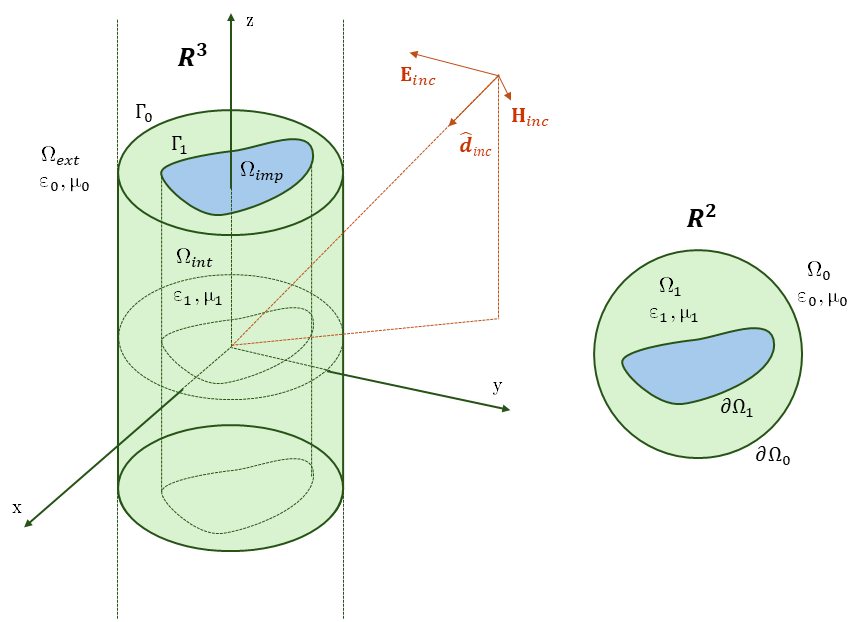}
        \label{fig:second}
    \caption{Electromagnetic scattering by a doubly-connected magneto-dielectric cylinder. 3D (left) and cross-section (right).}
    \label{fig1}
\end{figure}

The boundary of the cylinder, \(\Gamma = \Gamma_1 \cup \Gamma_0\) consists of  two disjoint surfaces: the interior boundary \(\Gamma_1\) and the exterior boundary \(\Gamma_0\). An illustration is shown in Figure \ref{fig1}, with the cross-section of this cylinder being a two-dimensional, doubly-connected domain.
This formulation can be used to model, for example, coated antennas or tubes. The scattering problem is governed by Maxwell's equations. The exterior fields \(\mathbf{E}_{ext}, \mathbf{H}_{ext} : \Omega_{ext} \to \mathbb{C}^3\) and the interior fields \(\mathbf{E}_{int}, \mathbf{H}_{int} : \Omega_{int} \to \mathbb{C}^3\) satisfy (assuming time dependence of the form $e^{-\mathrm{i}\omega t},$ where \(\omega > 0\) is the angular frequency):
\begin{equation}\label{eq_maxwell}
\begin{aligned}
\nabla \times \mathbf{E}_{ext} - \mathrm{i}\omega \mu_0 \mathbf{H}_{ext} &= 0, \quad 
\nabla \times \mathbf{H}_{ext} + \mathrm{i}\omega \varepsilon_0 \mathbf{E}_{ext} = 0, \quad \text{in } \Omega_{ext},   \\
\nabla \times \mathbf{E}_{int} - \mathrm{i}\omega \mu_1 \mathbf{H}_{int} &= 0, \quad 
\nabla \times \mathbf{H}_{int} + \mathrm{i}\omega \varepsilon_1 \mathbf{E}_{int} = 0, \quad \text{in } \Omega_{int}.
\end{aligned}
\end{equation}

The impedance cylinder is surrounded by the boundary \(\Gamma_1\) where we impose the Leontovich impedance boundary condition (IBC):
\begin{equation}\label{eq_impe}
(\hat{\mathbf{n}} \times \mathbf{E}_{int}) \times \hat{\mathbf{n}} = \lambda (\hat{\mathbf{n}} \times \mathbf{H}_{int}), \quad \text{on } \Gamma_1,    
\end{equation}
where \(\hat{\mathbf{n}}\) is the outward unit normal vector and \(\lambda \in C^1(\Gamma_1)\) is the impedance function. This boundary condition approximates the properties of lossy or dielectric-coated materials and describes electromagnetic fields that do not penetrate the region \(\Omega_{imp}\), see e.g., \cite{Alex, Tezel}. 

The cylinder is assumed to be penetrable, with the transmission conditions applied on the exterior boundary:
\begin{equation}\label{eq_trans}
\hat{\mathbf{n}} \times \mathbf{E}_{int} = \hat{\mathbf{n}} \times \mathbf{E}_{ext}, \quad 
\hat{\mathbf{n}} \times \mathbf{H}_{int} = \hat{\mathbf{n}}\times \mathbf{H}_{ext}, \quad \text{on } \Gamma_0.   
\end{equation}

The direct problem involves computing the scattered and interior fields given the incident wave, object geometry, and material properties. In contrast, the inverse problem, considered in this work, aims to reconstruct the impedance cylinder, characterized by $\Gamma_1$ and $\lambda$, from far-field data corresponding to one (in most of the cases) incident wave, given the material parameters. This framework has been used to analyze scattering in a doubly-connected domain for both the direct problem \cite{min19} and the inverse problem \cite{min23}.

\subsection{Numerical solution of the direct problem} \label{direct}

Given the cylindrical symmetry of the scatterer, the three-dimensional problem \eqref{eq_maxwell}--\eqref{eq_trans} can be reduced to a two-dimensional one involving only the $z$-components of the electric and magnetic fields~\cite{NakW}. More precisely, Maxwell's equations \eqref{eq_maxwell} are reformulated as a system of four Helmholtz equations (two in each domain). Furthermore, each boundary condition in \eqref{eq_impe}--\eqref{eq_trans} is correspondingly transformed into two scalar equations for the field components, involving both normal and tangential derivatives, formulated on $\partial\Omega_1$ and $\partial\Omega_0,$ the cross-sections of $\Gamma_1$ and $\Gamma_0,$ respectively. From this point onward, we refer to the cross-sectional profile of the impedance cylinder $\Omega_{imp}$ as the obstacle, whose boundary curve corresponds to $\partial\Omega_1$ (see Figure \ref{fig1}, right).

The incident wave admits the form
\[
e_{inc}(\hat{\bm d}; \bm x) = \frac{\sin \theta}{\sqrt{\varepsilon_0}}  e^{\mathrm{i} \kappa_0 \hat{\bm d} \cdot \bm x}, \quad h_{inc} (\hat{\bm d}; \bm x)=0, \quad \bm x \in \Omega_0 \subset \mathbb{R}^2,
\]
for $\theta \in (0,\pi),$ the wavenumber $\kappa^2_0 = \omega^2 \mu_0 \varepsilon_0 (1-\cos^2 \theta)$ and the incidence direction 
$\hat{\bm d} = ( \cos \phi,\, \sin \phi), \, \phi \in [0,2\pi].$

The exterior boundary is known, a circle with a given radius, with parametrization
\[
\partial \Omega_0 = \{\bm y(\tau) = 0.8 (\cos \tau, \, \sin \tau): \tau \in [0,2\pi] \}.
\]
For the interior boundary curve, we consider the following three types (classes) for $\tau\in[0,2\pi],$ centered around $\bm x_0 \in [-0.2,\, 0.2]^2 : $
\begin{itemize}
\item Peanut-shaped boundary: 
\begin{equation}
\bm{x}(\tau) = \rho(\tau) (\cos \tau, \sin \tau) + \bm{x}_0, \quad \rho(\tau) = \sqrt{\alpha \cos^2 \tau + \beta \sin^2 \tau}. \label{pean}
\end{equation}
\item Kite-shaped boundary: 
\begin{equation}
\bm{x}(\tau) = (\alpha \cos \tau + \beta \cos(\tau), \gamma \sin \tau) + \bm{x}_0. \label{kite}
\end{equation}
\item Star-shaped boundary: 
\begin{equation}
\bm{x}(\tau) = \rho(\tau) (\cos \tau, \sin \tau) + \bm{x}_0, \quad \rho(\tau) = \alpha_0 \left\{ 1 + \frac{1}{2Q} \sum_{q=1}^{Q} \left[ \alpha_q \cos(q\tau) + \beta_q \sin(\tau) \right] \right\}. \label{star}
\end{equation}
\end{itemize}
All the coefficients mentioned in the above boundary curve representations are random numbers, generated using the MATLAB function  $\texttt{rand},$ and are uniformly distributed within appropriate intervals to ensure that the interior boundary does not intersect with the exterior boundary. We consider equidistant grid points defined by $$\tau_k = \frac{ 2 \pi k}{\mathcal{T}}, \, k=0,...,\mathcal{T}-1,$$  on both boundary curves and we fix $\mathcal{T}=128$ in all experiments. 

 We set $\omega=5$ the angular frequency,  $\theta = \tfrac{\pi}{6},$ and $\phi \in \{ 0, \pi \}.$ The material parameters are given by $(\varepsilon_1, \, \mu_1) = (2,1)$ and $(\varepsilon_0, \, \mu_0) = (1,1)$ in the interior and exterior domain, respectively. The impedance function defined in \eqref{eq_impe} is a variable $\lambda\in[0.1,10]$.  For the star-shaped boundary we fix $Q=5$ (i.e. 11 total boundary parameters).

We solve the direct problem, that is, to compute the far-field patterns $(e^\infty(\boldsymbol{\hat{x}}),\,h^\infty(\boldsymbol{\hat{x}}))$ for $\boldsymbol{\hat{x}} = (\cos t, \sin t)$, $t \in [0,2\pi],$ of the scattered fields (satisfying radiation conditions) using the BIE formulation. This method was presented in~\cite{min19} and originally introduced in~\cite{KreRun01}. In \cite{min19}, a hybrid BIE method was proposed, in which the direct approach (based on Green’s formulas) was applied to the exterior fields, while a single-layer ansatz was employed for the interior fields.

The unit circle, where the far-field data are measured, is parametrized through
\[
t_k = \frac{ 2 \pi k}{T_0}, \, k=0,...,T_0-1,
\]
where we consider two cases: $T_0=128$ (dense data) and $T_0=32$ (sparse data), corresponding to the number of measurement angles. We assemble the data in the form
\[
\bigl(\Re\{e^\infty\}, \, \Im\{e^\infty\}, \, \Re\{h^\infty\}, \, \Im\{h^\infty\}\bigr) \in \R^{4T_0},
\]
if one incident field is used $(\phi =0),$ or
\[
\bigl(\Re\{e^\infty_0\}, \, \Im\{e^\infty_0\}, \, \Re\{h^\infty_0\}, \, \Im\{h^\infty_0\}, \Re\{e^\infty_{\pi}\}, \, \Im\{e^\infty_{\pi}\}, \, \Re\{h^\infty_{\pi}\}, \, \Im\{h^\infty_{\pi}\}\bigr) \in \R^{8T_0},
\]
for two $(\phi=0\ \textrm{and}\ \pi)$, where the subscript indicates dependence on the incidence direction. In some of the inverse problems considered in this study, information from the electric field alone is sufficient and the limited case
\[
\bigl(\Re\{e^\infty\}, \, \Im\{e^\infty\}\bigr) \in \R^{2T_0},
\]
will be considered. 

In what follows, the parameter $C_0 \in \{2,4,8\}$ denotes the input channel dimension,  which specifies the number of field components measured at each angular position and  determines the dimension of the $C_0\cdot T_0-$dimensional space.  The far-field data will be used to train the CNN models for solving the inverse problems that follow.

\subsection{Discussion on the inverse problem} 

The unique solvability of the inverse problem is a theoretical open problem. A relevant uniqueness result for the case of an impedance cylinder can be found in \cite{nak}. Moreover, the numerical solution of the inverse problem for reconstructing the interior boundary poses challenges even when the impedance $\lambda$ is known, as demonstrated in \cite{min23}. These difficulties arise from the oblique scattering conditions and the inherently ill-posed nature of the problem. To underscore this point, a successful reconstruction of the cross-section of the cylinder could indicate defects in the medium. However, since no theoretical result exists, there is no guarantee that the knowledge of the far-field data suffices for an accurate reconstruction. 

\begin{figure}[t!]
    \centering
    \begin{minipage}{0.42\textwidth}
        \centering
        \includegraphics[width=\textwidth]{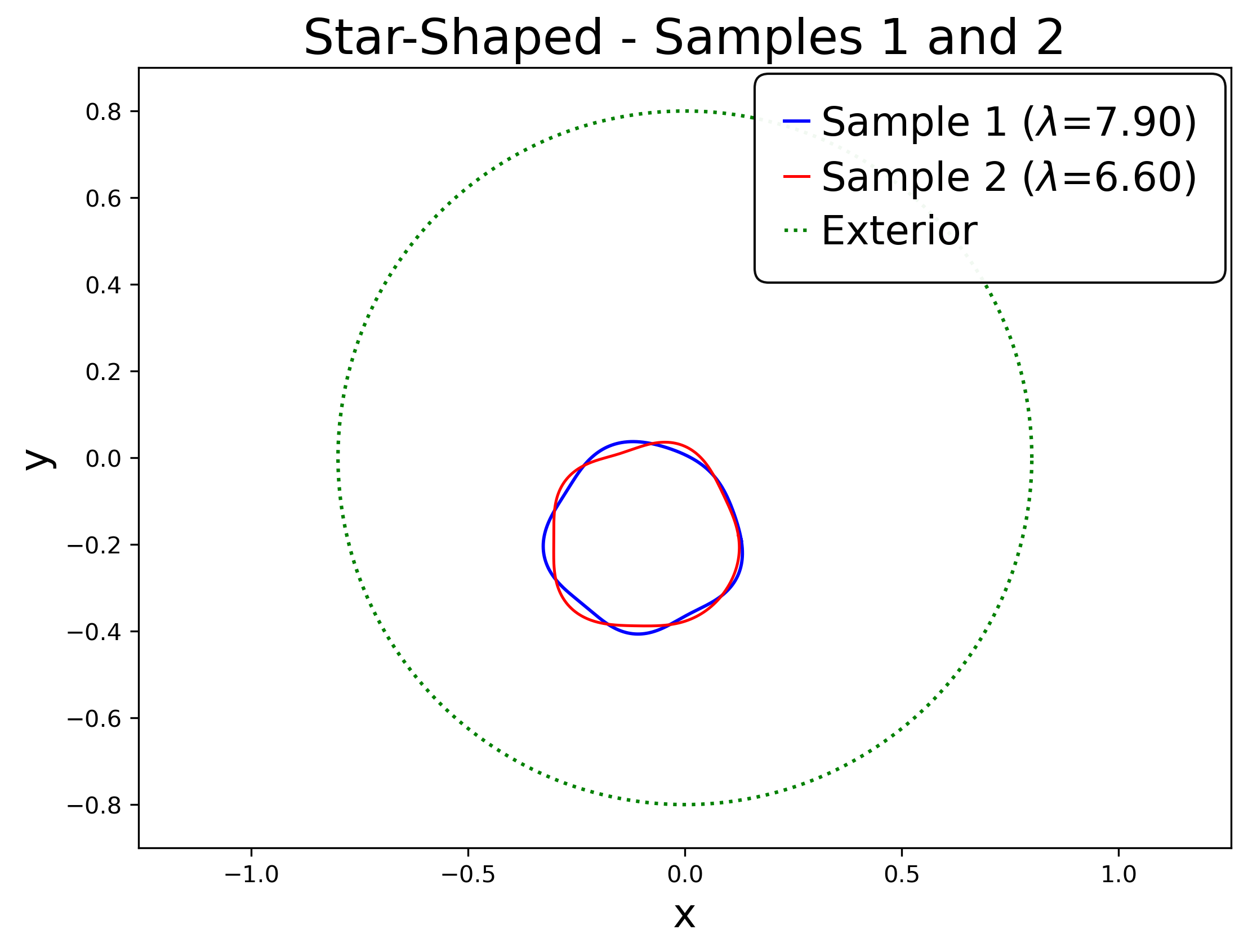}
    \end{minipage}
    \hfill
    \begin{minipage}{0.56\textwidth}
        \centering
        \includegraphics[width=\textwidth]{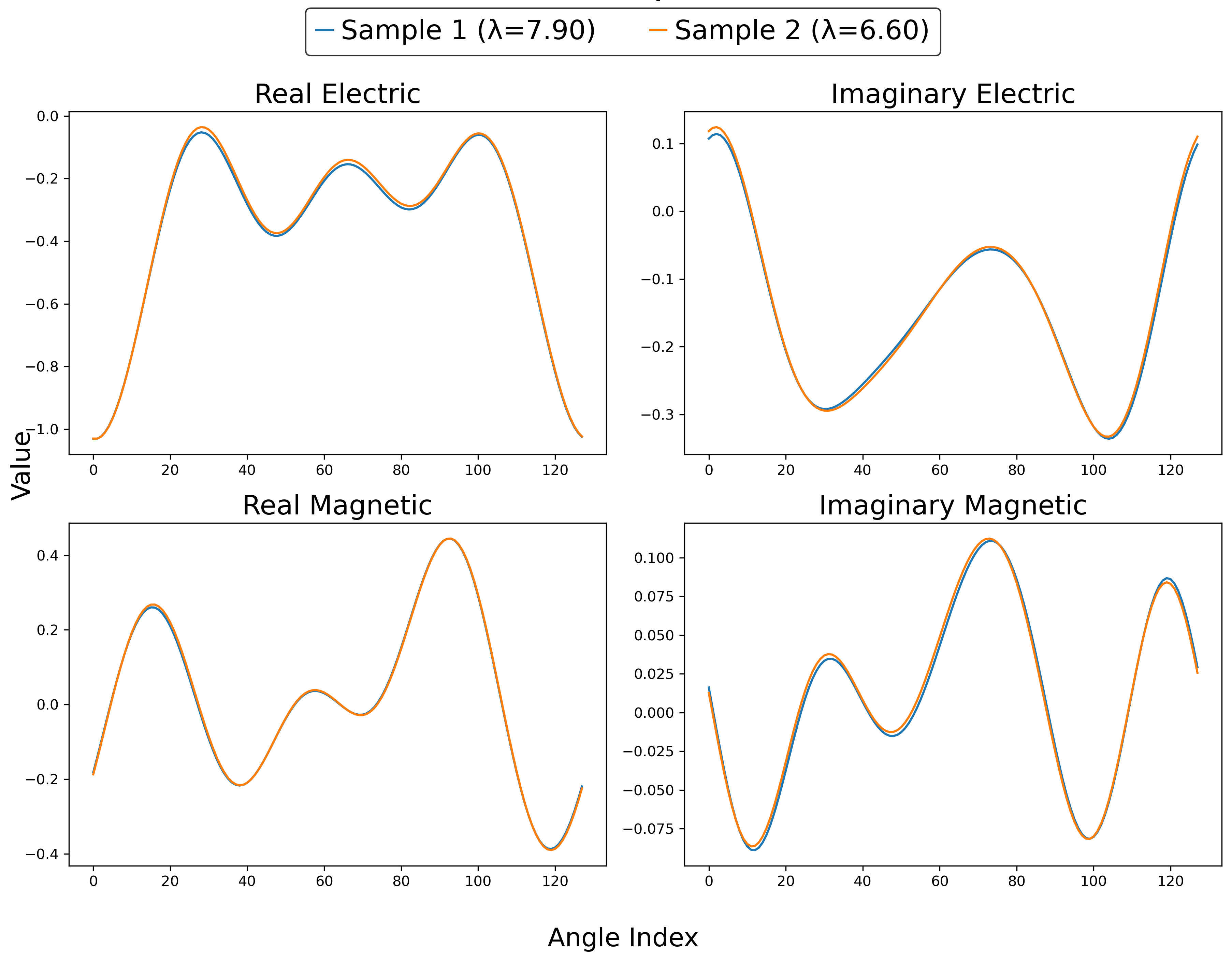}
        \label{fig3}
    \end{minipage}
    \caption{Representative example of the ill-posedness of the inverse problem: two different star-shaped obstacles (left) produce nearly identical far-field data (right).}
    \label{fig2}
\end{figure}

To illustrate these challenges, we present the following example. We examine two different star-shaped obstacles (Sample 1 and 2) with similar centers but distinct boundaries, where their constant impedance values $\lambda$ also differ, see Figure \ref{fig2} (left). As shown in Figure \ref{fig2} (right), the far-field data measured at 128 angles outisde the external cylinder demonstrate that both the real and imaginary components of the electric and magnetic fields are remarkably similar between the two samples. Consequently, solving the inverse problem--namely, recovering both the inner boundary curve and the impedance function from the far-field measurements--becomes a challenging and potentially ill-posed process.

In this work, we adopt a divide-and-conquer strategy to solve the inverse electromagnetic scattering problem. The approach begins with a classification/obstacle perception stage, where scatterers are categorized into three (but not limited to) distinct classes based on their boundary shapes: peanut, kite and star-shaped. Following classification, a regression/shape recognition model is applied within each class to reconstruct the coefficients appearing in equations (\ref{pean})–(\ref{star}) that characterize the unknown boundary shape, together with the unknown impedance $\lambda$. This two-step process enables a more targeted and efficient solution, thereby addressing the complexity of the inverse problem. We also examine the capability of the proposed method to reconstruct the boundary of a misclassified object, demonstrating its robustness.

\section{Convolutional Neural Network Methodology} \label{CNNmeth}

\begin{figure}[t!]
\centering
\includegraphics[width=1\textwidth]{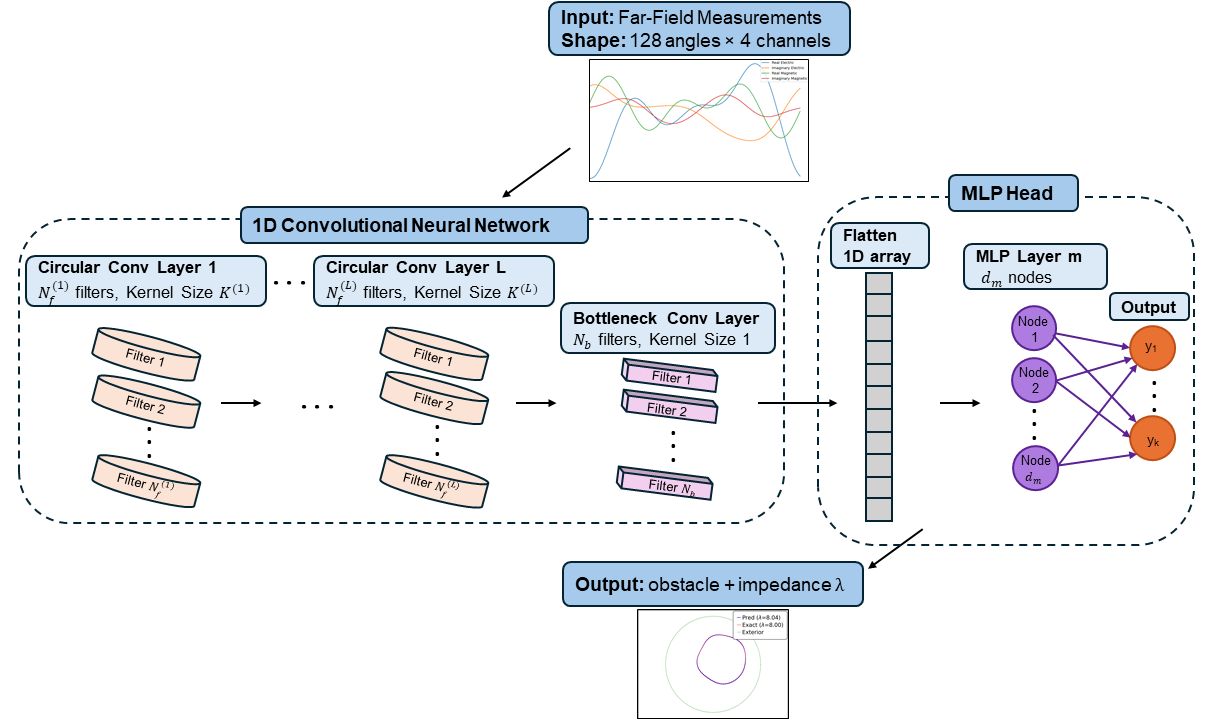}
\caption{Architecture of the circular padding CNN and MLP inversion head, forming an end‐to‐end network from raw input measurements to final outputs.}
\label{fig4}
\end{figure}

One-dimensional CNNs are widely used for processing 1D signals; we refer to the survey article \cite{1dcnn} for a comprehensive review. In this work, we propose a specialized 1D CNN to handle the inherent periodicity of angular electromagnetic measurements. Since the far-field data are collected at cyclical positions outside  the magneto-dielectric cylinder, the signal values repeat after a full $2\pi$ rotation. The architecture uses circular padding and multiscale feature extraction to preserve this $2\pi$–periodic nature of the far-field data while capturing both local and global angular dependencies. Additionally, the CNN processes the real and imaginary components of the electric and magnetic fields as separate input channels, allowing the convolutional filters to learn their inter-relationships at each angular position. After the CNN layers, a Multi-Layer Perceptron (MLP) head is appended as a nonlinear mapping from the extracted feature vector to the regression or classification output. We refer to the supplementary material of \cite{pallMLP} for a detailed review of the MLP algorithm. While circular CNNs have been applied to 2D image data \cite{schubert2019circular}, and to classification of periodic 1D sequential data \cite{cheng2022classification}, our approach adapts these principles specifically for 1D angular electromagnetic measurements. Figure \ref{fig4} shows the overall network architecture, and the details of each component are explained in the following subsections.

The framework presented herein formulates the end-to-end network as a sequence of well‐defined linear and nonlinear operators, thereby explicitly embedding the inverse problem structure into the network architecture and training objective.

\subsection{Problem Setting and Operator Formulation} \label{meth_set}

Let
\[
\mathcal{X} = \mathbb{R}^{T_0\times C_0}
\]
be the Euclidean (Hilbert) space of channel-wise angular measurements, where $T_0$ is 
the number of angular sampling points and $C_0$ the number of input field components. We equip $\mathcal{X}$ with the inner product
\[
\langle \mathbf{X},\mathbf{Y}\rangle
= \sum_{i=0}^{T_0-1}\sum_{c=0}^{C_0-1} X_{i,c}\,Y_{i,c},
\]
where the entries are real-valued. 

Denote the forward scattering operator
\[
\mathcal{F}: S \longrightarrow \mathcal{X},
\]
which maps an obstacle parameter vector \(\mathbf{s}\in S\subseteq\mathbb{R}^P\) (e.g.\ boundary curve coefficients and impedance) to far-field measurements \(\mathbf{X}=\mathcal{F}(\mathbf{s})\). The inverse problem is to obtain an estimate
\[
\mathbf{s}^\star \approx \mathcal{F}^{-1}(\mathbf{X}),\qquad \mathbf{X}\in\mathcal{X},
\]
which is generally ill-posed and nonlinear.

We approximate \(\mathcal{F}^{-1}\) by a parameterized operator (the end-to-end network)
\[
\mathcal{N}_{\boldsymbol{\theta}}:\mathcal{X}\longrightarrow\widehat{S}=\mathbb{R}^{d_{\mathrm{out}}},
\]
with learnable parameters \(\boldsymbol{\theta}\). Here
\[
  d_{\mathrm{out}}
  = 
  \begin{cases}
  K_{\mathrm{cls}}, & \text{classification (number of classes),}\\
    P, & \text{regression (number of boundary and impedance parameters).}    
  \end{cases}
\]
Training enforces \(\mathcal{N}_{\boldsymbol{\theta}}\approx\mathcal{F}^{-1}\).

\subsection{Input Data Representation} \label{meth_input}

Define the $n$-th training sample as a vector
\[
\mathbf{x}_n\in\mathbb{R}^{C_0 \cdot T_0},
\]
and the canonical isomorphism
\[
\mathrm{reshape}:\mathbb{R}^{C_0\cdot T_0}\;\xrightarrow{\sim}\;\mathbb{R}^{T_0\times C_0},
\quad
\mathbf{x}_n\mapsto \mathbf{X}_n,
\]
so that $X_{n; i,c}$ is the $c$-th channel at angular index $i$ of the sample $n$. For example, if $C_0=4$ (real/imaginary electric and magnetic fields) measured at $T_0 = 128$ positions, then $\mathbf{x}_n\in\mathbb{R}^{512}$ and $\mathbf{X}_n\in\mathbb{R}^{128\times 4}$.

This reshaping is essential to treat the angular and channel axes separately, allowing 1D convolutions to slide kernels along the angular dimension without mixing channels arbitrarily, as we explain in the following.  From a physical perspective, this multi-channel approach is important because the different fields encode distinct information about wave amplitude, phase, and interactions. By maintaining separate channels, the convolutional filters learn which combinations of these field components are physically meaningful at each angular position. This is conceptually analogous to processing RGB color channels in image CNNs, 
where kernels slide over spatial dimensions while learning cross-channel correlations.

\subsection{Circular Convolution} \label{meth_cnn}
Given the periodic conditions inherent in the angular far-field measurements, we define a circular convolution operation that preserves the $2\pi$-periodicity of the input signal. Throughout the network, both the angular dimension and channel dimension evolve 
through convolutional layers; we denote the angular dimension at layer $\ell$ as $T_\ell$ 
and the number of channels as $C_\ell$. 

\subsubsection{Circular Padding Operator}

To perform convolution along the angular dimension, we apply a circular padding 
strategy to the input of each convolutional layer. For layer $\ell=1,\dots, L$ 
with kernel size $K^{(\ell)}$ applied along the angular dimension, we define 
the layer-specific padding operator
\[
\mathcal{P}_{K^{(\ell)}}:\mathbb{R}^{T_{\ell-1}\times C_{\ell-1}}
\;\to\;\mathbb{R}^{(T_{\ell-1}+K^{(\ell)}-1)\times C_{\ell-1}},
\]
which extends the angular dimension by periodic wrap‐around with
\[
P_{\mathrm{left}}^{(\ell)}=\bigl\lfloor (K^{(\ell)}-1)/2 \bigr\rfloor,
\qquad
P_{\mathrm{right}}^{(\ell)}=(K^{(\ell)}-1)-P_{\mathrm{left}}^{(\ell)}.
\] 
Then the padded representation $\tilde{\mathbf{X}}^{(\ell)}=\mathcal{P}_{K^{(\ell)}}(\mathbf{X}^{(\ell-1)})$, where $\mathbf{X}^{(\ell-1)}$ denotes the input to layer $\ell$, satisfies
\[
\tilde{X}_{i,c}^{(\ell)}
= X_{(i \bmod T_{\ell-1}),\,c}^{(\ell-1)}
\quad
\text{for }i=-P_{\mathrm{left}}^{(\ell)},\dots,T_{\ell-1}-1+P_{\mathrm{right}}^{(\ell)}.
\]
The padded dimension is $T_{\ell-1}+K^{(\ell)}-1$ positions with $C_{\ell-1}$ channels 
preserved. This ensures that all channels are maintained through the padding step.

For example, for layer 1 with $T_0 = 128$ angular positions, $C_0 = 4$ input 
channels, and kernel size $K^{(1)} = 5$, we have $P_{\mathrm{left}}^{(1)} = \lfloor 
(5-1)/2 \rfloor = 2$ and $P_{\mathrm{right}}^{(1)} = (5-1) - 2 = 2$. The padded 
sequence becomes
\[
\tilde{\mathbf{X}}^{(1)}
= \{\;\mathbf{X}^{(0)}_{126},\,\mathbf{X}^{(0)}_{127},\,\mathbf{X}^{(0)}_0,\,\ldots,\,\mathbf{X}^{(0)}_{127},\,\mathbf{X}^{(0)}_0,\,\mathbf{X}^{(0)}_1\}.
\]
This ensures that convolution at the boundary positions (i.e., 
$\mathbf{X}^{(0)}_0$ and $\mathbf{X}^{(0)}_{127}$) can access neighboring angular 
measurements through the periodic wrapping.

This approach guarantees that a convolutional filter sliding along the angular dimension of $\tilde{\mathbf{X}}^{(\ell)}$ operates on a continuous periodic signal, thereby preserving the $2\pi$-periodicity and avoiding boundary artifacts (unlike zero-padding or constant-value padding, which introduces artificial edge values at the boundaries).

\subsubsection{Convolution Operator}  

The circular convolution operation with $N_f^{(\ell)}$ filters, kernel size $K^{(\ell)}$, 
and stride $S^{(\ell)}$ can be viewed as applying a sliding window of length $K^{(\ell)}$ 
(with wrap‐around) along the angular dimension.

For each layer $\ell$, we define
\[
\mathcal{C}_{\mathbf{W},\mathbf{b}}^{(\ell)}:
\mathbb{R}^{(T_{\ell-1}+K^{(\ell)}-1)\times C_{\ell-1}}
\longrightarrow
\mathbb{R}^{T_\ell\times C_\ell}
\]

by
\begin{eqnarray*}
&\bigl[\mathcal{C}_{\mathbf{W},\mathbf{b}}^{(\ell)}(\tilde{\mathbf{X}})\bigr]_{i,j}
=
\sum_{k=0}^{K^{(\ell)}-1}\sum_{c=0}^{C_{\ell-1}-1}
W_{j,k,c}^{(\ell)}\,\tilde{X}_{\,iS^{(\ell)}+k,\,c}
\;+\;b_j^{(\ell)},\\
&\text{for}\quad j=0,1,\dots,N_f^{(\ell)}-1 \  \textrm{(filter)} \quad \text{and} \quad i=0,1,\dots,\left\lceil\frac{T_{\ell-1}}{S^{(\ell)}}\right\rceil-1 \ \textrm{(angular position).}
\end{eqnarray*}
Here $\mathbf{W}^{(\ell)}\in\mathbb{R}^{N_f^{(\ell)}\times K^{(\ell)}\times C_{\ell-1}}$ 
are the learnable filter weights and $\mathbf{b}^{(\ell)}\in\mathbb{R}^{N_f^{(\ell)}}$ 
are the bias terms. Each of the $N_f^{(\ell)}$ filters independently processes the 
padded input and produces one output channel, so the number of output channels is 
$C_\ell = N_f^{(\ell)}$. The output angular dimension is 
$T_\ell = \left\lceil T_{\ell-1}/S^{(\ell)} \right\rceil$. This hierarchical 
mechanism, where successive layers combine the previous layer's channels into new 
learned features, enables the network to capture multi-scale angular patterns.

The full circular‐convolution‐activation operator for layer $\ell$ is
\[
\mathrm{CircularConv1D}^{(\ell)}(\mathbf{X}^{(\ell-1)};\mathbf{W}^{(\ell)},\mathbf{b}^{(\ell)},\phi)
=
\phi\!\bigl(\mathcal{C}_{\mathbf{W},\mathbf{b}}^{(\ell)}\bigl(\mathcal{P}_{K^{(\ell)}}(\mathbf{X}^{(\ell-1)})\bigr)\bigr),
\]
for $\phi(\cdot)$ being the activation function. The output lies in $\mathbb{R}^{T_\ell\times C_\ell}$.

\begin{figure}[t!]
\centering
\includegraphics[width=0.58\textwidth]{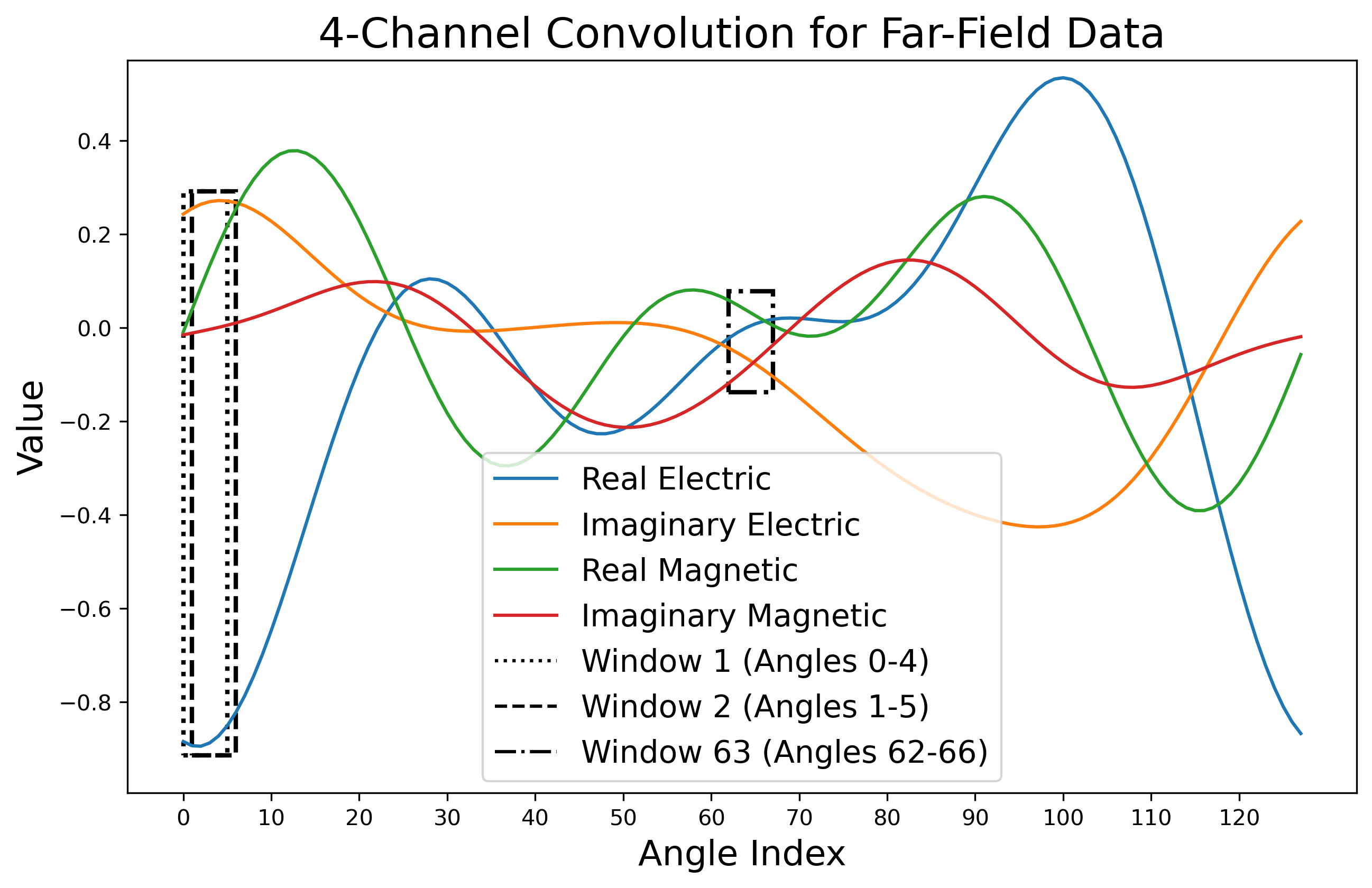}
\caption{Illustration of a sliding window of length $K^{(1)}=5$
moving along the $T_0=128$ angular measurement positions of the $C_0=4$ input channels.}
\label{fig5}
\end{figure}

To illustrate how the convolutional mechanism operates on the far-field dataset, we 
present two representative figures. Figure \ref{fig5} depicts a sliding window of 
length $K^{(1)}=5$ moving along the $T_0=128$ angular measurements of a four‐channel 
input (real and imaginary parts of the electric and magnetic far-fields). Each window 
extracts a $5 \times 4$ patch from the input.

\begin{figure}[t!]
\centering
\includegraphics[width=1\textwidth]{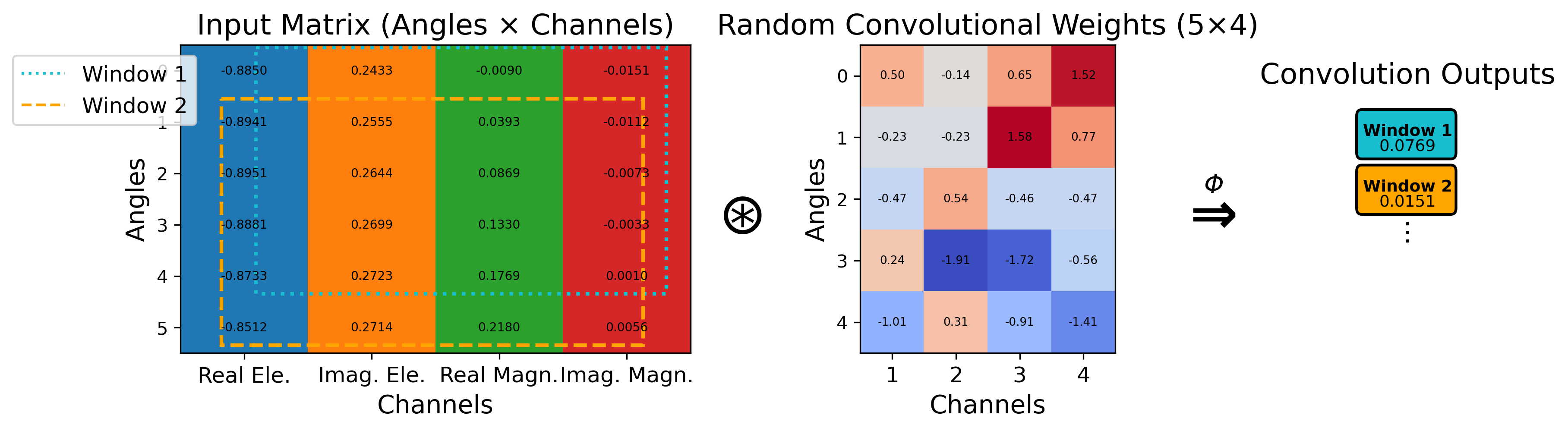}
\caption{Element‐wise convolution for a single filter: a $5\times
4$ kernel multiplies a 5-angle, 4-channel patch.}
\label{fig6}
\end{figure}

Figure \ref{fig6} details the element‐wise convolution: for a kernel size $K^{(1)}=5$, 
the $5\times 4$ weight matrix $W_{j,:,:}^{(1)}$ is element‐wise multiplied with the 
corresponding $5\times 4$ patch, summed over the angle index $k=0,\dots,4$ and channel 
index $c=0,\dots,3$, then added to bias $b_j^{(1)}$ and passed through activation $\phi$. 
The window then shifts by stride $S^{(1)}=1$ to compute the next output, repeating 
until all positions are processed.

\subsubsection{Multiscale Feature Extraction}

Our architecture employs a sequence of circular convolutional layers with progressively 
increasing kernel sizes to capture angular features at multiple scales.

Let $L$ be the number of convolutional layers. For each layer $\ell=1,\dots,L$, 
we define
\[
\mathbf{H}^{(0)} = \mathbf{X} \in \mathbb{R}^{T_0 \times C_0},
\quad
\mathbf{H}^{(\ell)} 
= \mathrm{CircularConv1D}^{(\ell)}\bigl(\mathbf{H}^{(\ell-1)};\mathbf{W}^{(\ell)},\mathbf{b}^{(\ell)},\phi\bigr) 
\in \mathbb{R}^{T_\ell\times C_\ell}.
\]
Each layer $\ell$ is configured with $(N_f^{(\ell)},K^{(\ell)},S^{(\ell)})$. An example 
configuration with 4 layers could be
\begin{align*}
&\ell=1:\;N_f^{(1)}=128 \text{ filters},\;K^{(1)}=5,\;S^{(1)}=1,\quad (\text{local, short-range angular features})\\
&\ell=2: \;N_f^{(2)}=128 \text{ filters},\;K^{(2)}=5,\;S^{(2)}=2,\quad (\text{downsampling by 2})\\
&\ell=3: \;N_f^{(3)}=128 \text{ filters},\;K^{(3)}=15,\;S^{(3)}=1, \quad(\text{medium-scale angular patterns})\\
&\ell=4: \;N_f^{(4)}=128 \text{ filters},\;K^{(4)}=31,\;S^{(4)}=1,\quad (\text{global, long-range angular dependencies}).
\end{align*}

We employ the Swish activation:
\[
\phi(z) = z\,\sigma(z) = \frac{z}{1 + e^{-z}},
\]
which offers smooth, non-monotonic behavior that can improve gradient flow compared to 
ReLU \cite{swish}.

\subsubsection{Channel Attention for Angular Data} \label{attention}

Attention mechanisms have proven effective in deep learning by allowing models to focus on relevant features while suppressing less important information \cite{vaswani2017attention}. For angular electromagnetic data, we employ a lightweight attention module inspired by the Squeeze-and-Excitation (SE) block \cite{hu2018squeeze}, which adaptively recalibrates channel-wise feature responses. For a related application of attention-based CNN architectures in inverse scattering problems, we refer the reader to \cite{ChengGuo}.

Let $T_L$ denote the angular sequence length after the final convolutional layer, 
and let $C_L = N_f^{(L)}$ denote the number of channels. Define the input tensor 
to attention as
\[
\mathbf{H}^{(\mathrm{in})}\in\mathbb{R}^{T_L\times C_L}.
\]

The attention operates in two stages:

Spatial Mixing: First, we apply a small circular convolution 
to mix information from neighboring angular positions:
\[
\mathbf{H}^{(\text{mix})}
= \phi\bigl(\mathcal{C}_{\mathbf{W}_{\text{mix}},\mathbf{b}_{\text{mix}}}
(\mathcal{P}_{K_{\text{mix}}}(\mathbf{H}^{(\mathrm{in})}))\bigr),
\quad
\mathbf{H}^{(\text{mix})}\in\mathbb{R}^{T_L\times C_L},
\]
where $K_{\text{mix}}$ is a small kernel size (e.g.\ 3) and the stride is 1 
(preserving the angular dimension). This is followed by layer normalization:

\[
\mathbf{H}^{(\text{norm})} = \mathrm{LayerNorm}(\mathbf{H}^{(\text{mix})}),
\quad
\mathbf{H}^{(\text{norm})}\in\mathbb{R}^{T_L\times C_L}.
\]

Channel Attention: We then compute channel-wise attention 
weights using global average pooling along the angular dimension:
\[
\mathbf{g}
= \frac{1}{T_L}\sum_{i=0}^{T_L-1}\mathbf{H}^{(\text{norm})}_{i,:}
\in\mathbb{R}^{C_L}.
\]

The attention weights are computed through a two-layer MLP with reduction factor $r$
\[
\mathbf{a}
= \sigma\bigl(\mathbf{W}_2\,\phi(\mathbf{W}_1\,\mathbf{g}+\mathbf{b}_1)+\mathbf{b}_2\bigr)
\in\mathbb{R}^{C_L},
\]
where $\mathbf{W}_1 \in \mathbb{R}^{(C_L/r) \times C_L}$, 
$\mathbf{W}_2 \in \mathbb{R}^{C_L \times (C_L/r)}$, and $\sigma(\cdot)$ is the sigmoid function.
The final attended features are obtained by channel-wise multiplication:
\[
\mathbf{H}^{(\mathrm{att})}
= \mathbf{H}^{(\text{norm})}\odot(\mathbf{1}_{T_L}\otimes\mathbf{a}),
\quad
\mathbf{H}^{(\mathrm{att})}\in\mathbb{R}^{T_L\times C_L},
\]
where \(\odot\) denotes element‐wise multiplication, \(\otimes\) represents 
outer product broadcasting, and \(\mathbf{1}_{T_L}\) is a vector of ones 
with length $T_L$. The outer product \(\mathbf{1}_{T_L}\otimes\mathbf{a}\) 
broadcasts the channel attention weights $\mathbf{a} \in \mathbb{R}^{C_L}$ 
to all $T_L$ angular positions, so that each channel is equally weighted 
across all angular positions.

\subsubsection{Bottleneck Convolution}  

To reduce the number of output channels before the MLP head, we apply a pointwise 
$(1 \times 1)$ convolution to the attended (or final convolutional) features:
\[
\mathbf{H}^{(\mathrm{b})}
= \phi\bigl(\mathbf{H}^{(\mathrm{in})}\,\mathbf{W}^{(\mathrm{b})}+\mathbf{b}^{(\mathrm{b})}\bigr),
\quad
\mathbf{H}^{(\mathrm{b})}\in\mathbb{R}^{T_L\times N_{\mathrm{b}}},
\]
where \(\mathbf{H}^{(\mathrm{in})} \in \mathbb{R}^{T_L \times C_L}\) is either 
\(\mathbf{H}^{(\mathrm{att})}\) (if attention is used) or \(\mathbf{H}^{(L)}\) (otherwise). 
The learnable weights are \(\mathbf{W}^{(\mathrm{b})}\in\mathbb{R}^{C_L\times N_{\mathrm{b}}}\) 
and bias \(\mathbf{b}^{(\mathrm{b})}\in\mathbb{R}^{N_{\mathrm{b}}}\), where 
$N_{\mathrm{b}} < C_L$ is the reduced number of output channels. 

\subsection{MLP Inversion Head} \label{meth_mlp}

After the bottleneck convolution, the resulting tensor \(\mathbf{H}^{(\mathrm{b})}\in\mathbb{R}^{T_L\times N_{\mathrm{b}}}\) is flattened
\[
\mathrm{vec}:\mathbb{R}^{T_L\times N_{\mathrm{b}}}\;\to\;\mathbb{R}^{T_L N_{\mathrm{b}}},
\quad
\mathbf{v} = \mathrm{vec}\bigl(\mathbf{H}^{(\mathrm{b})}\bigr),
\]
to generate the input for the subsequent dense layers.

Each dense layer \(m=1,\dots,M\) with weights \(\mathbf{W}^{(m)}\in\mathbb{R}^{d_m\times d_{m-1}}\) 
and bias \(\mathbf{b}^{(m)}\in\mathbb{R}^{d_m}\) defines the affine map
\[
\mathcal{A}^{(m)}:\mathbb{R}^{d_{m-1}}\to\mathbb{R}^{d_m},\quad
[\mathcal{A}^{(m)}(\mathbf{h})]_i
=\sum_{j=0}^{d_{m-1}-1}W^{(m)}_{i,j}\,h_j + b^{(m)}_i,
\quad i=0,1,\dots,d_m-1.
\]
where \(d_{m-1}\) is the number of inputs to layer \(m\) and \(d_m\) the number of neurons (nodes) in that layer. Including activation \(\phi\), layer normalization, and dropout \(p^{(m)}\) gives  
\[
\mathcal{M}^{(m)}(\mathbf{h})
=\mathrm{Dropout}\bigl(\mathrm{LayerNorm}(\phi(\mathcal{A}^{(m)}(\mathbf{h}))),\,p^{(m)}\bigr).
\]
With \(\mathbf{h}^{(0)}=\mathbf{v}\in\mathbb{R}^{d_0}\) where \(d_0=T_LN_{\mathrm{b}}\), we compute sequentially
\[
\mathbf{h}^{(m)} = \mathcal{M}^{(m)}\bigl(\mathbf{h}^{(m-1)}\bigr),
\quad m=1,\dots,M.
\]
The MLP operator is then given by applying the final affine map to the last hidden state:
\[
\mathcal{A}^{(\mathrm{out})}\bigl(\mathbf{h}^{(M)}\bigr)
= \mathbf{W}^{(\mathrm{out})}\,\mathbf{h}^{(M)} + \mathbf{b}^{(\mathrm{out})},
\] 
mapping \(\mathbb{R}^{d_0}\to\mathbb{R}^{d_{\mathrm{out}}}\).
The MLP head thus consists of these \(M\) fully connected layers that map the flattened features to either class probabilities or regression targets. 

To prevent overfitting in the MLP head, $L2$ regularization can be applied to the weights of the dense layers
\[
\mathcal{R}(\boldsymbol{\theta})
=
\lambda_{\mathrm{reg}}
\sum_{m=1}^{M}
\lVert \mathbf{W}^{(m)}\rVert_F^2,
\]
where \(\lambda_{\mathrm{reg}}\) is the regularization parameter and \(\|\cdot\|_F\) denotes the Frobenius norm. Recall that $\boldsymbol{\theta}$ denotes the collection of all trainable parameters of the network.

The final output layer produces the class probabilities 
\[
\widehat{\mathbf{p}}^\mathrm{class}
=
\mathrm{softmax}\bigl(\mathbf{W}^{(\mathrm{out})}\,\mathbf{h}^{(M)} + \mathbf{b}^{(\mathrm{out})}\bigr)\ \in\Delta^{K_{\mathrm{cls}}-1},
\]
for the classification problem or 
\[
\widehat{\mathbf{s}}^\mathrm{reg}
=
\mathbf{W}^{(\mathrm{out})}\,\mathbf{h}^{(M)} + \mathbf{b}^{(\mathrm{out})}\ \in\mathbb{R}^P,
\]
composed of the coefficients representing the boundary curve and the impedance function for the regression problem. Here,
\[
\Delta^{K_{\mathrm{cls}}-1}
= \Bigl\{\mathbf{p}\in[0,1]^{K_{\mathrm{cls}}} : \sum_{k=1}^{K_{\mathrm{cls}}}p_k = 1\Bigr\}\subset\mathbb{R}^{K_{\mathrm{cls}}},
\]
denotes the probability simplex of dimension \(K_{\mathrm{cls}}-1\), i.e.\ the set of all valid class‐probability vectors.

\subsection{Training as Inverse‐Operator Approximation} \label{meth_train}

Having described all components of the network, we now define the complete model as a sequence of transformations mapping the input measurements $\mathbf{X}$ to the outputs: class probabilities $\widehat{\mathbf{p}}^{\mathrm{class}}$ or regression estimates $\widehat{\mathbf{s}}^{\mathrm{reg}}$:

\begin{enumerate}
  \item Multiscale Circular Convolution\\
    \[
      \mathbf{H}^{(0)} = \mathbf{X},\quad
      \mathbf{H}^{(\ell)}
      = \mathrm{CircularConv1D}\bigl(\mathbf{H}^{(\ell-1)};\mathbf{W}^{(\ell)},\mathbf{b}^{(\ell)},\phi\bigr),
      \;\ell=1,\dots,L.
    \]

  \item Angular Attention (optional)\\
    \[
      \mathbf{H}^{(\mathrm{att})}
      = \mathrm{Attn}\bigl(\mathbf{H}^{(L)}\bigr).
    \]

  \item Bottleneck Convolution\\
    \[
      \mathbf{H}^{(\mathrm{b})}
      = \phi\bigl(\mathbf{H}^{(\mathrm{in})}\,\mathbf{W}^{(\mathrm{b})}
        + \mathbf{b}^{(\mathrm{b})}\bigr),
      \quad
      \mathbf{H}^{(\mathrm{in})} = 
        \begin{cases}
          \mathbf{H}^{(\mathrm{att})}, & \text{if attention used},\\
          \mathbf{H}^{(L)},           & \text{otherwise}.
        \end{cases}
    \]

  \item Flatten to Vector
  \[
      \mathbf{v}
      = \mathrm{vec}\bigl(\mathbf{H}^{(\mathrm{b})}\bigr).
\]
\item MLP Head\\
\[
\mathbf{h}^{(0)} = \mathbf{v},\quad
\mathbf{h}^{(m)} = \mathcal{M}^{(m)}\bigl(\mathbf{h}^{(m-1)}\bigr),
\quad m=1,\dots,M.
\]
The final affine operator \(\mathcal{A}^{(\mathrm{out})}\) yields
\[
\widehat{\mathbf{p}}^\mathrm{class}
= \mathrm{softmax}\!\bigl(\mathcal{A}^{(\mathrm{out})}(\mathbf{h}^{(M)})\bigr),
\quad \textrm{or} \quad \widehat{\mathbf{s}}^\mathrm{reg}
= \mathcal{A}^{(\mathrm{out})}\bigl(\mathbf{h}^{(M)}\bigr)
= \mathbf{W}^{(\mathrm{out})}\,\mathbf{h}^{(M)} + \mathbf{b}^{(\mathrm{out})}.
\]
\end{enumerate}
Putting it all together, the end‐to‐end operator is given by
\[
\mathcal{N}_{\boldsymbol{\theta}}(\mathbf{X})
= \mathcal{A}^{(\mathrm{out})}\Bigl[
  \mathcal{M}^{(M)}\Bigl(\cdots\mathcal{M}^{(1)}
    \Bigl(\mathrm{vec}\bigl(
      \mathrm{Bottl}\bigl(
        \mathrm{Attn}\bigl(
          \mathrm{CircConv}^{(L)}(\cdots\mathrm{CircConv}^{(1)}(\mathbf{X}))
        \bigr)
      \bigr)
    \bigr)\Bigr)
  \Bigr)
\Bigr],
\]
where each component operator is defined as above. The overall mapping can be decomposed into

\begin{equation*}
\mathcal{X}=\mathbb{R}^{T_0\times C_0}
\;\xrightarrow{\mathrm{CircConv1D}}\;
\mathbb{R}^{T_L\times C_L}
\;\xrightarrow{\mathrm{Attn}}\;
\mathbb{R}^{T_L\times C_L}\\
\;\xrightarrow{\mathrm{Bottl}}\;
\mathbb{R}^{T_L\times N_{\mathrm{b}}}
\;\xrightarrow{\mathrm{vec}}\;
\mathbb{R}^{d_0} 
\;\xrightarrow{\mathrm{MLP}}\;
 \mathbb{R}^{d_{\mathrm{out}}}.
\end{equation*}
The network parameters (weights and biases) 
\[
\begin{aligned}
\boldsymbol{\theta} ={}& \{\mathbf{W}^{(\ell)},\mathbf{b}^{(\ell)}\}_{\ell=1}^{L}
\;\cup\;
\{\mathbf{W}_{\mathrm{mix}},\mathbf{b}_{\mathrm{mix}}\}
\;\cup\;
\{\mathbf{W}_1,\mathbf{b}_1,\mathbf{W}_2,\mathbf{b}_2\}\\
&\cup\;
\{\mathbf{W}^{(\mathrm{b})},\mathbf{b}^{(\mathrm{b})}\}
\;\cup\;
\{\mathbf{W}^{(m)},\mathbf{b}^{(m)}\}_{m=1}^{M}
\;\cup\;
\{\mathbf{W}^{(\mathrm{out})},\mathbf{b}^{(\mathrm{out})}\}
\end{aligned}
\]
are optimized by minimizing the appropriate loss function. The operator $\mathcal{N}_{\boldsymbol{\theta}}$ is trained to minimize either the Categorical Cross-Entropy or the Mean Squared Error (MSE):
\[
\mathcal{L}_{\mathrm{class}}(\boldsymbol{\theta})
=
-\frac{1}{N}\sum_{n=1}^N\sum_{k=1}^{K_{\mathrm{cls}}}\mathrm{y}_{n,k}\log (\widehat{p}_{n,k}), \quad  \quad 
\mathcal{L}_{\mathrm{reg}}(\boldsymbol{\theta})
=
\frac{1}{N}\sum_{n=1}^N\bigl\|\widehat{\mathbf{s}}_n-\mathbf{s}_n\bigr\|_2^2,
\] 
for the classification and regression problem, respectively. Here $N$ is the number of samples, $\widehat{p}_{n,k}$ is 
the predicted probability of the class $k$ in sample $n$, and $\mathrm{y}_{n,k}$ is the  corresponding one-hot encoded label. The parameters
$\widehat{\mathbf{s}}_n$ and $\mathbf{s}_n$ 
denote the vector of the predicted and true values of the sample $n,$ respectively. This is performed using Adam optimizer with learning rate $\alpha$ and gradient clipping threshold $\gamma$
\[
\boldsymbol{\theta}_{t+1}
=
\boldsymbol{\theta}_t
-\alpha\,\mathrm{clip}\bigl(\nabla_{\boldsymbol{\theta}}\mathcal{L},\gamma\bigr).
\]
 
This formulation frames training as an approximation of the inverse scattering operator $\mathcal{F}^{-1}$ by fitting the mapping $\mathcal{N}_{\boldsymbol{\theta}}$ to the training data.

\subsection{On approximation and convergence}

The proposed network is trained on large synthetic datasets obtained by solving the forward scattering problem and therefore realizes an empirical approximation to the inverse operator $\mathcal{F}^{-1}$ by fitting the input–output mapping. Intuitively, this approach is based on the representational capacity of the network, sufficient training data, and successful loss minimization.

Classical universal approximation results guarantee that sufficiently large feedforward networks can approximate continuous mappings on compact sets \cite{cybenko1989approximation}. Recent work on operator learning shows how neural architectures can be used to approximate maps between function spaces from input–output pairs \cite{kovachki2021neural, lu2019deeponet}.  For inverse problems, theoretical advances establish that neural networks and operator learning architectures can approximate inverse mappings with stability guarantees under suitable regularization and sufficient training data \cite{buskulic2024convergence, nelsen2024operator}. Specifically, recent theoretical frameworks combine neural operators with 
classical regularization theory, establishing stability and error control for nonlinear 
inverse problems \cite{ scherzer2025regularization, scherzer2025neural}.
 In practice, ill-posedness or non-injectivity implies that the fitted network recovers a regularized estimator and therefore architectural inductive bias and explicit regularization are important for stability.
 
A rigorous operator-theoretic convergence analysis specific to our circular-CNN architecture is beyond the present scope; here we validate the approximation empirically through learning curves, sample-size scaling experiments, noise-stability tests, and architectural experiments (see Section \ref{inverse} and Appendix \ref{A}).

\section{Numerical implementation of the inverse problem} \label{inverse}

In this section, we solve in two steps the inverse problems: we first classify the unknown obstacles and then we reconstruct their boundary curves and the impedance function.

All models are developed using the Keras library \cite{keras}, built on the TensorFlow open-source platform \cite{tensorflow}. We first preprocess the data by scaling using StandardScaler to ensure each feature has zero mean and unit standard deviation. Data preprocessing and evaluation metrics are implemented using scikit-learn built-in functions (version 0.23.2). Each dataset is randomly split with 80\% used for training, 10\% for validation, and 10\% for testing.

All scripts are written in Python 3.8. Initial computations and model prototyping are performed on a standard Intel Core i7 (CPU) computer with 8GB RAM. For the computational demands of the deeper models, we utilize an HPC system with GPU acceleration. The HPC configuration employs a single NVIDIA Tesla GPU (based on the CUDA 10.1.168 module) on one node with 56GB allocated memory per job, running TensorFlow 2.4.1.

 The architectures and hyperparameters (e.g., number of circular‐convolution layers and filters, kernel sizes, bottleneck dimension, dense layer depth and node sizes, learning rate, $L2$ regularization, and dropout rates) were selected and tuned via extensive trial-and-error to balance model capacity and task complexity.

For the reader’s convenience, Table \ref{table_all} provides a summary of the dataset structure and size employed across the different inverse problems, along with the corresponding number of output variables. 

\begin{table}
\caption{\label{table_all}Input-output characteristics.}
\centering
\begin{tabular}{@{}lllll|l}
\br
Inverse Problem               & $C_0$   & $T_0$ & $\phi$ & $N$ & Output   \\
\mr
Classification           & 2 & 32 & 0 & 90000 & $\{1,2,3\}$ \\
Regression (Peanut)      &  2 & 32 & 0 & 30000 & $(\alpha, \beta, x_0, y_0, \lambda)$ \\
Regression (Kite)        &  2 & 32 & 0 & 30000 & $(\alpha, \beta, \gamma, x_0, y_0, \lambda)$ \\
Regression (Star-fixed $\lambda$) & 4 & 128 & 0 & 80000 & $(\alpha_0, ...,\alpha_5, \beta_1, ..., \beta_5, x_0, y_0 )$ \\
Regression (Star-variable $\lambda$) & 8 & 128 & $\{0, \pi \}$ & 120000 & $(\alpha_0, ...,\alpha_5, \beta_1, ..., \beta_5, x_0, y_0, \lambda )$ \\
\br
\end{tabular}
\end{table}

\subsection{Classification of obstacles} \label{class}

We propose a classification framework to address the first step--obstacle perception--of the inverse problem. Specifically, we train a model to accurately identify the class of an impedance obstacle (e.g., peanut-shaped, kite-shaped, or star-shaped) using far-field measurement data.

We consider only the real and imaginary electric fields, at a single incidence ($\phi = 0$), measured at $T_0=32$ equidistant angle positions on the unit circle. This corresponds to $C_0=2$ field components (channels) in our notation. We solve the direct problem to generate 30000 samples for every class, as discussed in Section \ref{direct}. Following the data representation from Section \ref{meth_input}, each sample $\mathbf{x}_n\in\mathbb{R}^{64}$ is reshaped as $\mathbf{X}_n\in\mathbb{R}^{32\times 2}$, for $n=1,\ldots,90000$. The class labels belong to the set $\{1,2,3\}$ where $1=$ peanut-shaped, $2=$ kite-shaped and $3=$ star-shaped.

\begin{figure}[h!]
    \centering
    \begin{minipage}{0.48\textwidth}
        \centering
        \includegraphics[width=\textwidth]{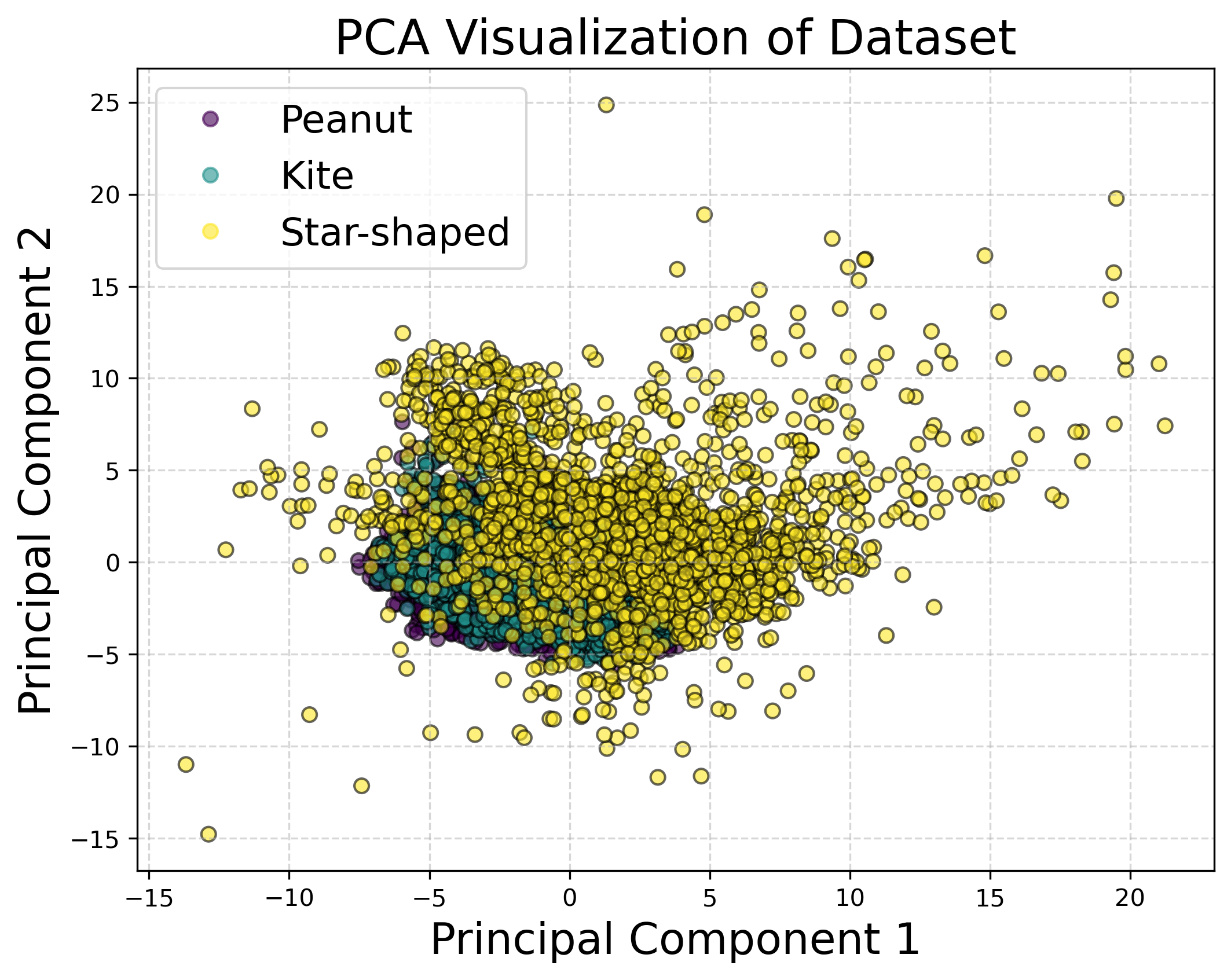} 
    \end{minipage} \hfill
    \begin{minipage}{0.48\textwidth}
        \centering
        \includegraphics[width=\textwidth]{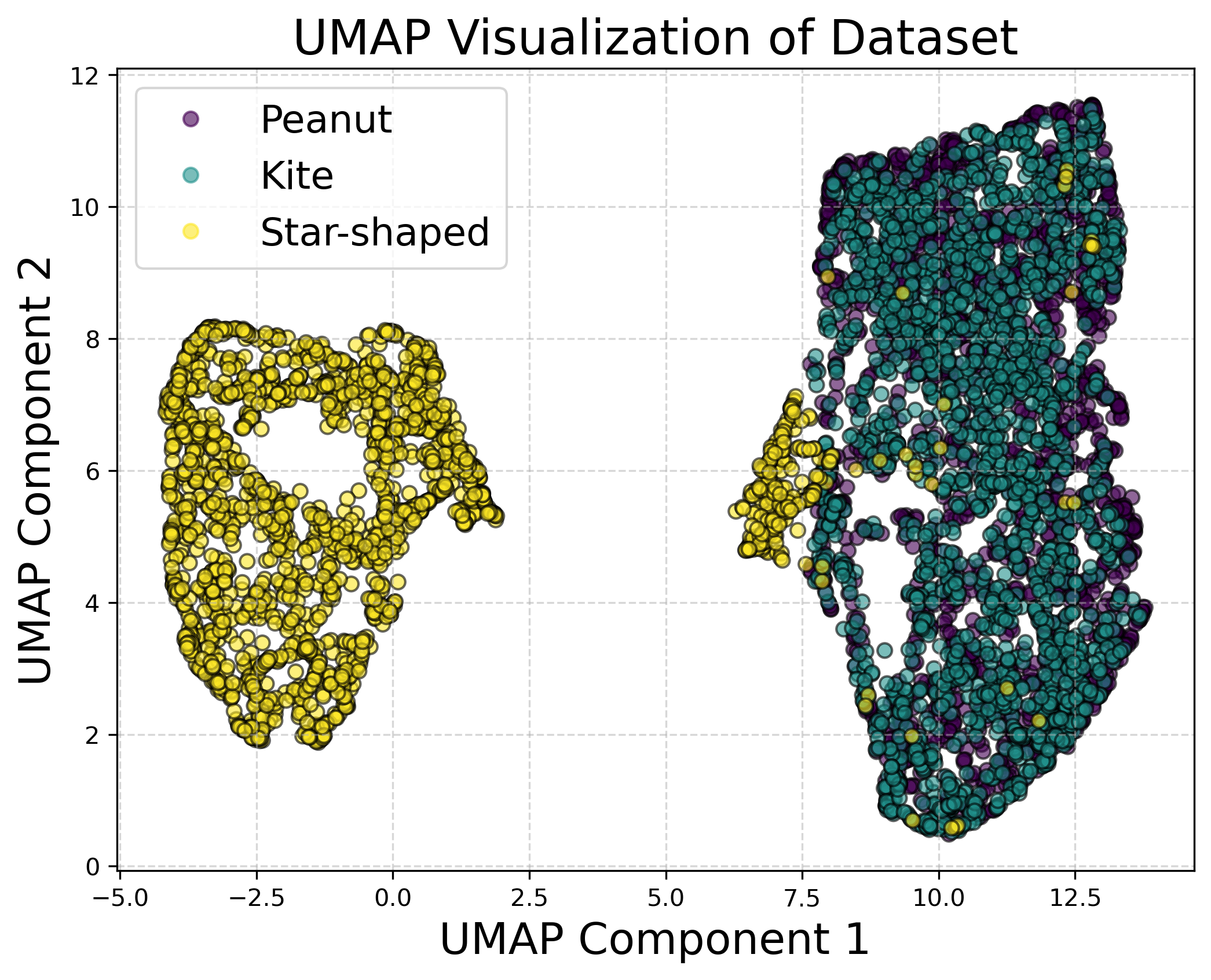} 
    \end{minipage}
    \caption{Comparison of PCA (left) and UMAP (right) for 2D dimension reduction of far-field data.}
    \label{fig7}
\end{figure}

As a pre-processing step, to gain insight into how the obstacles are represented in the feature space, we apply two different unsupervised dimension reduction methods: PCA and UMAP \cite{pca,umap}. The 2D projections for randomly chosen $3\times 2000$ samples of the far-field dataset are shown in Figure \ref{fig7}. In the PCA plot (left), we observe strong overlap between all the classes, while the UMAP plot (right) provides a clearer separation between the star-shaped and the other two classes, which appear to fully overlap. This suggests that distinguishing between obstacle classes may not be straightforward, at least in the 2D projected scattering data.

The performance of the classification model is typically assessed using the accuracy metric
\[
\text{accuracy} = \frac{\text{number of correct predictions}}{\text{total number of predictions}},
\]
and for each class \(k\in K_{\mathrm{cls}}\), the recall is calculated as
\[
\text{recall}_{k} = \frac{\text{true positives for } k}{\text{true positives for } k + \text{false negatives for } k}.
\]

The CNN architecture is given in Table \ref{tabAp1} while the loss and accuracy plots are shown in Figure \ref{figAp1}, of Appendix \ref{A1}. The results are presented in Table \ref{tab1}. Additionally, the confusion matrix corresponding to the test set is shown in Figure \ref{fig8}, where the diagonal displays the recall for each class, with values scaled from 0 to 1. It is demonstrated that the model accurately predicts the class of the unknown obstacle based only on the (real and imaginary parts) electric measurements from one incident wave. We notice that the kite class is more challenging to classify correctly compared to the other two classes. Star-shaped obstacles are classified with perfect accuracy, consistent with the preceding unsupervised pre-processing analysis. 

\noindent
\begin{minipage}{0.42\textwidth}
  \captionsetup{type=table}
  \captionof{table}{\label{tab1}Metrics for the classification problem.}
  \centering
\begin{tabular}{@{}llll}
    \br
    Metric               & Train   & Valid & Test   \\
    \mr
    Accuracy             & 0.992   & 0.988      & 0.988  \\
    Recall (Peanut)      & 0.998   & 0.995      & 0.994  \\
    Recall (Kite)        & 0.978   & 0.969      & 0.972  \\
    Recall (Star-shaped) & 1.000   & 1.000      & 1.000  \\
    \br
    \end{tabular}
\end{minipage}%
\hfill
\begin{minipage}{0.44\textwidth}
  \centering
  \includegraphics[width=\textwidth]{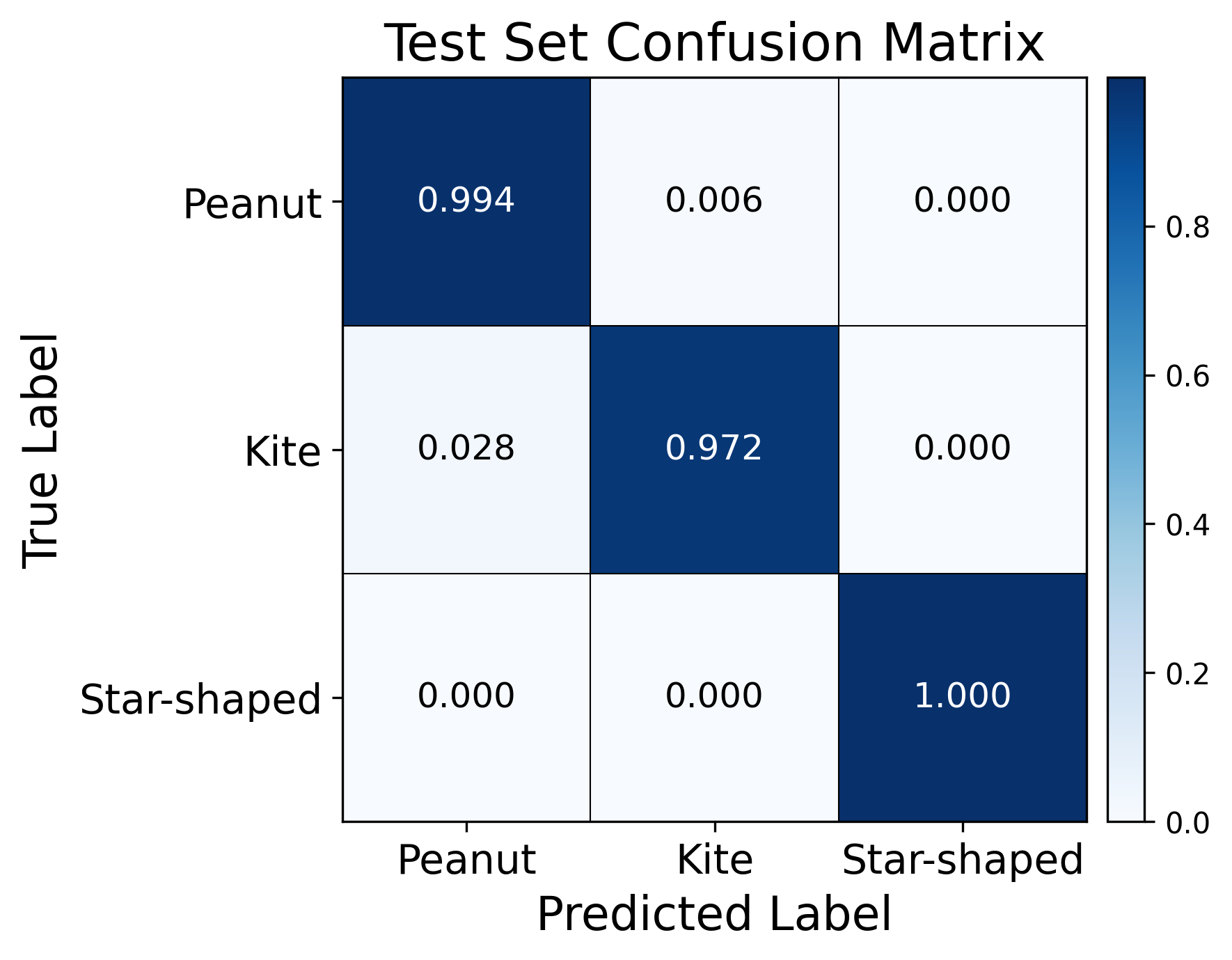}
  \captionof{figure}{The confusion matrix of the test set in the classification problem.}
  \label{fig8}
\end{minipage}
\medskip

Finally, we evaluate the model's robustness by introducing Gaussian noise to the scaled test dataset at varying levels ($0.5\%$ to $5\%$). The results are  included in Table \ref{tab2}. We observe that as noise increases, overall accuracy decreases but remains high, demonstrating strong resilience. Notably, the peanut class recall degrades most rapidly with noise while kite and star-shaped recalls remain stable, reflecting that peanuts become easily misclassified under noisy conditions. 

\begin{table}
\caption{\label{tab2}Classification performance on noisy test data.}
\centering
\begin{tabular}{@{}lllll}
\br
Noise level & Accuracy & Recall (Peanut) & Recall (Kite) & Recall (Star-shaped) \\
\mr
0.5\%  & 0.986 & 0.984 & 0.974 & 1.000 \\
1.0\%  & 0.967 & 0.931 & 0.971 & 1.000 \\
2.0\%  & 0.926 & 0.813 & 0.967 & 1.000 \\
5.0\%  & 0.842 & 0.580 & 0.948 & 1.000 \\
\br
\end{tabular}
\end{table}

Having demonstrated successful obstacle classification, we next turn to the regression-based shape reconstruction problem in the following sections, beginning with the peanut-shaped obstacle.

The performance of the following regression models is assessed using the root mean squared error (RMSE) and the coefficient of determination $R^2$, defined respectively as
\begin{equation*}
\begin{aligned}
  \mathrm{RMSE} &= \sqrt{\frac{1}{N}\sum_{n=1}^N \bigl\| \widehat{\mathbf{s}}_n - \mathbf{s}_n \bigr\|_2^2}, \\
R^2 &= 1 - \frac{\sum_{n=1}^N \bigl\| \widehat{\mathbf{s}}_n - \mathbf{s}_n \bigr\|_2^2}{\sum_{n=1}^N \bigl\| \mathbf{s}_n - \overline{\mathbf{s}} \bigr\|_2^2},
\end{aligned}
\end{equation*}
where $\overline{\mathbf{s}}$ denotes the mean of the exact parameter vectors over all samples.

\subsection{Inverse problem for peanut-shaped obstacles} \label{peanCNN} 

The inverse problem involves reconstructing both the coefficients of the boundary parametrization and the impedance constant~$\lambda$. According to (\ref{pean}), this requires estimating a total of 5 parameters: $\alpha$, $\beta$, the center coordinates $\bm{x}_0 = (x_0, y_0)$, and the variable $\lambda$. 

We use the dataset described in the previous classification section, consisting of 30000 peanut-shaped obstacle samples, where each input $\mathbf{x}_n\in\mathbb{R}^{64}$ is reshaped to $\mathbf{X}_n\in\mathbb{R}^{32\times 2}$ (i.e, using only the real and imaginary electric measurements at $32$ angle positions). The network architecture is detailed in Table~\ref{tabAp2}, and the corresponding loss evolution is shown in Figure~\ref{figAp2} of the Appendix \ref{A2}. 

\begin{figure}[ht!]
    \centering
    \begin{minipage}[b]{0.48\textwidth}
        \centering
        \includegraphics[width=\textwidth]{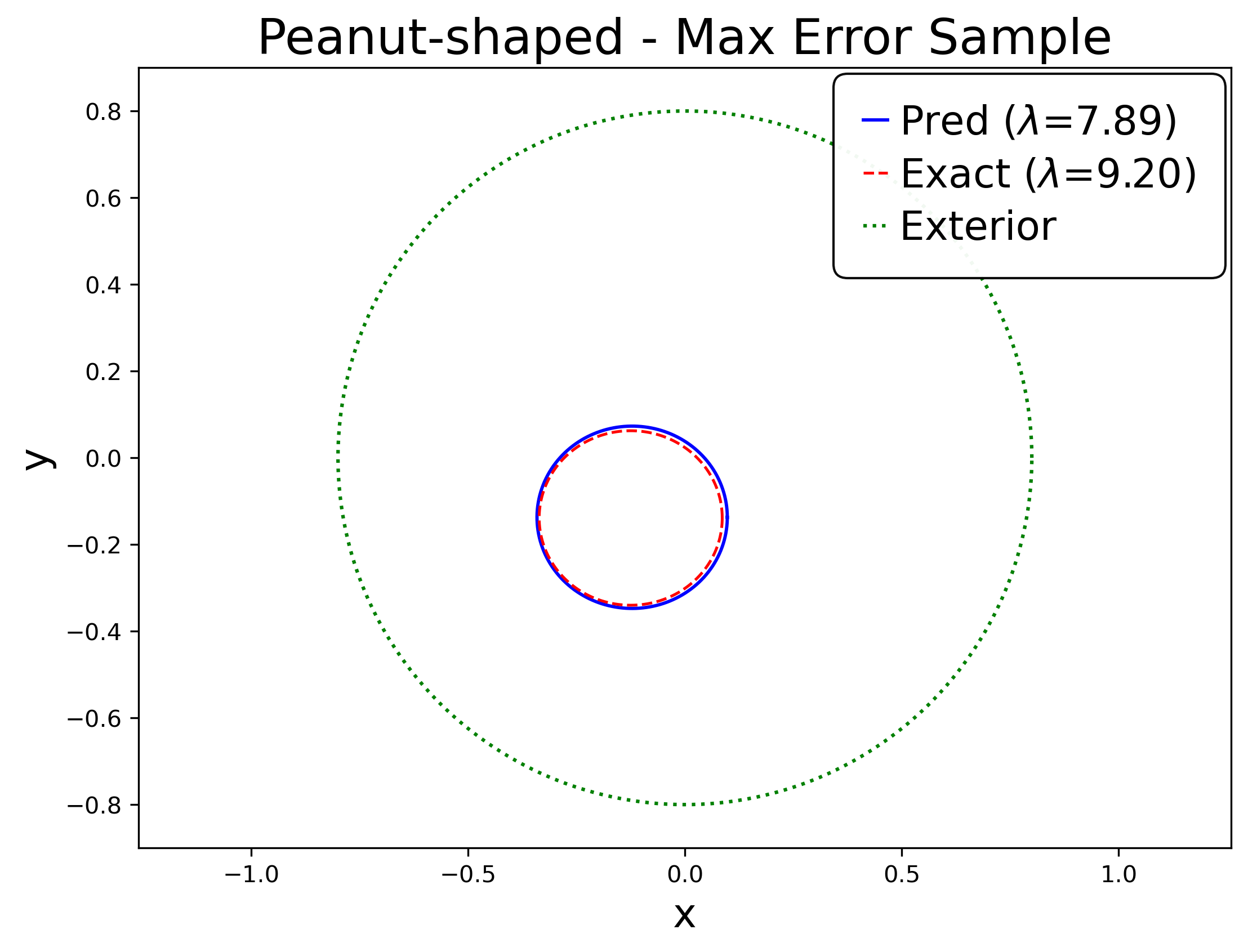} 
        \\[0.5ex]
            \end{minipage}\hfill
    \begin{minipage}[b]{0.48\textwidth}
        \centering
        \includegraphics[width=\textwidth]{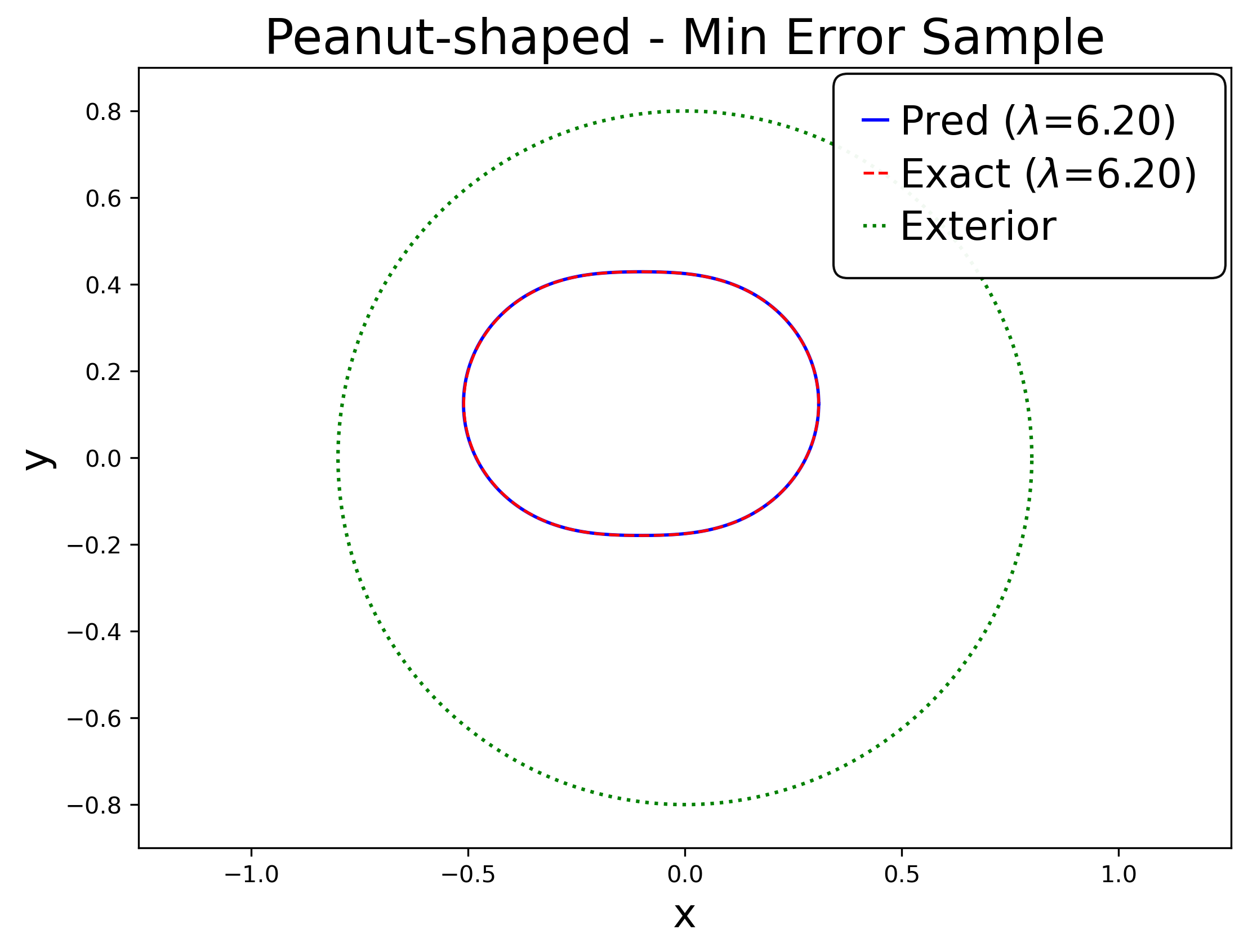} 
        \\[0.5ex]
         \end{minipage}
    \begin{minipage}[b]{0.48\textwidth}
        \centering
        \includegraphics[width=\textwidth]{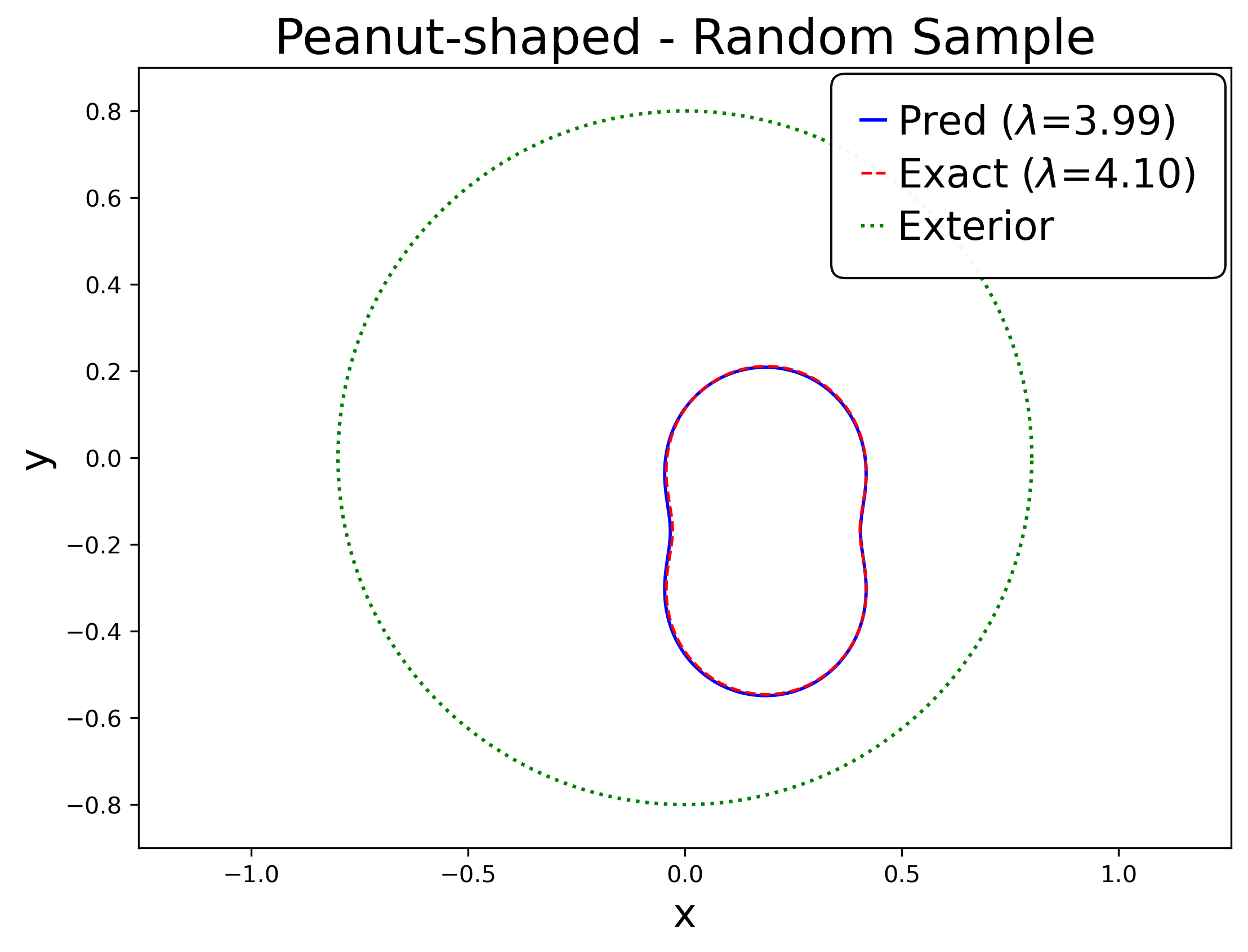} 
        \\[0.5ex]
    \end{minipage}

    \caption{Peanut-shaped reconstructions: max error sample (top left),  min error sample (top right) and randomly chosen (bottom).}
    \label{fig9}
\end{figure}

Table~\ref{tab3} summarizes the overall regression performance. Detailed per-parameter results on the test set appear in Table~\ref{tabAp3}, with the regression plot for $\lambda$ presented in Figure~\ref{figAp3} (Appendix \ref{A2}). Peanut reconstructions for the test set are shown in Figure~\ref{fig9}: the samples with maximum error (top left), minimum error (top right), and a randomly chosen one (bottom). The histogram of test set RMSE values is provided in Figure~\ref{figAp4}. Overall, the regression model demonstrates excellent performance, achieving high $R^2$ score and very low RMSE, indicating highly accurate parameter recovery for the peanut-shaped obstacles. 

\begin{table}
\caption{\label{tab3} Peanut-shaped regression performance.}
\centering
\begin{tabular}{@{}lccc}
\br
Metric             & Train & Valid & Test   \\
\mr
$R^2$ score (\%)   & 99.93    & 99.91      & 99.92  \\
RMSE               & 0.0528   & 0.0615     & 0.0594 \\
\br
\end{tabular}
\end{table}

Finally, we evaluate the model’s robustness to additive measurement noise. Table~\ref{tab4} summarizes the regression performance under varying noise levels, and Figure~\ref{fig10} illustrates the reconstruction for the randomly chosen sample. As demonstrated, the model maintains stable results in the presence of noise. 

\noindent
\begin{minipage}{0.45\textwidth}
  \captionsetup{type=table}
  \captionof{table}{\label{tab4}Peanut-shaped regression performance on noisy test data.}
  \centering
\begin{tabular}{@{}lcc}
      \br
      Noise level & $R^2$ score (\%) & RMSE    \\
      \mr
      0.5\%       & 99.90             & 0.0631  \\
      1.0\%       & 99.86             & 0.0726  \\
      2.0\%       & 99.69             & 0.1013  \\
      5.0\%       & 98.51             & 0.2128  \\
      \br
    \end{tabular}
\end{minipage}%
\hfill
\begin{minipage}{0.48\textwidth}
  \centering
  \includegraphics[width=\textwidth]{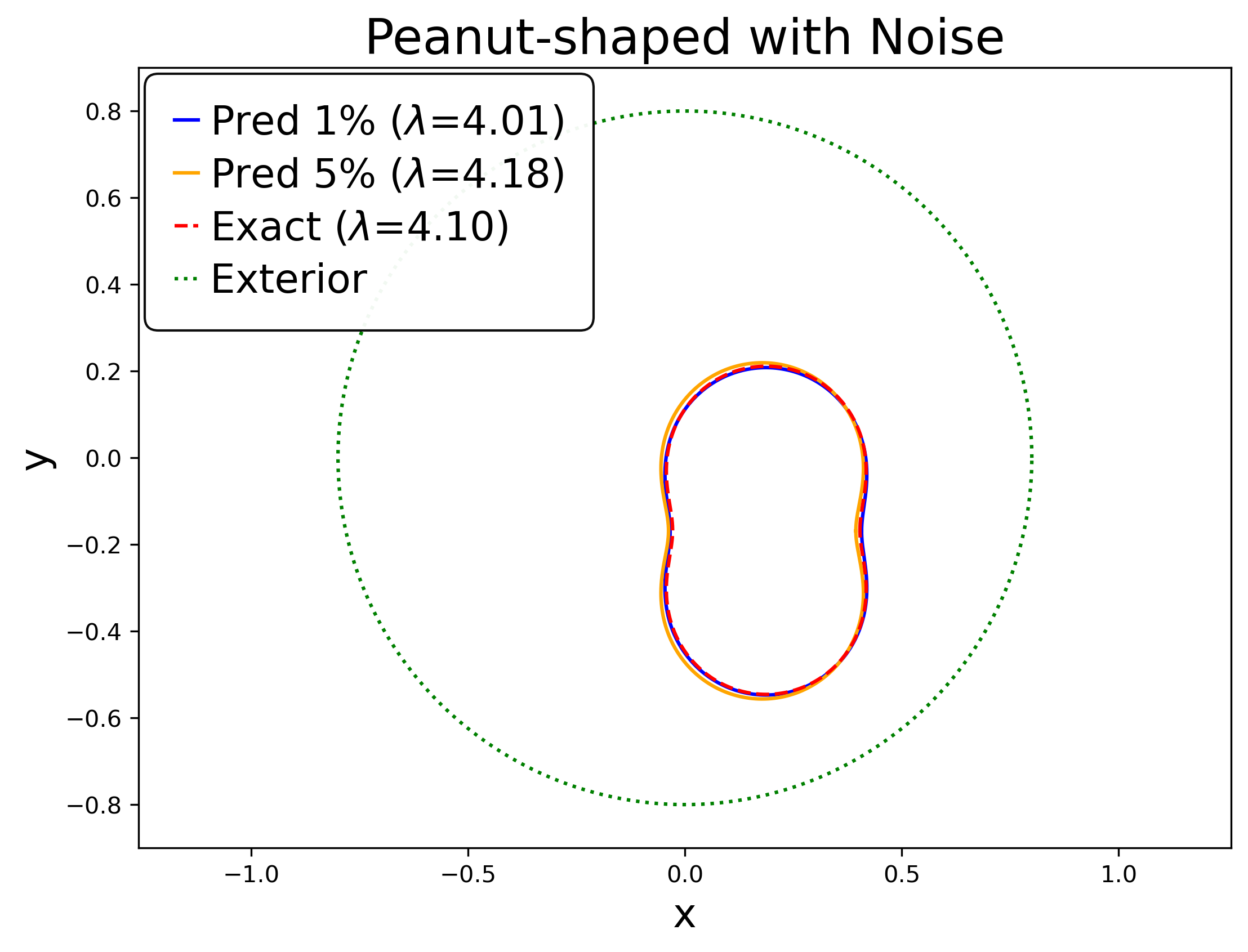}
  \captionof{figure}{Peanut-shaped reconstructions for the randomly chosen sample in presence of increasing levels of noise.}
  \label{fig10}
\end{minipage}

\subsection{Inverse problem for kite-shaped obstacles} \label{kiteCNN} 

Next, we proceed with the inverse problem for the kite-shaped obstacles. This is to recover the boundary coefficients given in (\ref{kite}) and the impedance function~$\lambda$. Thus, we have a regression problem with 6 output parameters: $\alpha$, $\beta$, $\gamma$, the center $\bm{x}_0 = (x_0, y_0)$, and the constant $\lambda$. 

In a manner analogous to the peanut-shaped inverse problem, we employ a dataset of 30000 kite-shaped obstacle samples, each of which is reshaped to $\mathbf{X}_n\in\mathbb{R}^{32\times 2}$. Table \ref{tabAp4} details the CNN architecture, and the loss is presented in Figure \ref{figAp5} (Appendix \ref{A3}).

\begin{figure}[ht!]
    \centering
    \begin{minipage}[b]{0.48\textwidth}
        \centering
        \includegraphics[width=\textwidth]{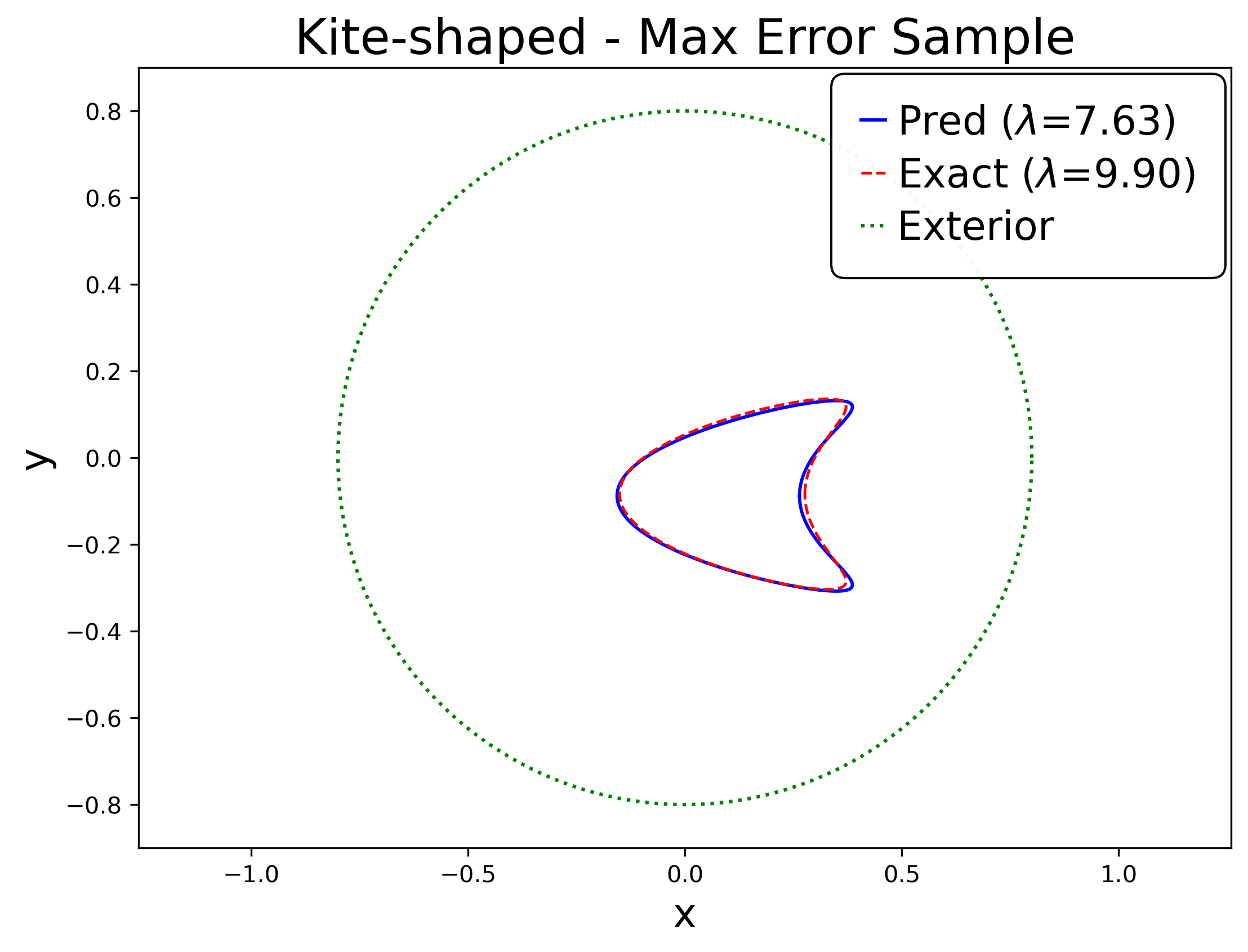} 
        \\[0.5ex]
            \end{minipage}\hfill
    \begin{minipage}[b]{0.48\textwidth}
        \centering
        \includegraphics[width=\textwidth]{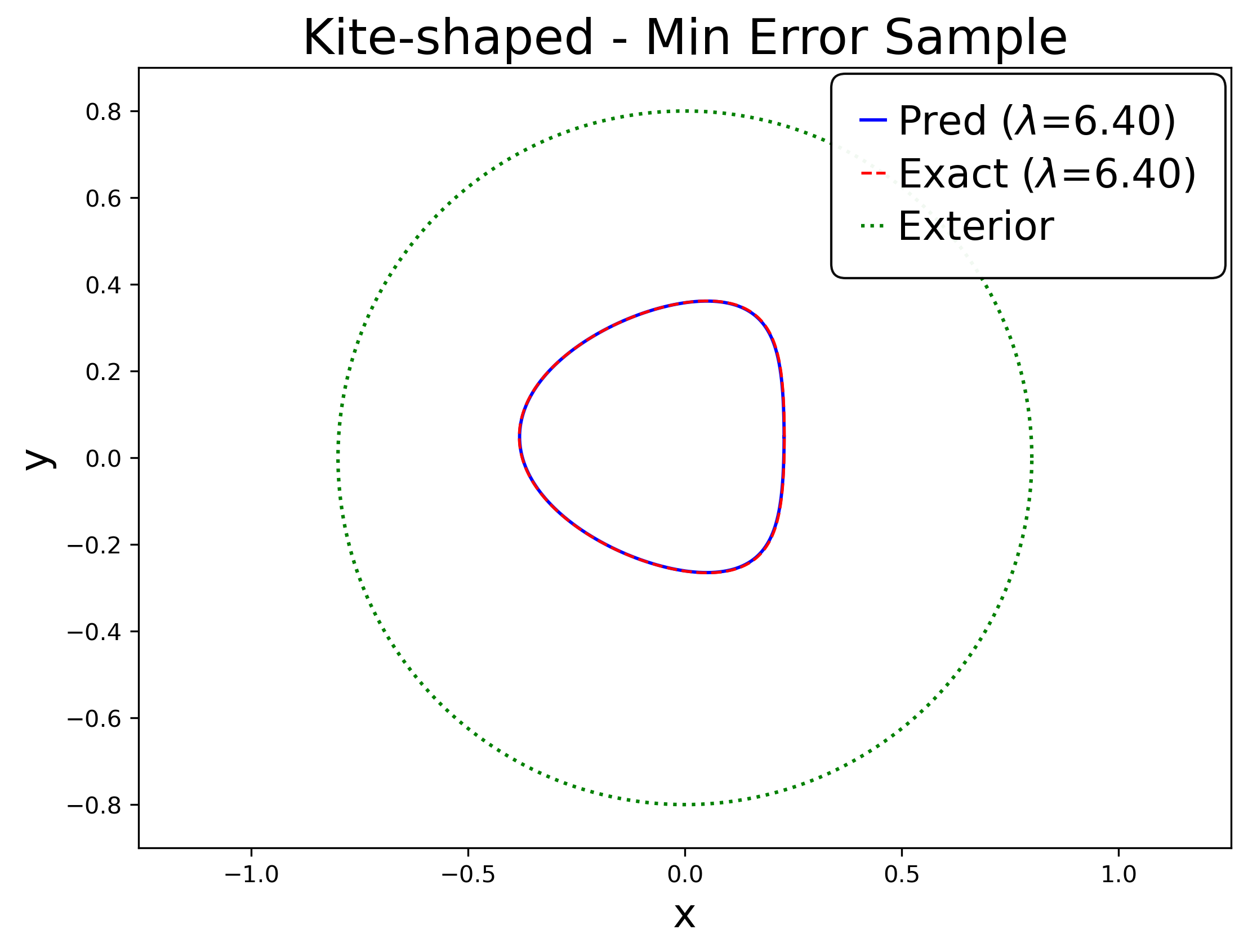} 
        \\[0.5ex]
         \end{minipage}
    \begin{minipage}[b]{0.48\textwidth}
        \centering
        \includegraphics[width=\textwidth]{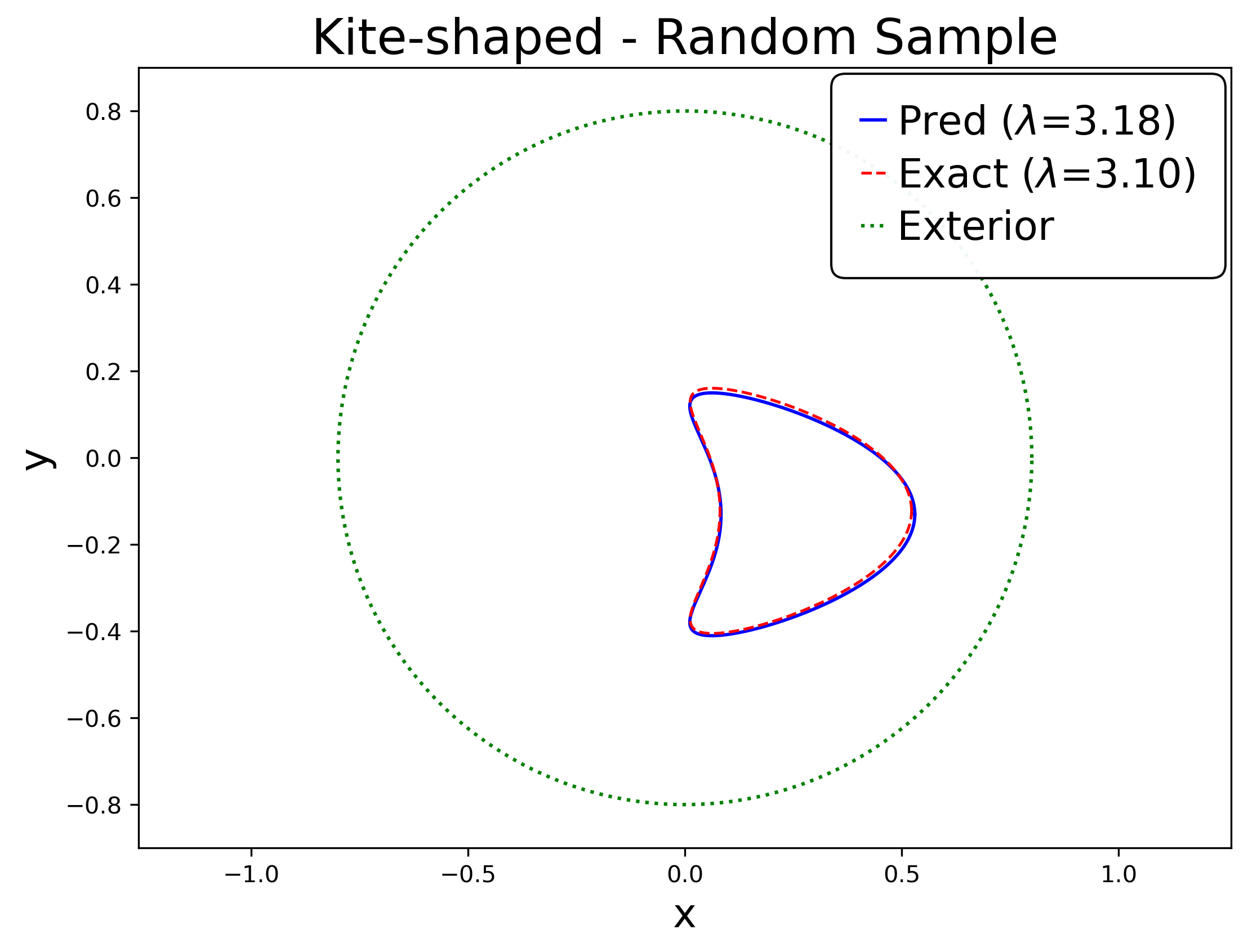} 
        \\[0.5ex]
    \end{minipage}

    \caption{Kite-shaped reconstructions: max error sample (top left),  min error sample (top right) and randomly chosen (bottom).}
    \label{fig11}
\end{figure}

Table \ref{tab5} presents the overall regression metrics for kite-shaped obstacles, and Table \ref{tabAp5} provides a detailed breakdown of test set performance by parameter. The regression of the impedance $\lambda$ values is illustrated in Figure \ref{figAp6}. Example reconstructions for the kite-shaped geometry are presented in Figure \ref{fig11}, showing the cases with the largest and smallest reconstruction errors, along with a representative randomly selected sample. To visualize the predictive accuracy of the model, Figure \ref{figAp7} shows the distribution of RMSE values over all test samples. As demonstrated, the model performs very well on kite-shaped obstacles, with very high $R^2$ scores   and low RMSE, showing it accurately recovers the regression parameters.

\begin{table}
\caption{\label{tab5}Kite-shaped regression performance.}
\centering
\begin{tabular}{@{}lccc}
\br
Metric               & Train & Valid & Test   \\
\mr
$R^2$ score (\%)     & 99.79    & 99.70      & 99.70  \\
RMSE                 & 0.0826   & 0.0988     & 0.0970 \\
\br
\end{tabular}
\end{table}

Last, to verify robustness against noise, Table \ref{tab6} includes the results under successive noise levels, and Figure \ref{fig12} depicts noisy reconstructions for the random sample. The modest decline in accuracy confirms that the network maintains reliable parameter recovery even when measurements are perturbed.

\noindent
\begin{minipage}{0.45\textwidth}
  \captionsetup{type=table}
  \captionof{table}{\label{tab6}Kite-shaped regression performance on noisy test data.}
  \centering
    \begin{tabular}{@{}lcc}
      \br
      Noise level & $R^2$ score (\%) & RMSE     \\
      \mr
      0.5\%       & 99.63             & 0.1045   \\
      1.0\%       & 99.44             & 0.1245   \\
      2.0\%       & 98.65             & 0.1832   \\
      5.0\%       & 92.94             & 0.4195   \\
      \br
    \end{tabular}
\end{minipage}%
\hfill
\begin{minipage}{0.48\textwidth}
  \centering
  \includegraphics[width=\textwidth]{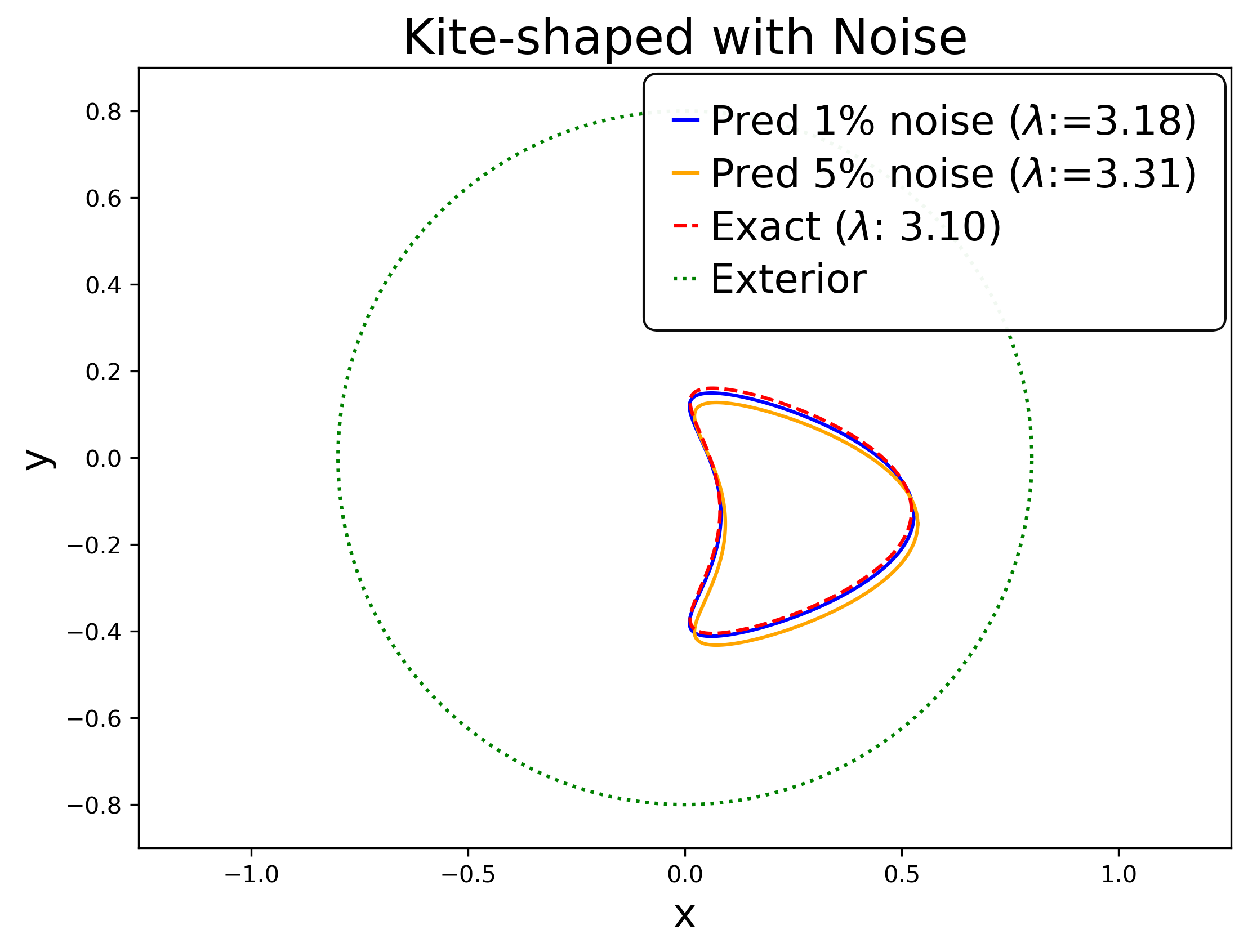}
  \captionof{figure}{Kite-shaped reconstructions for the randomly chosen sample in presence of increasing levels of noise.}
  \label{fig12}
\end{minipage}

\subsection{Inverse problem for star-shaped obstacles}

We now consider the star-shaped obstacles. As defined in (\ref{star}), their reconstruction requires estimating 14 parameters: the boundary coefficients \(\{\alpha_i\}_{i=0}^5\) and \(\{\beta_i\}_{i=1}^5\), the center coordinates \((x_0,y_0)\), and the impedance \(\lambda\). Compared to the peanut- and kite-shaped cases, this increased complexity requires deeper network architectures and larger training datasets, as we demonstrate in the experiments below.

To systematically investigate this setup, we divide our study into two sub-problems. First, we consider the star-shaped boundary reconstruction with fixed \(\lambda\), and then we examine the more general case of varying \(\lambda\).

\subsubsection{Fixed impedance}

We begin our analysis with star-shaped obstacles characterized by a fixed impedance $\lambda=2$. This choice simplifies the study for two reasons: first and most importantly, a fixed impedance makes the scattered far-fields dependent solely on the obstacle’s geometry; second, the inverse regression problem reduces to recovering 13 parameters instead of 14.

As expected, we initially considered a setting similar to the peanut- and kite-shaped problems: a dataset consisting of 30000 samples, with real and imaginary electric fields measured at 32 angle positions and a CNN model as in Tables \ref{tabAp2} and \ref{tabAp4}. However, this approach proved inadequate due to the increased complexity of the star-shaped geometry, producing poor regression performance as shown in Table \ref{tabAp6}, in the Appendix \ref{A41}. In particular, we discovered that the lowest-order coefficients, $\alpha_1$ and $\beta_1$ exhibited significant errors, causing the overall approach to fail. 

Therefore, we pursued two complementary strategies. First, we increased the model’s representational capacity to better tackle the problem complexity. Second, we expanded the training data by generating larger datasets with far-field measurements taken at more angular positions and by incorporating both electric and magnetic field components. In particular, we considered both the real and imaginary parts of the electric and magnetic fields, measured at \(T_0=128\) equidistant angles. This yields \(C_0=4\) field components (channels). We solved the direct problem to generate 80000 samples, each represented as \(\mathbf{x}_n\in\mathbb{R}^{512}\) and reshaped to \(\mathbf{X}_n\in\mathbb{R}^{128\times 4}\).
 
The deeper architecture consists of more convolutional layers and deeper dense layers in the MLP head. The details are given in Table \ref{tabAp7}, and the corresponding loss is shown in Figure \ref{figAp8}, (Appendix \ref{A41}). 

\begin{figure}[ht!]
    \centering
    \begin{minipage}[b]{0.48\textwidth}
        \centering
        \includegraphics[width=\textwidth]{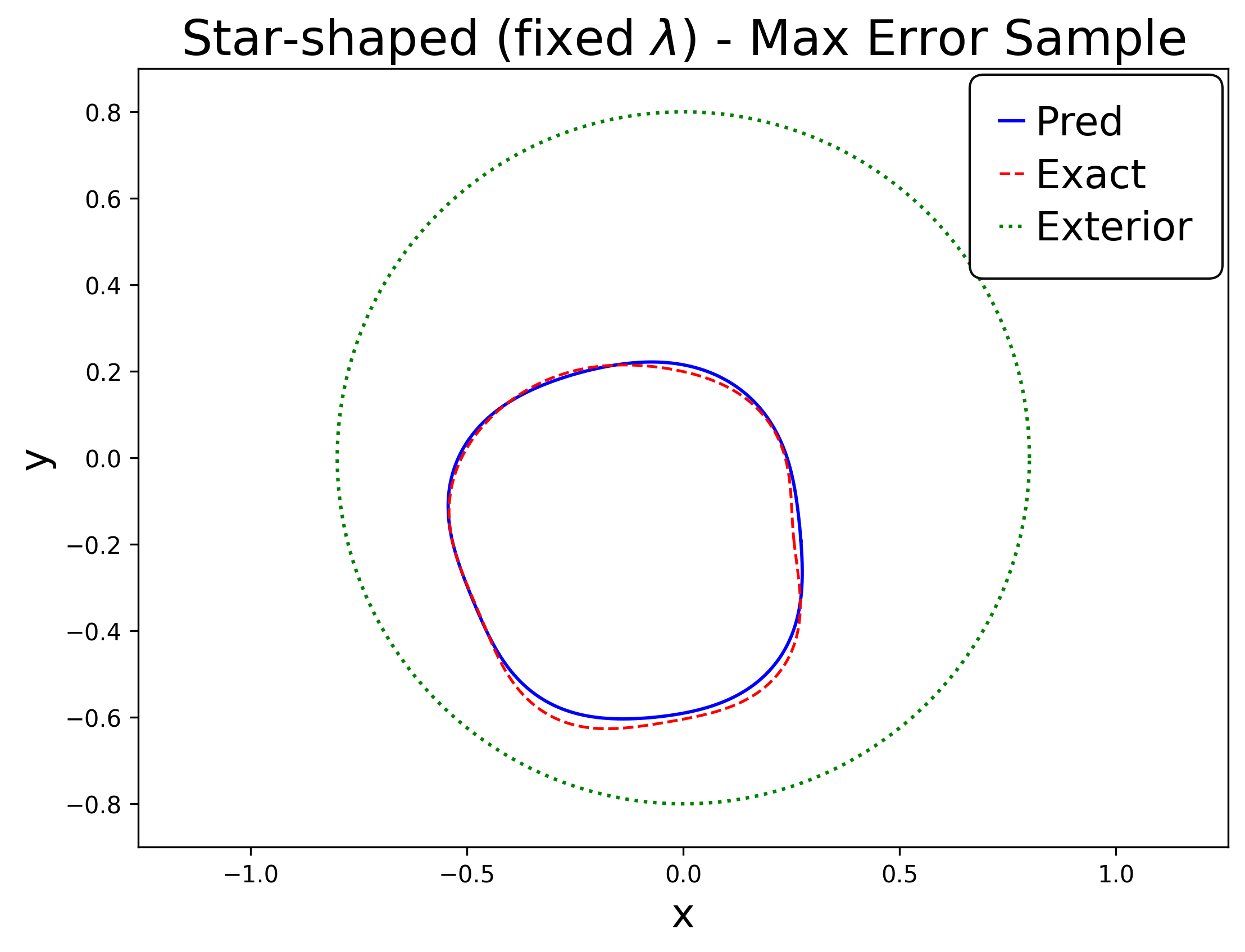} 
        \\[0.5ex]
            \end{minipage}\hfill
    \begin{minipage}[b]{0.48\textwidth}
        \centering
        \includegraphics[width=\textwidth]{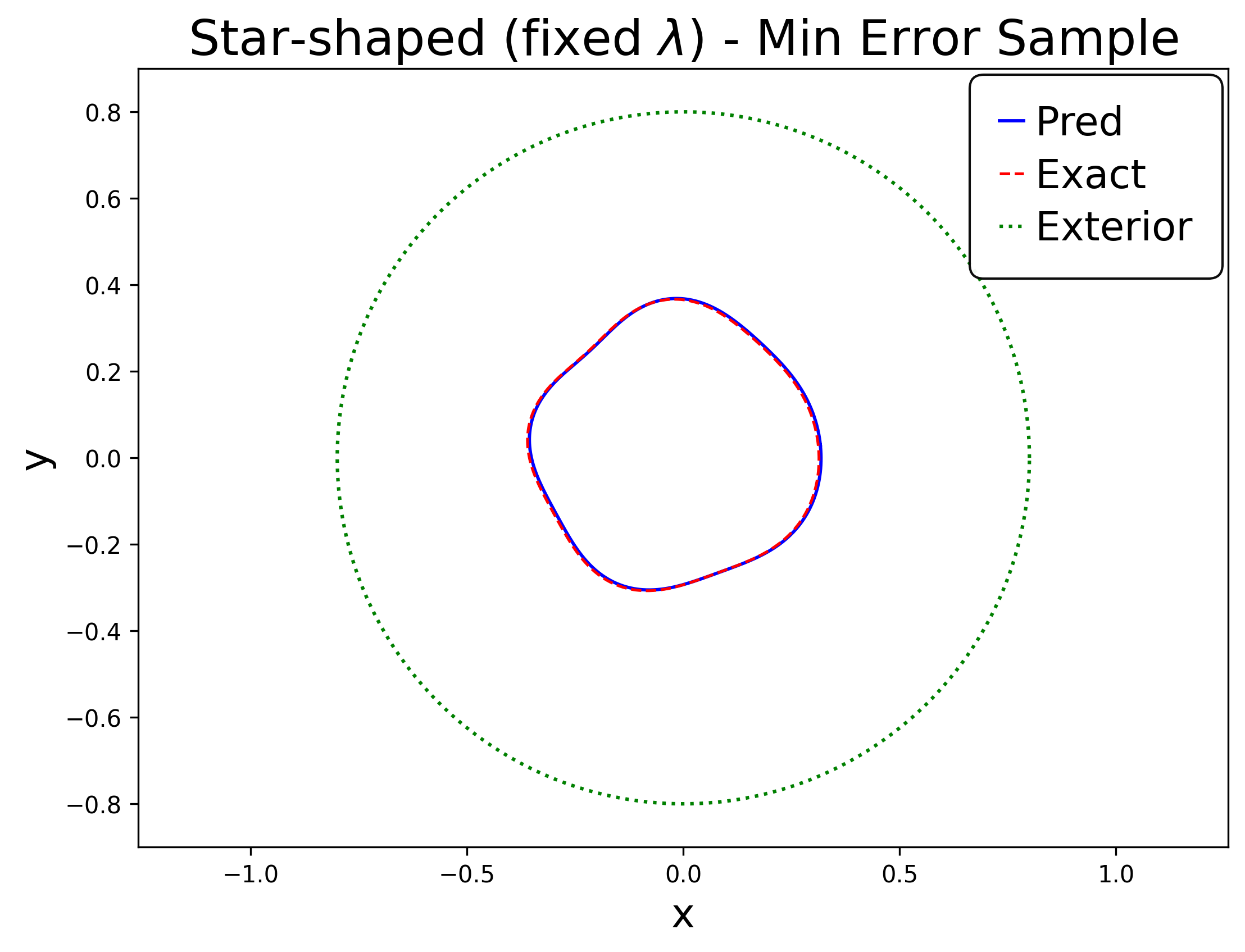} 
        \\[0.5ex]
         \end{minipage}
    \begin{minipage}[b]{0.48\textwidth}
        \centering
        \includegraphics[width=\textwidth]{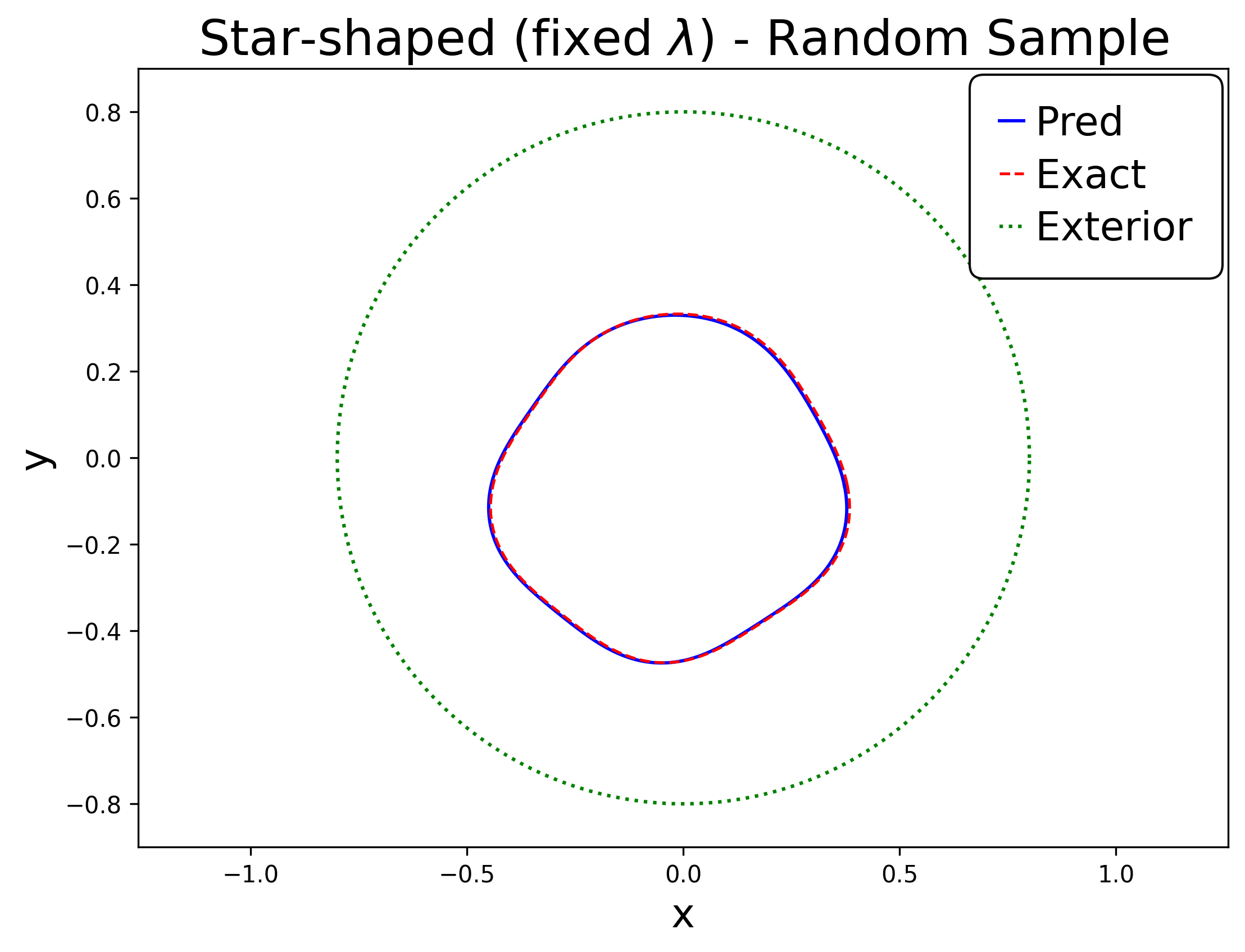} 
        \\[0.5ex]
    \end{minipage}

    \caption{Star-shaped (fixed $\lambda$) reconstructions: max error sample (top left),  min error sample (top right) and randomly chosen (bottom).}
    \label{fig13}
\end{figure}

We summarize the regression results across the star-shaped obstacles in Table \ref{tab7}, while Table \ref{tabAp9} breaks down test set performance for every individual parameter. We notice that the larger errors correspond to the  coefficients $\alpha_1$ and $\beta_1$.  In Figure \ref{fig13}, example boundary reconstructions are shown, presenting the samples with the largest and smallest prediction errors as well as a random case. To further illustrate model accuracy, Figure \ref{figAp9} plots the RMSE distribution over the test set. Overall, the high $R^2$ scores and low RMSE values confirm that the network reliably recovers the star-shaped boundary parameters.

\begin{table}
\caption{\label{tab7}Star-shaped (fixed $\lambda$) regression performance.}
\centering
\begin{tabular}{@{}lccc}
\br
Metric               & Train & Valid & Test     \\
\mr
$R^2$ score (\%)     & 99.62    & 99.15      & 99.14    \\
RMSE                 & 0.0105   & 0.0159     & 0.0159   \\
\br
\end{tabular}
\end{table}

As we did in the previous inverse problems, we evaluate robustness by introducing progressively higher levels of noise in the input measurements. Table \ref{tab8} reports the resulting regression metrics, and Figure \ref{fig14} illustrates the boundary reconstructions for the randomly chosen sample under noisy conditions. The decline in accuracy shows that the model recover parameters reliably up to a certain noise level.

\noindent
\begin{minipage}{0.45\textwidth}
  \captionsetup{type=table}
  \captionof{table}{\label{tab8}Star-shaped (fixed $\lambda$) regression performance on noisy test data.}
  \centering
    \begin{tabular}{@{}lcc}
      \br
      Noise level & $R^2$ score (\%) & RMSE     \\
      \mr
      0.5\%       & 98.66            & 0.0199   \\
      1.0\%       & 97.26            & 0.0285   \\
      2.0\%       & 92.02            & 0.0488   \\
      5.0\%       & 65.06            & 0.1021   \\
      \br
    \end{tabular}
\end{minipage}%
\hfill
\begin{minipage}{0.48\textwidth}
  \centering
  \includegraphics[width=\textwidth]{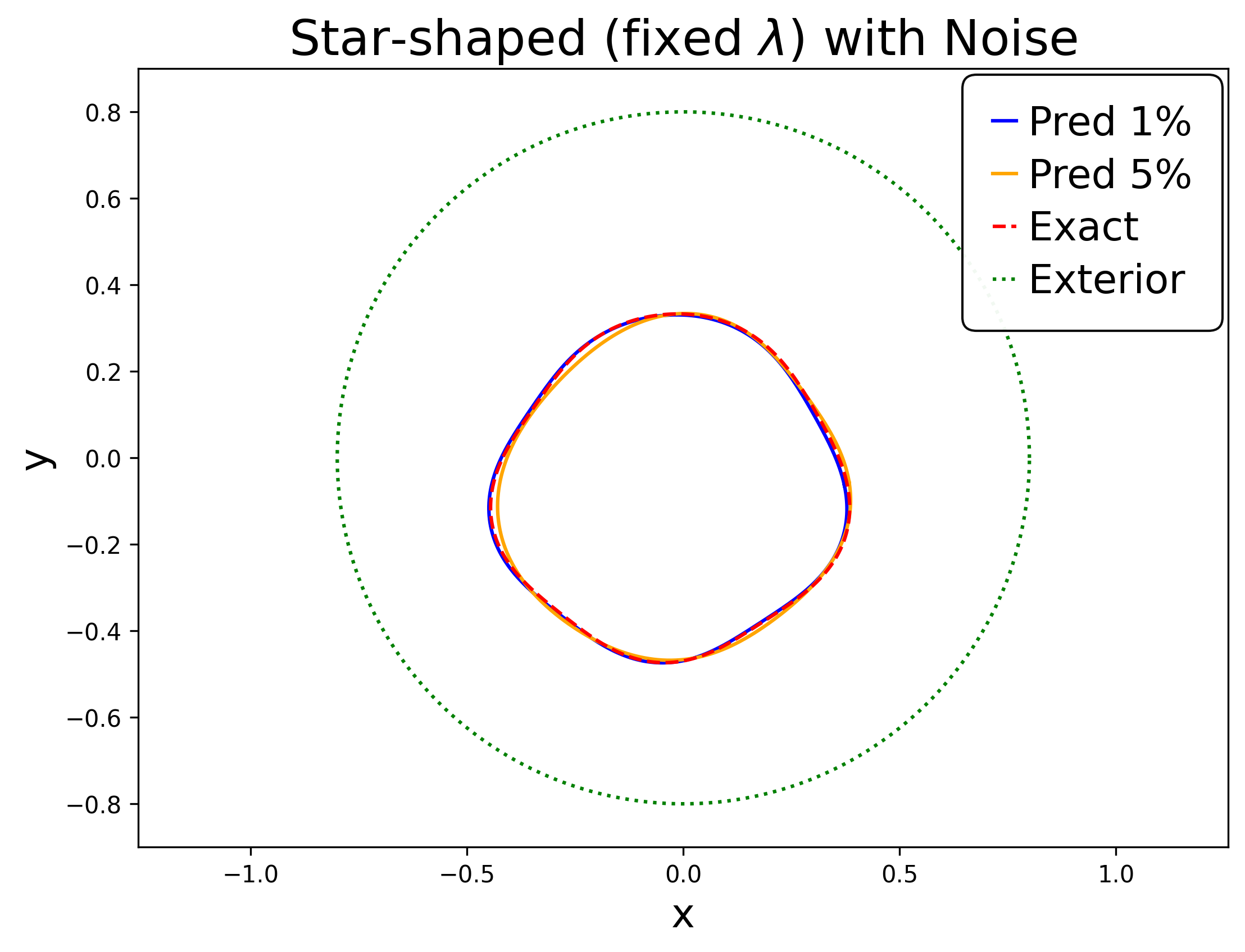}
  \captionof{figure}{Star-shaped (fixed $\lambda$) reconstructions for the randomly chosen sample in presence of increasing levels of noise.}
  \label{fig14}
\end{minipage}

\subsubsection{Variable impedance} 

Finally, we tackle the star-shaped obstacle inverse problem with variable impedance \(\lambda \in [0.1,10]\). Initially, we apply the same approach that succeeded in the fixed-\(\lambda\) case: a dataset of 80000 samples with \(T_0=128\) measurement angles for the \(C_0=4\) real and imaginary components of the electric and magnetic far-fields, under a single incident wave. However, the resulting parameter predictions using  the same CNN architecture as in fixed-$\lambda$ case, see Table \ref{tabAp7}, are poor in precision, as summarized in Table \ref{tabAp8} (Appendix \ref{A42}).  Analysis reveals that the lowest-order coefficients are particularly inaccurate, leading to high boundary reconstruction errors and degrading the overall model accuracy. 

As a result, we increase both the dataset size and the model’s complexity. Specifically, we generate 120000 samples from far-field measurements at \(T_0=128\) angles, of the real and imaginary parts of both electric and magnetic fields under two opposite incident directions (\(\phi=0\) and \(\phi=\pi\)). This yields \(C_0=8\) channels, so each input vector \(\mathbf{x}_n\in\mathbb{R}^{1024}\) is reshaped to \(\mathbf{X}_n\in\mathbb{R}^{128\times 8}\).

To better extract relevant features, we integrate an attention mechanism into our CNN (see Section \ref{attention}), inserting the attention layer immediately before the bottleneck convolution. The resulting architecture is detailed in Table \ref{tabAp10}, and its training loss curve appears in Figure \ref{figAp10}.

\begin{figure}[ht!]
    \centering
    \begin{minipage}[b]{0.48\textwidth}
        \centering
        \includegraphics[width=\textwidth]{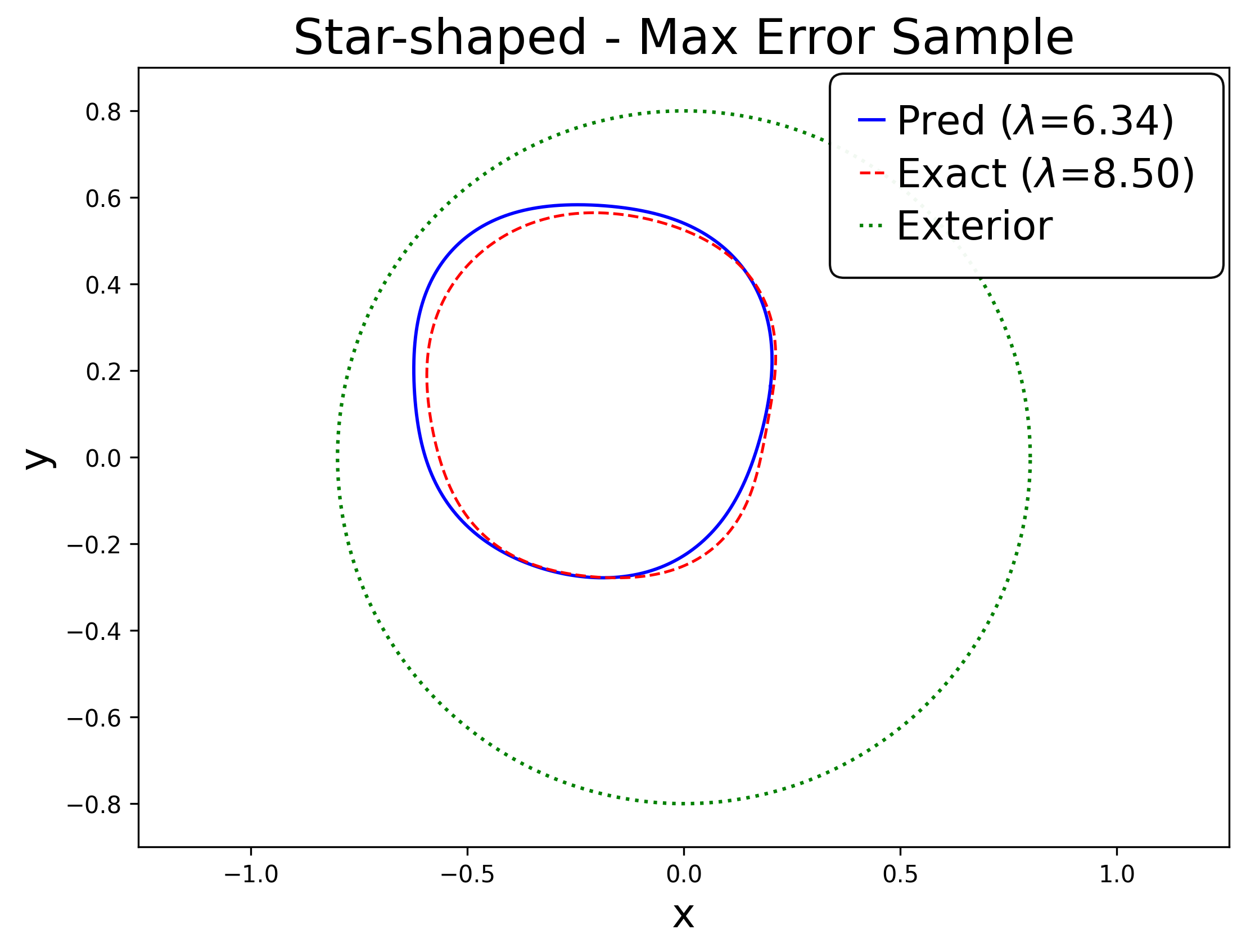} 
        \\[0.5ex]
            \end{minipage}\hfill
    \begin{minipage}[b]{0.48\textwidth}
        \centering
        \includegraphics[width=\textwidth]{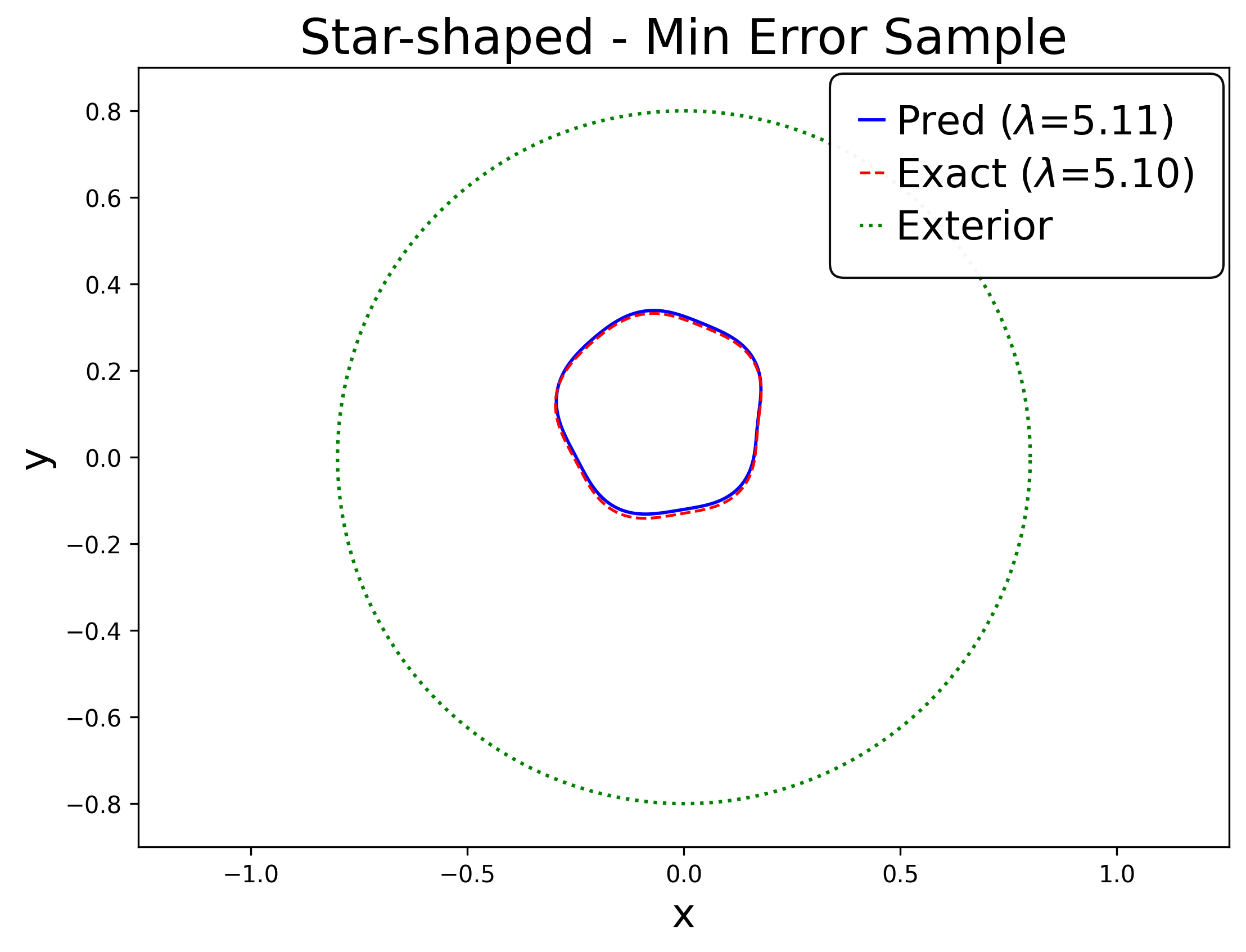} 
        \\[0.5ex]
         \end{minipage}
    \begin{minipage}[b]{0.48\textwidth}
        \centering
        \includegraphics[width=\textwidth]{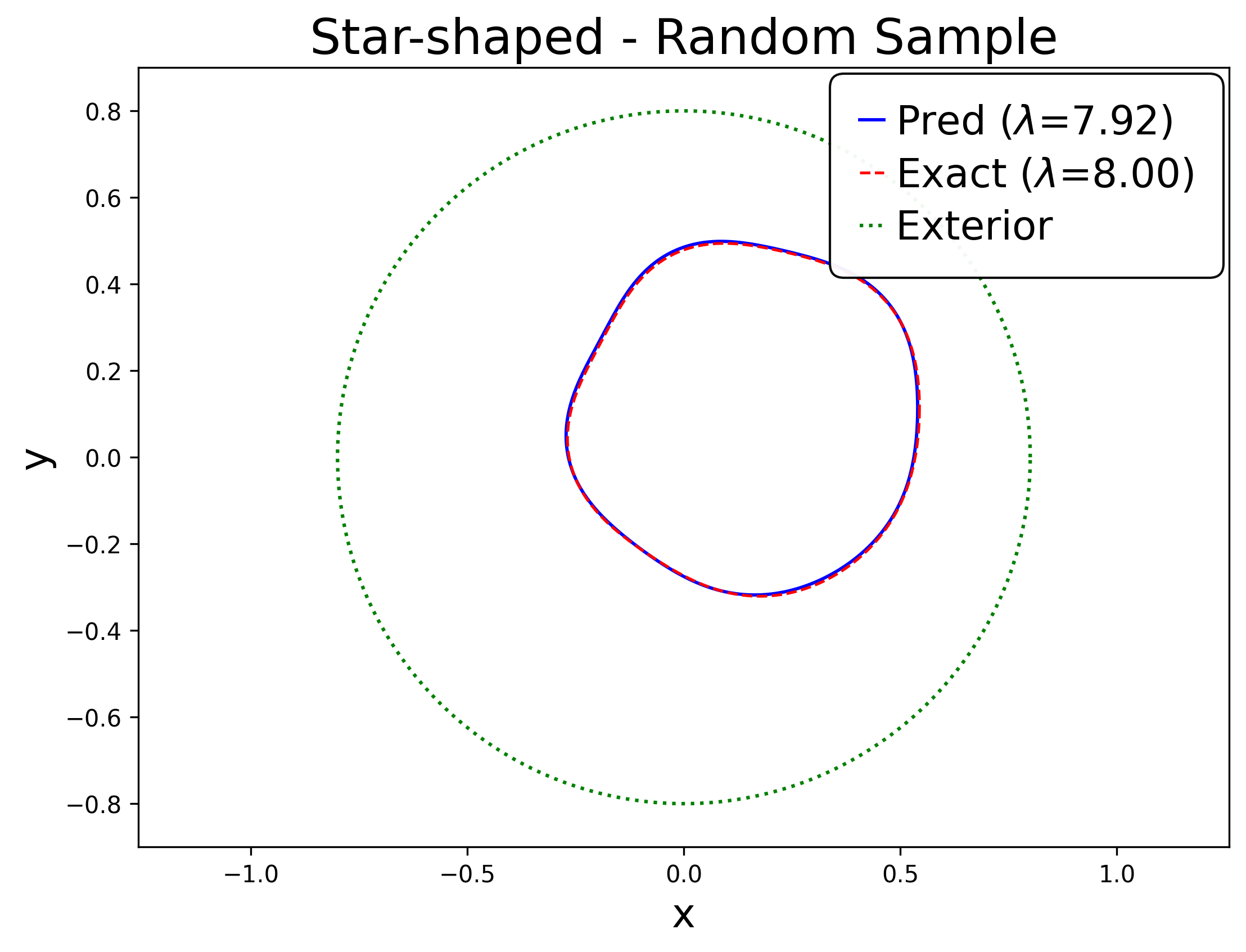} 
        \\[0.5ex]
    \end{minipage}

    \caption{Star-shaped reconstructions: max error sample (top left),  min error sample (top right) and randomly chosen (bottom).}
    \label{fig15}
\end{figure}

The regression results are summarized in Table \ref{tab9}, with per-parameter metrics detailed in Table \ref{tabAp11}. As in the fixed-$\lambda$ case, the first-order coefficients $\alpha_1$ and $\beta_1$ remain the most challenging to recover accurately. Example boundary reconstructions for the max-error, min-error, and a representative random sample appear in Figure \ref{fig15}. Figure \ref{figAp11} plots the predicted versus exact values for 
$\lambda$, and Figure \ref{figAp12} shows the RMSE histogram across the test set. Although the $R^2$
scores are lower than in the fixed-impedance case (suggesting mild overfitting), the model still achieves sufficient accuracy to solve the full star-shaped inverse problem reliably.

\begin{table}
\caption{\label{tab9}Star-shaped regression performance.}
\centering
\begin{tabular}{@{}lccc}
\br
Metric               & Train & Valid & Test   \\
\mr
$R^2$ score (\%)     & 99.05    & 96.32      & 96.02  \\
RMSE                 & 0.0469   & 0.0618     & 0.0617 \\
\br

\end{tabular}
\end{table}

We also evaluate the robustness of the model under noisy measurement conditions. Table \ref{tab10} summarizes the regression metrics, and Figure \ref{fig16} presents the boundary reconstruction for the random sample. The results demonstrate that the model loses stability more rapidly as noise increases compared to previous cases, highlighting the greater complexity of this inverse problem.

\noindent
\begin{minipage}{0.45\textwidth}
  \captionsetup{type=table}
  \captionof{table}{\label{tab10}Star-shaped regression performance on noisy test data.}
  \centering
    \begin{tabular}{@{}lcc}
      \br
      Noise level & $R^2$ score (\%) & RMSE     \\
      \mr
      0.5\%       & 95.25             & 0.0643   \\
      1.0\%       & 93.15             & 0.0722   \\
      2.0\%       & 86.93             & 0.1071   \\
      5.0\%       & 64.83             & 0.2541   \\
      \br
    \end{tabular}
\end{minipage}%
\hfill
\begin{minipage}{0.48\textwidth}
  \centering
  \includegraphics[width=\textwidth]{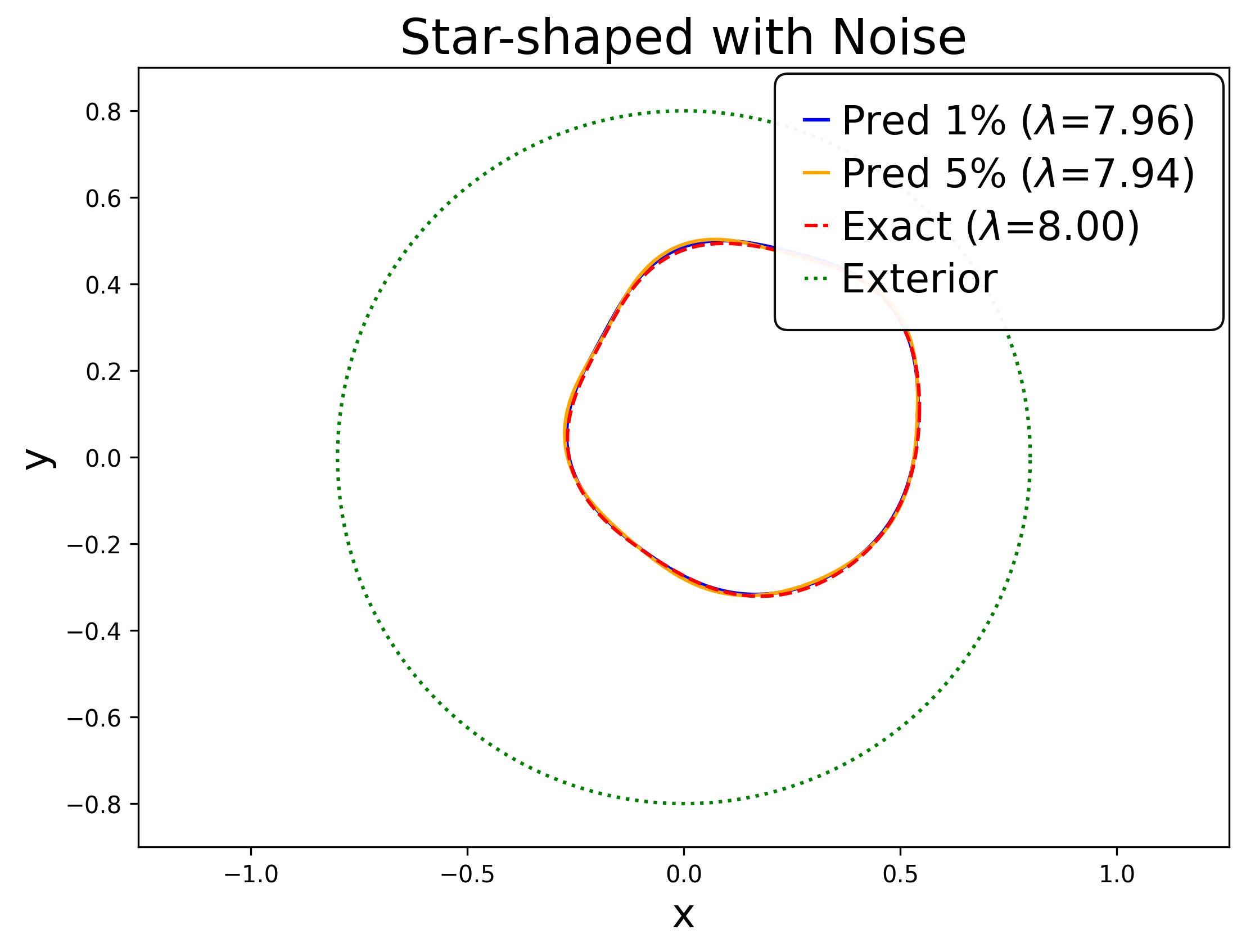}
  \captionof{figure}{Star-shaped reconstructions for the randomly chosen sample in presence of increasing levels of noise.}
  \label{fig16}
\end{minipage}

\subsection{Inverse problem for misclassified obstacles}

We now examine how misclassified obstacles behave when processed through the reconstruction pipeline. This idea, analyzing misclassified samples within a divide-and-conquer framework, has also been explored in recent works such as \cite{MENG2025116525}. 

As shown by the classification metrics in Table \ref{tab1}, the accuracy of the test set is equal to 0.988. With a test set of 9000 samples (10\% of the 90000–sample dataset considered in Section \ref{class}), this corresponds to 85 kites misclassified as peanuts and 18 peanuts misclassified as kites. 

Our primary question is how much misclassifications influence the subsequent boundary and impedance reconstructions. Our analysis shows that most misclassified cases arise when the true obstacle shape lies near the decision boundary--i.e., when the geometry is sufficiently smooth or ``intermediate$"$ so that multiple classes can be chosen.

\begin{figure}[ht!]
    \centering
    \begin{minipage}{0.48\textwidth}
        \centering
        \includegraphics[width=\textwidth]{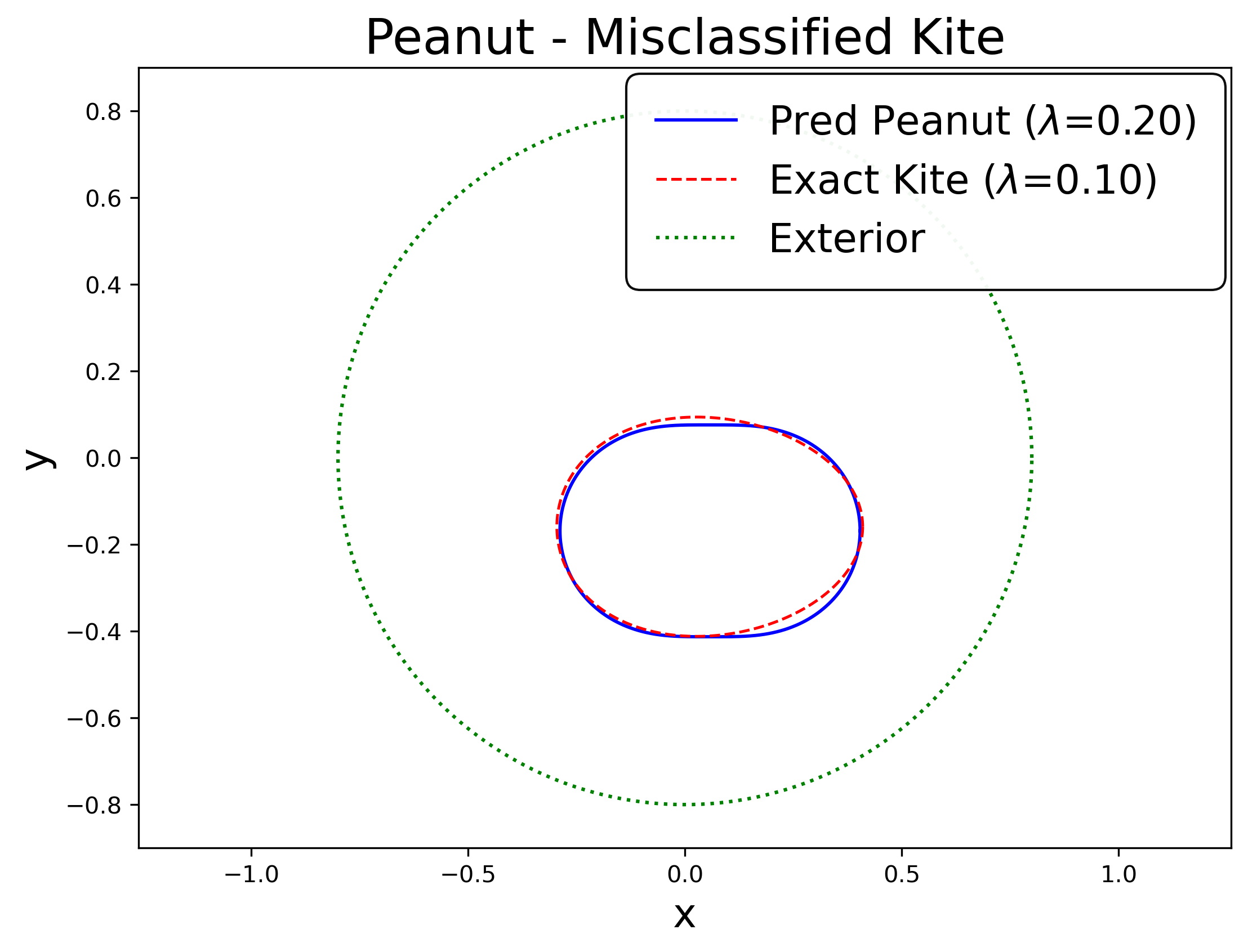}
    \end{minipage}
    \hfill
    \begin{minipage}{0.48\textwidth}
        \centering
        \includegraphics[width=\textwidth]{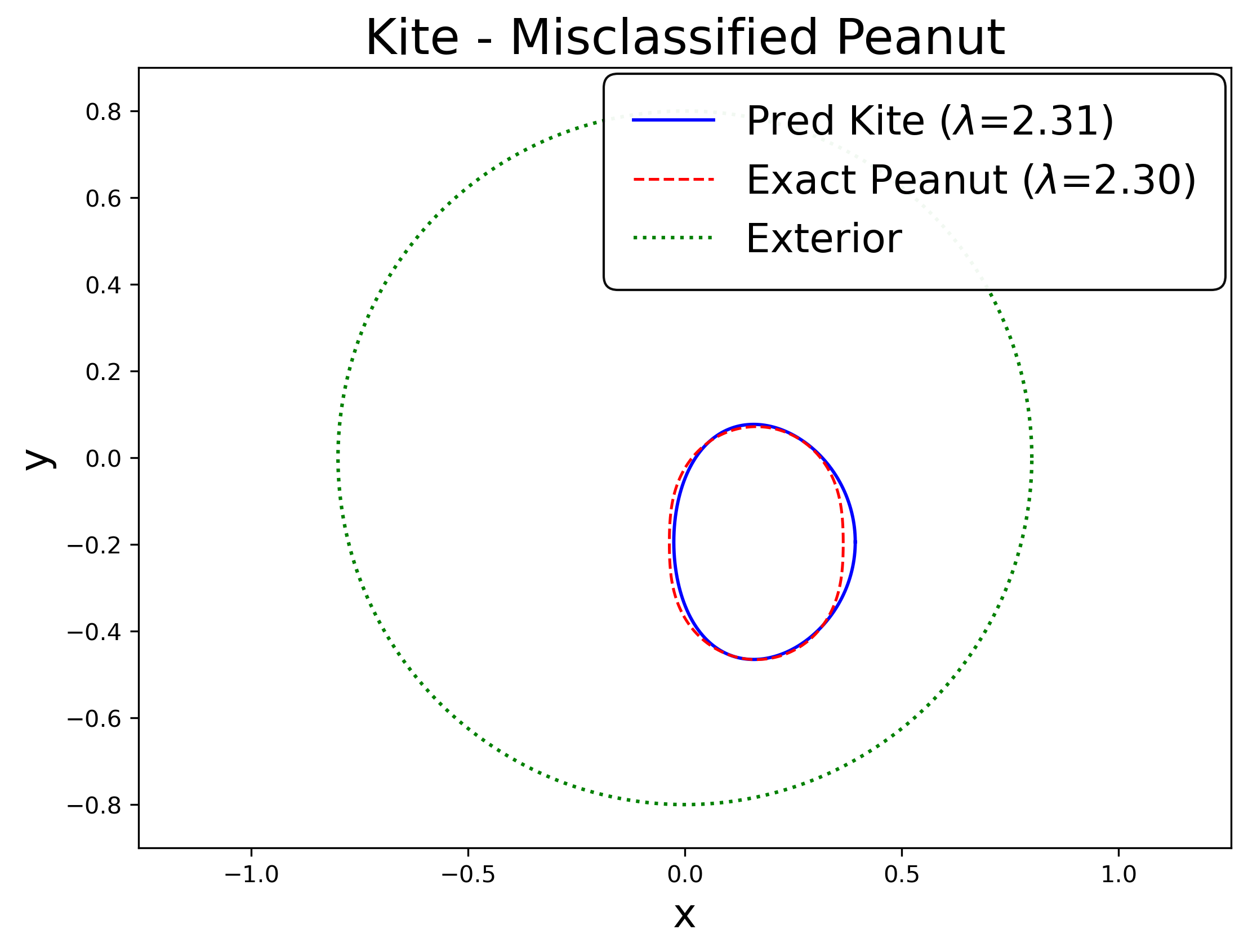}
    \end{minipage}
    \caption{Misclassification cases: kite reconstructed as peanut (left) and peanut reconstructed as kite (right).}
    \label{fig17}
\end{figure}

Figure \ref{fig17} (left) illustrates one such case: a true kite that the classifier labeled it as a peanut. Applying the peanut‐model from Section \ref{peanCNN}, we recover the obstacle center and area accurately, and both the predicted impedance 
$\lambda$ and reconstructed boundary closely match the true values. Similarly, Figure \ref{fig17} (right) shows a peanut misclassified as a kite. Using the kite‐model from Section \ref{kiteCNN}, the reconstructed center, area, and $\lambda$ again match the ground truth closely, and the boundary approximation is acceptable. 

These observations indicate that, even when classification errors occur, the regression-based shape reconstruction model can still yield reliable results, provided the true geometry lies near the classifier’s decision boundary.  

\section{Conclusions} \label{concl}

This work introduced a novel divide-and-conquer methodology for the inverse electromagnetic scattering problem to reconstruct the characteristics of an impedance cylinder, coated by a magneto-dielectric cylinder. Using multi-channel CNNs with circular padding, we first classified obstacles according to their shape and then regressed both the boundary coefficients and the impedance constant $\lambda$. This approach shared the following two main advantages: it was non-invasive, since it was based only on limited information from far-field data; and went beyond simply connected domains, addressing the additional challenge of jointly recovering geometry and impedance.

Training data were generated via direct scattering formulations, providing measurements of real and imaginary parts of the electric and magnetic far-fields. In most cases, measurements from a single incident field were sufficient.  Classification results demonstrated high accuracy across all classes, with slight degradation under increasing noise levels. Regression on peanut- and kite-shaped obstacles achieved excellent $R^2$ scores and low RMSE, maintaining robustness to moderate noise. For star-shaped obstacles, the fixed-$\lambda$ case required deeper networks and more extensive sampling to attain reliable recovery of 13 boundary coefficients. Allowing  $\lambda$ to vary further increased complexity, necessitating additional training samples from two incident waves and an attention mechanism. Despite a modest drop in precision compared to the fixed-impedance scenario, the final model still delivered satisfactory boundary and impedance estimates. In all star-shaped models considered, we observed that the lowest-order boundary coefficients $\alpha_1$ and $\beta_1$ (corresponding to the first sinusoidal mode) were more challenging to recover accurately. This indicated that, since these coefficients change the object's boundary smoothly, their influence on the measured data is subtle, making them inherently harder to estimate than sharper, higher-order boundary variations. Finally, we noticed that, even when the classifier mislabels an obstacle, the corresponding shape reconstruction model accurately recovered both the geometry and the impedance for samples near decision boundaries. 

Further research could explore the inverse problem of a non-constant impedance $\lambda \in C^1(\Gamma_1)$. This extension demands deeper, more expressive DNN models, able to manage the additional complexity. Incorporating further unknown material parameters, such as the permittivity and permeability of the magneto-dielectric cylinder, or  scenarios where only phase-less data are available, would be an interesting direction. Finally, developing a rigorous operator-theoretic framework to analyze the approximation properties and convergence of the proposed CNN approach is a valuable task, though it lies beyond the scope of the present manuscript and will be addressed in future work. 

\section*{Data availability statement}
The source codes for the implementation of the CNN models will be available upon acceptance, in the GitHub repository of the corresponding author \href{https://github.com/npallikarakis}{https://github.com/npallikarakis}. Data will be available upon reasonable request.

\section*{Acknowledgments}
This work was supported by computational time granted to the corresponding author from the National Infrastructures for Research and Technology S.A. (GRNET  S.A.) in the National HPC facility – ARIS under the project ``Application of Machine Learning to Inverse Problems$"$ (AMALIP). 

\bibliographystyle{siam}
\bibliography{refs}

\clearpage

\appendix

\section{Computational details}
\label{A}

This appendix presents the implementation details and performance metrics of all computational experiments. Each subsection corresponds to one of the inverse problems discussed in Section \ref{inverse}, covering the classification and regression problems for peanut-, kite-, and star-shaped obstacle reconstructions.

\subsection{Classification of obstacles}
\label{A1}

\begin{table}[h!]
\caption{Circular CNN architecture for obstacle classification.}
\centering
\begin{tabular}{@{}lp{0.7\textwidth}}
\br
Layer & Configuration \\
\mr
Input & $(T_0=32, C_0=2)$ tensor of real/imaginary electric fields \\
CircularConv1D$^{(1)}$ & $N_f^{(1)}=64$ filters, $K^{(1)}=5$, $S^{(1)}=1$, Swish activation; output: $(T_1=32, C_1=64)$ \\
CircularConv1D$^{(2)}$ & $N_f^{(2)}=64$ filters, $K^{(2)}=5$, $S^{(2)}=2$, Swish activation; output: $(T_2=16, C_2=64)$ (downsampling)\\
CircularConv1D$^{(3)}$ & $N_f^{(3)}=64$ filters, $K^{(3)}=7$, $S^{(3)}=1$, Swish activation; output: $(T_3=16, C_3=64)$ \\
Pointwise Conv1D & $N_{\mathrm{b}}=16$ filters, $K=1$, Swish activation; output: $(T_3, N_{\mathrm{b}}=16)$ (bottleneck) \\
Flatten & $d_0 = T_3 \times N_{\mathrm{b}} = 256$ \\
Dense & $d_1=128$ nodes, Swish activation; LayerNorm; Dropout ($p=0.2$) \\
Dense & $d_2=64$ nodes, Swish activation; LayerNorm; Dropout ($p=0.1$) \\
Output Dense & $d_{\mathrm{out}}=3$ units, Softmax activation \\
\mr
Task Type & Classification\\
Optimizer & Adam, learning rate $\alpha=1\times10^{-5}$ \\
Loss Function & Categorical Cross-entropy (\(\mathcal{L}_{\mathrm{class}}\)) \\
Metrics & Accuracy \\
Training & Batch size $=64$; EarlyStopping: $\Delta=10^{-3}$, Patience$=150$; restore best weights; 7 sec/epoch on GPU \\
\br
\end{tabular}
\label{tabAp1}
\end{table}

\begin{figure}[h!]
    \centering
      \begin{minipage}{0.48\textwidth}
        \centering
    \includegraphics[width=\textwidth]{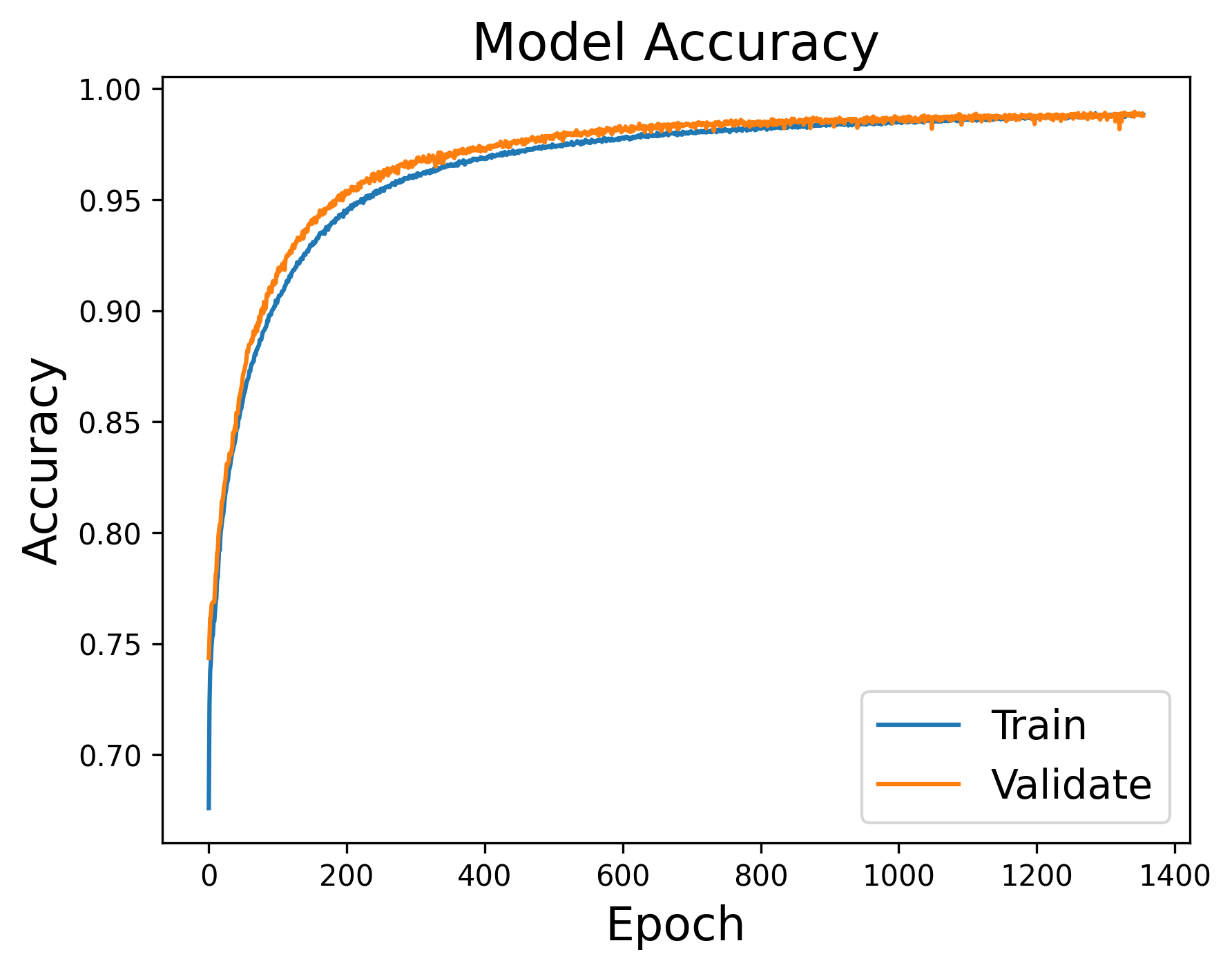} 
    \end{minipage}\hfill
    \begin{minipage}{0.48\textwidth}
        \centering
    \includegraphics[width=\textwidth]{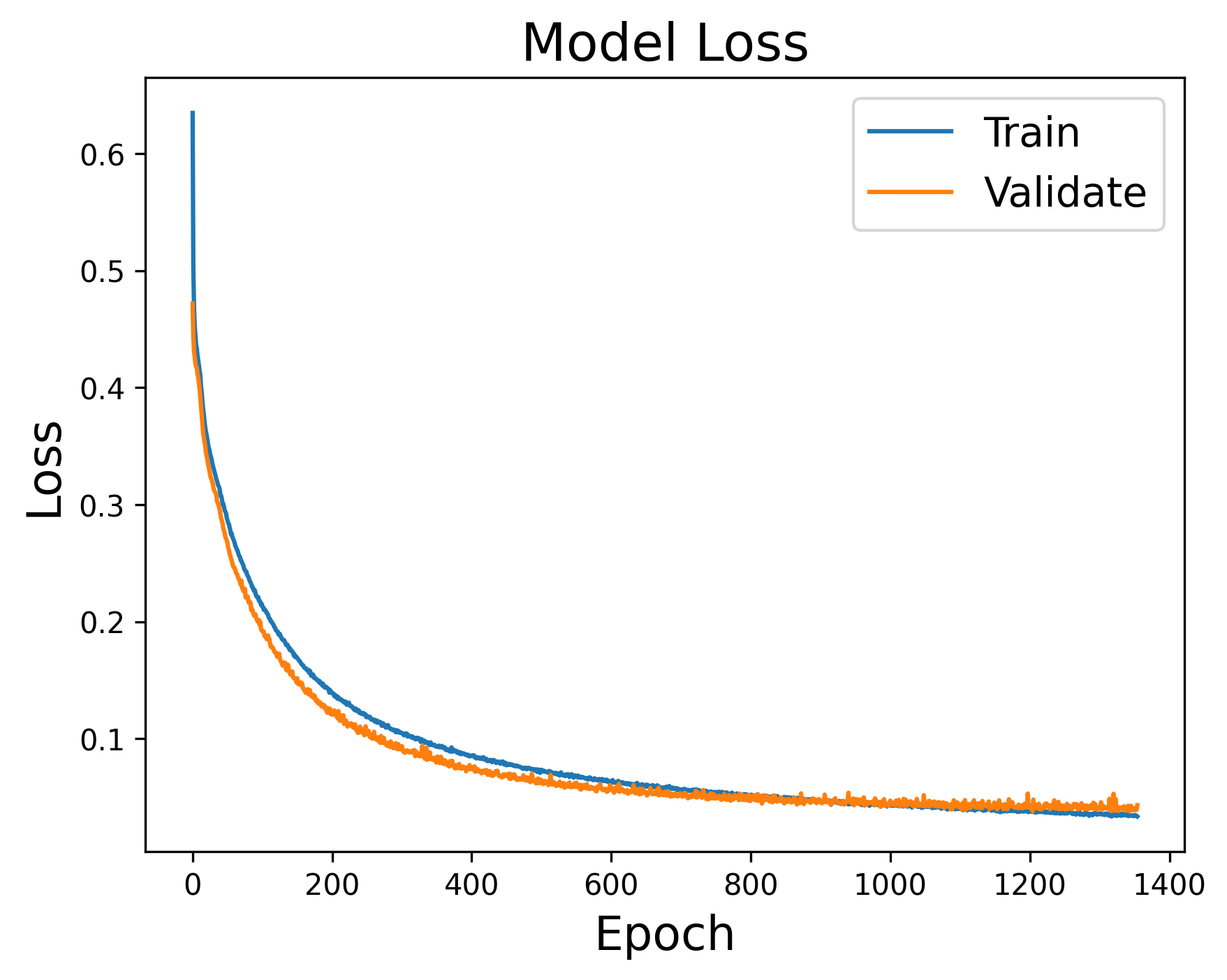} 
    \end{minipage}
    \caption{The training and validation accuracy (left) and loss (right) of the classification network. }
    \label{figAp1}
\end{figure}

\FloatBarrier
\subsection{Inverse problem for peanut-shaped obstacles} 
\label{A2}

\begin{table}[h!]
\caption{Circular CNN architecture for peanut-shaped obstacles.}
\begin{tabular}{@{}lp{0.7\textwidth}}
\br
Layer & Configuration \\
\mr
Input & $(T_0=32, C_0=2)$ tensor of real/imaginary electric fields \\
CircularConv1D$^{(1)}$ & $N_f^{(1)}=64$ filters, $K^{(1)}=5$, $S^{(1)}=1$, Swish activation; output: $(T_1=32, C_1=64)$ \\
CircularConv1D$^{(2)}$ & $N_f^{(2)}=64$ filters, $K^{(2)}=5$, $S^{(2)}=2$, Swish activation; output: $(T_2=16, C_2=64)$ (downsampling) \\
Pointwise Conv1D & $N_{\mathrm{b}}=16$ filters, $K=1$, Swish activation; output: $(T_2, N_{\mathrm{b}}=16)$ (bottleneck) \\
Flatten & $d_0 = T_2 \times N_{\mathrm{b}} = 256$ \\
Dense & $d_1=64$ nodes, Swish activation; LayerNormalization \\
Output Dense & $d_{\mathrm{out}}=5$ units, no activation\\
\mr
Task Type & Regression\\
Optimizer & Adam, learning rate $\alpha=1\times10^{-4}$ \\
Loss Function & Mean Squared Error (\(\mathcal{L}_{\mathrm{reg}}\)) \\
Metrics & Root Mean Squared Error (RMSE) \\
Training & Batch size $=128$; EarlyStopping: $\Delta=10^{-4}$, Patience$=80$; restore best weights; 3 sec/epoch on CPU \\
\br
\end{tabular}
\label{tabAp2}
\end{table}

\noindent
\begin{minipage}[t]{0.48\textwidth}
    \centering        
    \includegraphics[width=\textwidth]{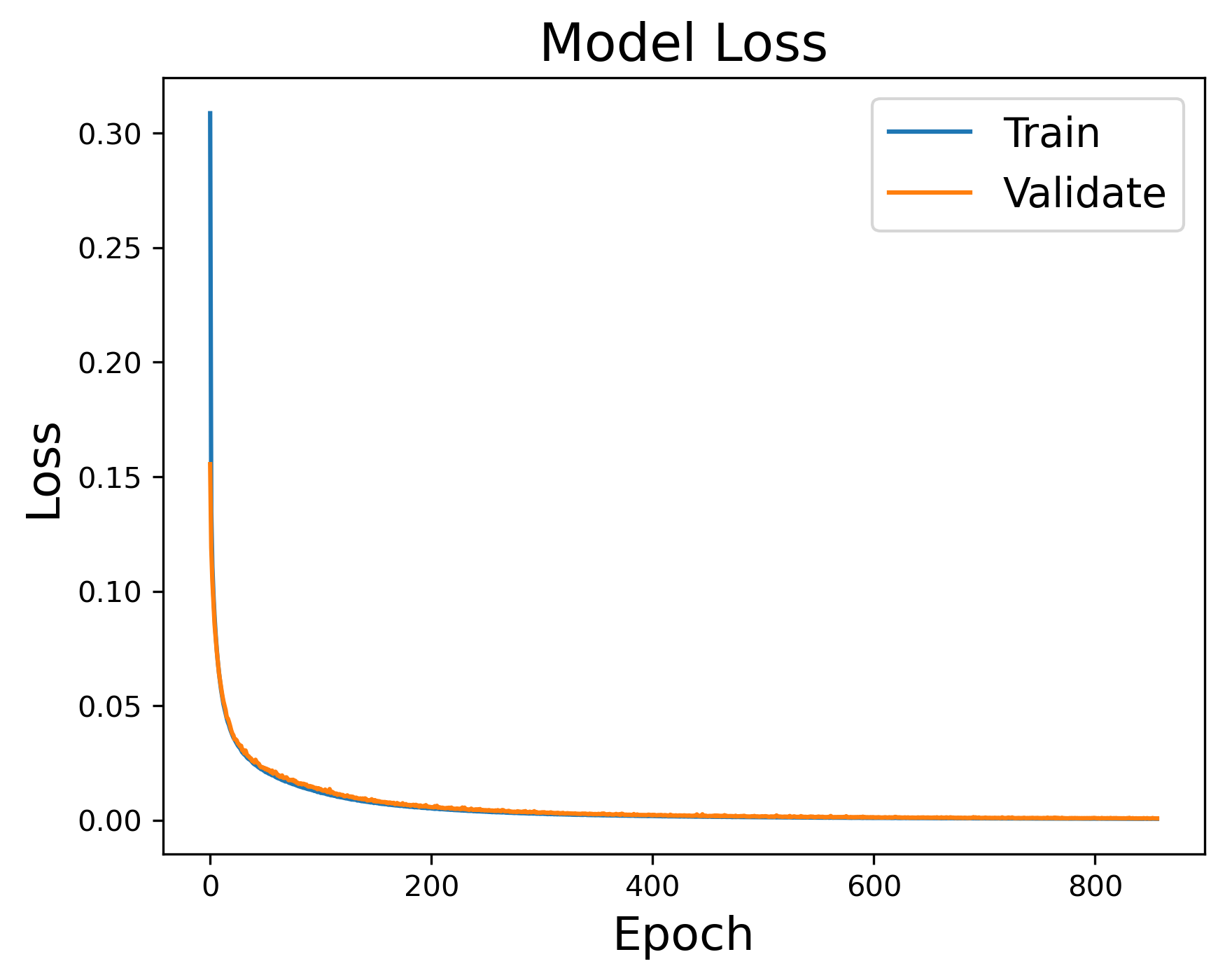} 
    \captionof{figure}{The training and validation loss of the network on peanut-shaped obstacles.}
    \label{figAp2}
\end{minipage}%
\hfill
\begin{minipage}[t]{0.48\textwidth}
    \centering        
    \includegraphics[width=\textwidth]{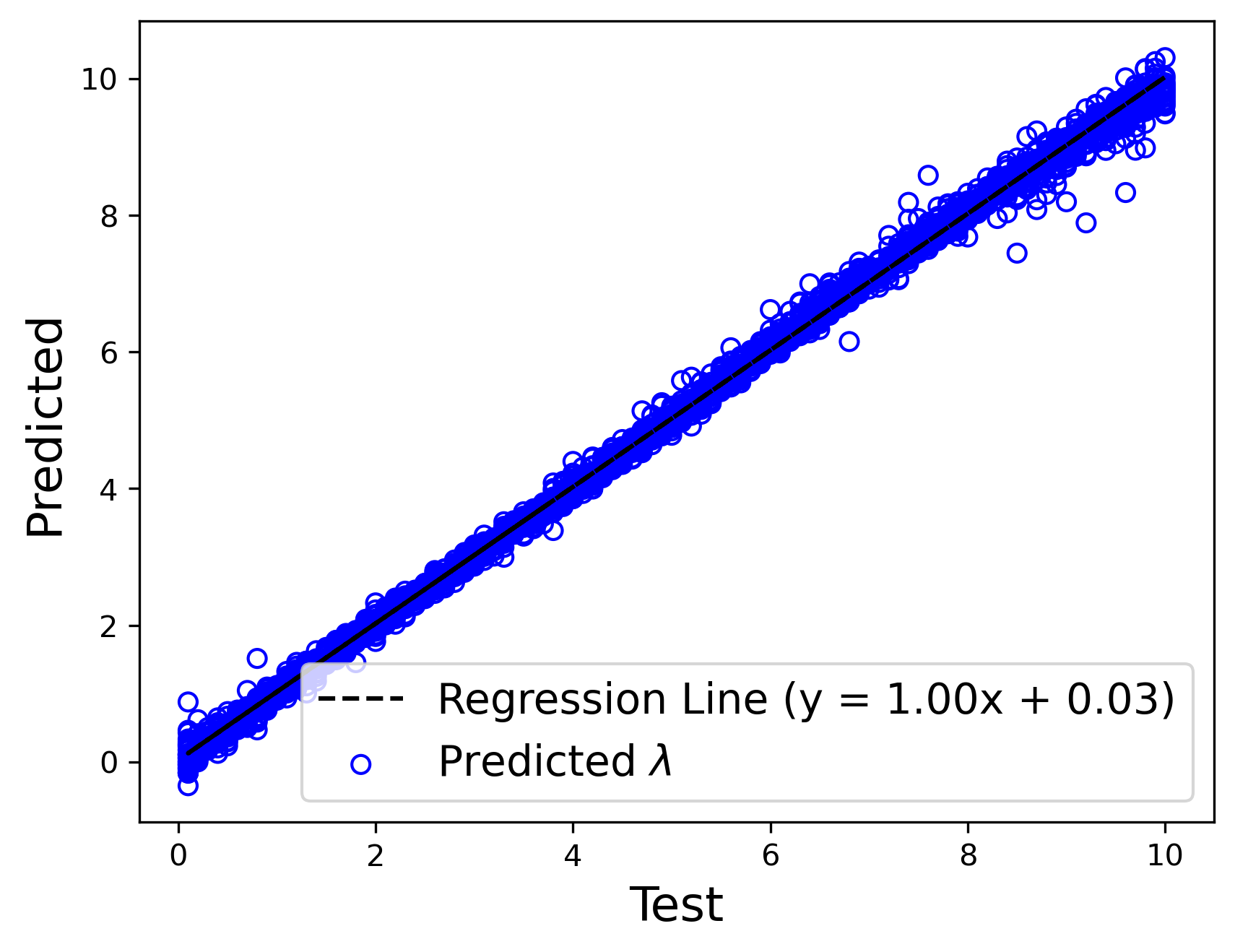} 
    \captionof{figure}{Regression results for the estimated impedance \(\lambda\) of peanut-shaped obstacles.}
    \label{figAp3}
\end{minipage}

\noindent
\begin{minipage}{0.45\textwidth}
  \captionsetup{type=table}
  \captionof{table}{\label{tabAp3}Test set performance by coefficient for the peanut-shaped obstacles.}
  \centering
    \begin{tabular}{@{}lcc}
      \br
      Coefficient & $R^2$ score (\%) & RMSE    \\
      \mr
      $\alpha$    & 99.94             & 0.0012  \\
      $\beta$     & 99.93             & 0.0012  \\
      $x_0$       & 99.96             & 0.0022  \\
      $y_0$       & 99.97             & 0.0020  \\
      $\lambda$   & 99.79             & 0.1329  \\
      \br
    \end{tabular}
\end{minipage}%
\hfill
\begin{minipage}{0.48\textwidth}
  \centering        
  \includegraphics[width=\textwidth]{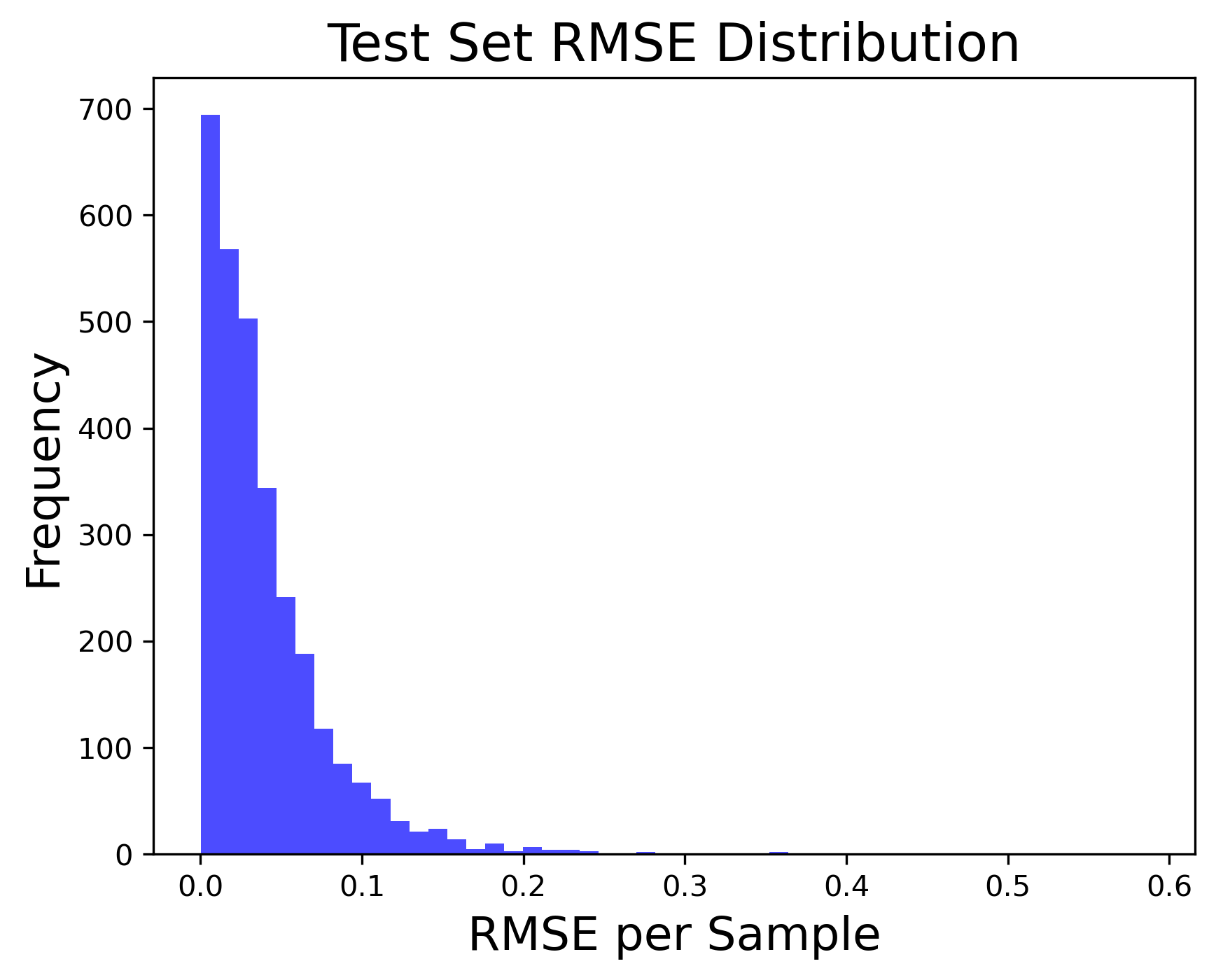} 
  \captionof{figure}{Distribution of the RMSE per sample in the test set of the peanut-shaped obstacles.}
  \label{figAp4}
\end{minipage}
\newpage
\subsection{Inverse problem for kite-shaped obstacles}
\label{A3}


\begin{table}[h!]
\caption{Circular CNN architecture for kite-shaped obstacles.}
\begin{tabular}{@{}lp{0.7\textwidth}}
\br
Layer & Configuration \\
\mr
Input & $(T_0=32, C_0=2)$ tensor of real/imaginary electric fields \\
CircularConv1D$^{(1)}$ & $N_f^{(1)}=64$ filters, $K^{(1)}=5$, $S^{(1)}=1$, Swish activation; output: $(T_1=32, C_1=64)$ \\
CircularConv1D$^{(2)}$ & $N_f^{(2)}=64$ filters, $K^{(2)}=5$, $S^{(2)}=2$, Swish activation; output: $(T_2=16, C_2=64)$ (downsampling) \\
Pointwise Conv1D & $N_{\mathrm{b}}=16$ filters, $K=1$, Swish activation; output: $(T_2, N_{\mathrm{b}}=16)$ (bottleneck) \\
Flatten & $d_0 = T_2 \times N_{\mathrm{b}} = 256$ \\
Dense & $d_1=64$ nodes, Swish activation; LayerNormalization \\
Output Dense & $d_{\mathrm{out}}=6$ units, no activation\\
\mr
Task Type & Regression \\
Optimizer & Adam, learning rate $\alpha=1\times10^{-4}$ \\
Loss Function & Mean Squared Error (\(\mathcal{L}_{\mathrm{reg}}\)) \\
Metrics & Root Mean Squared Error (RMSE) \\
Training & Batch size $=128$; EarlyStopping: $\Delta=10^{-4}$, Patience$=80$; restore best weights; 3 sec/epoch on CPU \\
\br
\end{tabular}
\label{tabAp4}
\end{table}

\noindent
\begin{minipage}[t]{0.48\textwidth}
    \centering        
    \includegraphics[width=\textwidth]{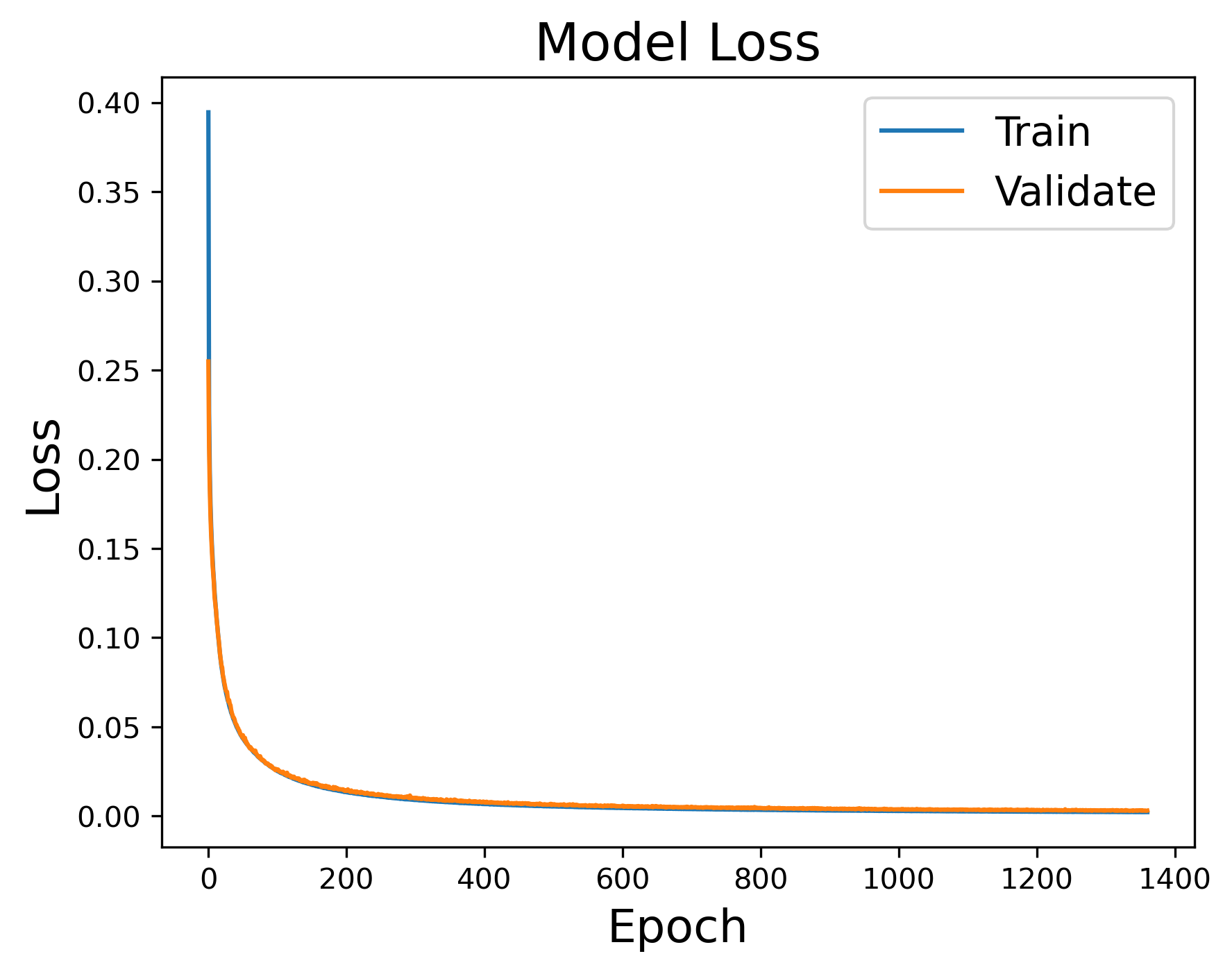} 
    \captionof{figure}{The training and validation loss of the network on kite-shaped obstacles.}
    \label{figAp5}
\end{minipage}%
\hfill
\begin{minipage}[t]{0.48\textwidth}
    \centering        
    \includegraphics[width=\textwidth]{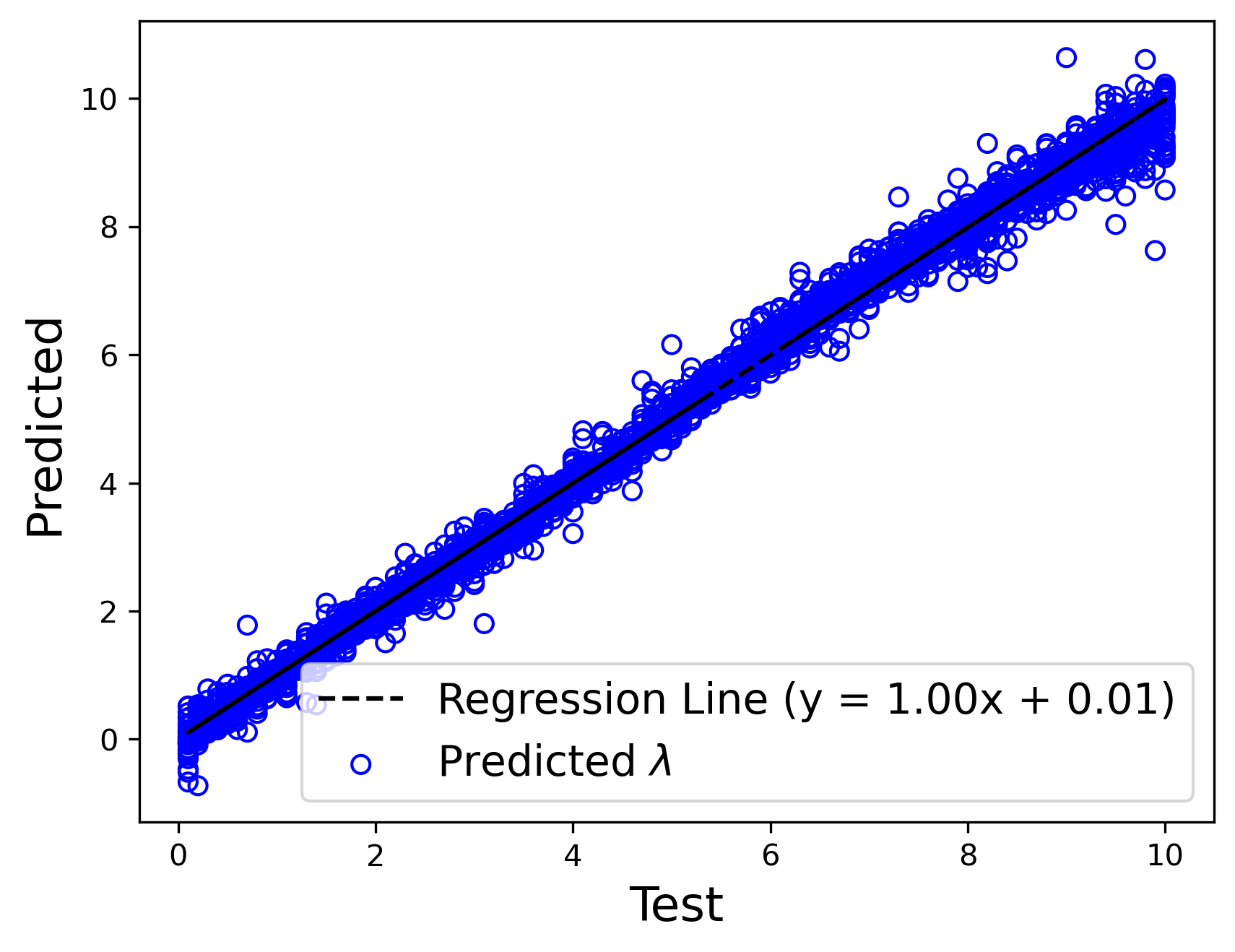} 
    \captionof{figure}{Regression results for the estimated impedance \(\lambda\) of kite-shaped obstacles.}
    \label{figAp6}
\end{minipage}
\medskip
\begin{minipage}{0.45\textwidth}
  \captionsetup{type=table}
  \captionof{table}{\label{tabAp5}Test set performance by coefficient for the kite-shaped obstacles.}
  \centering
    \begin{tabular}{@{}lcc}
      \br
      Coefficient & $R^2$ score (\%) & RMSE    \\
      \mr
      $\alpha$    & 99.64             & 0.0037  \\
      $\beta$     & 99.73             & 0.0042  \\
      $\gamma$    & 99.72             & 0.0032  \\
      $x_0$       & 99.87             & 0.0043  \\
      $y_0$       & 99.92             & 0.0032  \\
      $\lambda$   & 99.32             & 0.2375  \\
      \br
    \end{tabular}
\end{minipage}%
\hfill
\begin{minipage}{0.48\textwidth}
  \centering        
  \includegraphics[width=\textwidth]{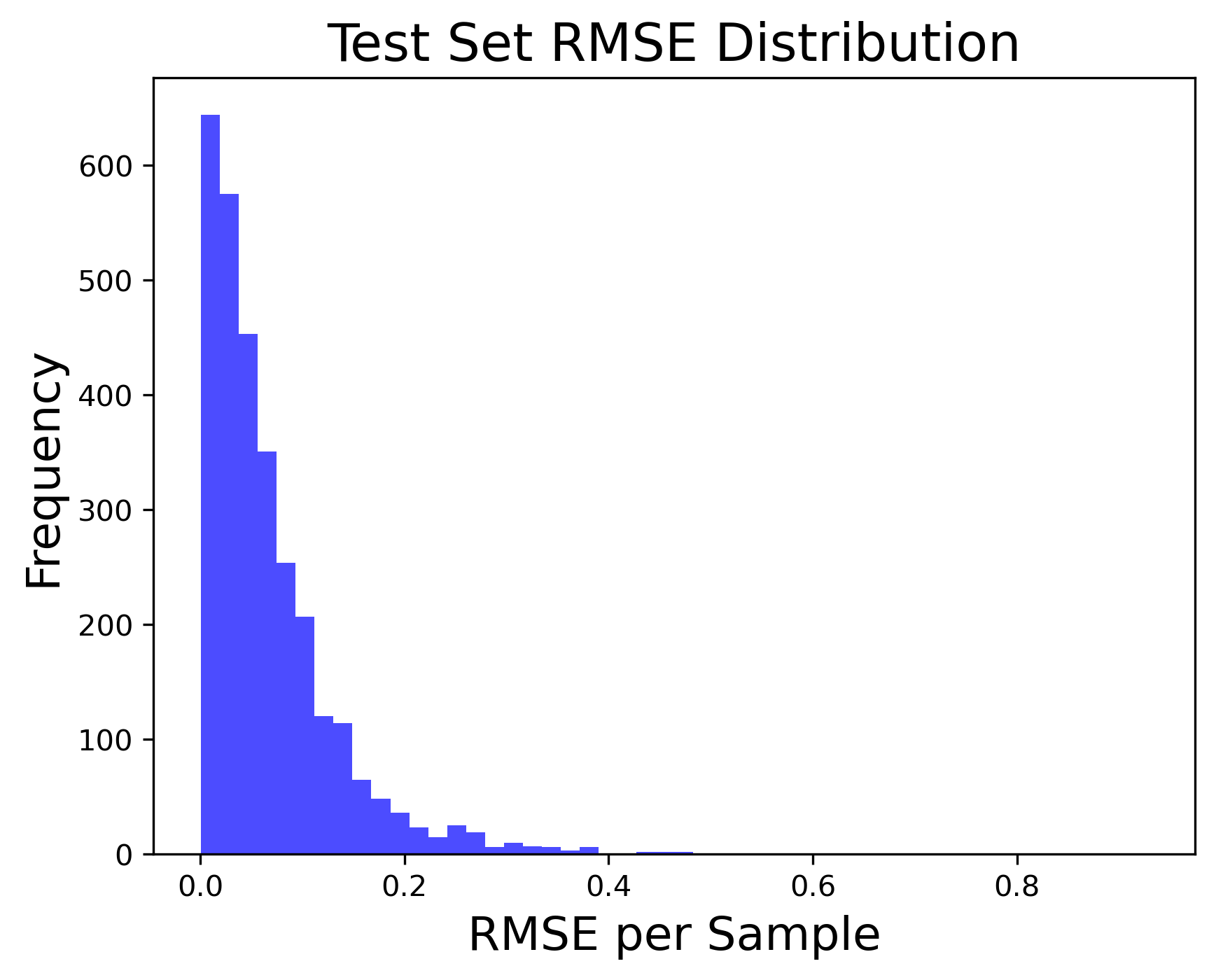} 
  \captionof{figure}{Distribution of the RMSE per sample in the test set of the kite-shaped obstacles.}
  \label{figAp7}
\end{minipage}

\subsection{Inverse problem for star-shaped obstacles}
\label{A4}

\subsubsection{Fixed impedance}
\label{A41}

\medskip
\begin{table}[h!]
\caption{Star-shaped (fixed $\lambda$) obstacle regression: Initial approach with poor performance.}
\centering
\begin{tabular}{@{}lccc}
\br
Metric               & Train & Valid & Test     \\
\mr
$R^2$ score (\%)     & 87.18    & 85.21      & 85.63    \\
RMSE                 & 0.0617   & 0.0667     & 0.0657   \\
\br
\end{tabular}
\label{tabAp6}
\end{table}

\begin{table}[h!]
\caption{Circular CNN architecture for star-shaped (fixed $\lambda$) obstacles.}
\begin{tabular}{@{}lp{0.7\textwidth}}
\br
Layer & Configuration \\
\mr
Input & $(T_0=128, C_0=4)$ tensor of real/imaginary electric/magnetic fields \\
CircularConv1D$^{(1)}$ & $N_f^{(1)}=128$ filters, $K^{(1)}=5$, $S^{(1)}=1$, Swish activation; output: $(T_1=128, C_1=128)$ \\
CircularConv1D$^{(2)}$ & $N_f^{(2)}=128$ filters, $K^{(2)}=5$, $S^{(2)}=2$, Swish activation; output: $(T_2=64, C_2=128)$ (downsampling) \\
CircularConv1D$^{(3)}$ & $N_f^{(3)}=128$ filters, $K^{(3)}=15$, $S^{(3)}=1$, Swish activation; output: $(T_3=64, C_3=128)$ \\
CircularConv1D$^{(4)}$ & $N_f^{(4)}=128$ filters, $K^{(4)}=31$, $S^{(4)}=1$, Swish activation; output: $(T_4=64, C_4=128)$ \\
Pointwise Conv1D & $N_{\mathrm{b}}=64$ filters, $K=1$, Swish activation; output: $(T_4, N_{\mathrm{b}}=64)$ (bottleneck) \\
Flatten & $d_0 = T_4 \times N_{\mathrm{b}} = 4096$ \\
Dense & $d_1=256$ nodes, Swish activation; LayerNormalization, Dropout ($p=0.1$) \\
Dense & $d_2=128$ nodes, Swish activation; LayerNormalization \\
Output Dense & $d_{\mathrm{out}}=13$ units, no activation\\
\mr
Task Type & Regression\\
Optimizer & Adam, learning rate $\alpha=1\times10^{-4}$ \\
Loss Function & Mean Squared Error (\(\mathcal{L}_{\mathrm{reg}}\)) \\
Metrics & Root Mean Squared Error (RMSE) \\
Training & Batch size $=128$; EarlyStopping: $\Delta=10^{-4}$, Patience$=200$; restore best weights; 33 sec/epoch on GPU \\
\br
\end{tabular}
\label{tabAp7}
\end{table}

\begin{figure}[h!]
    \centering        
    \includegraphics[width=0.48\textwidth]{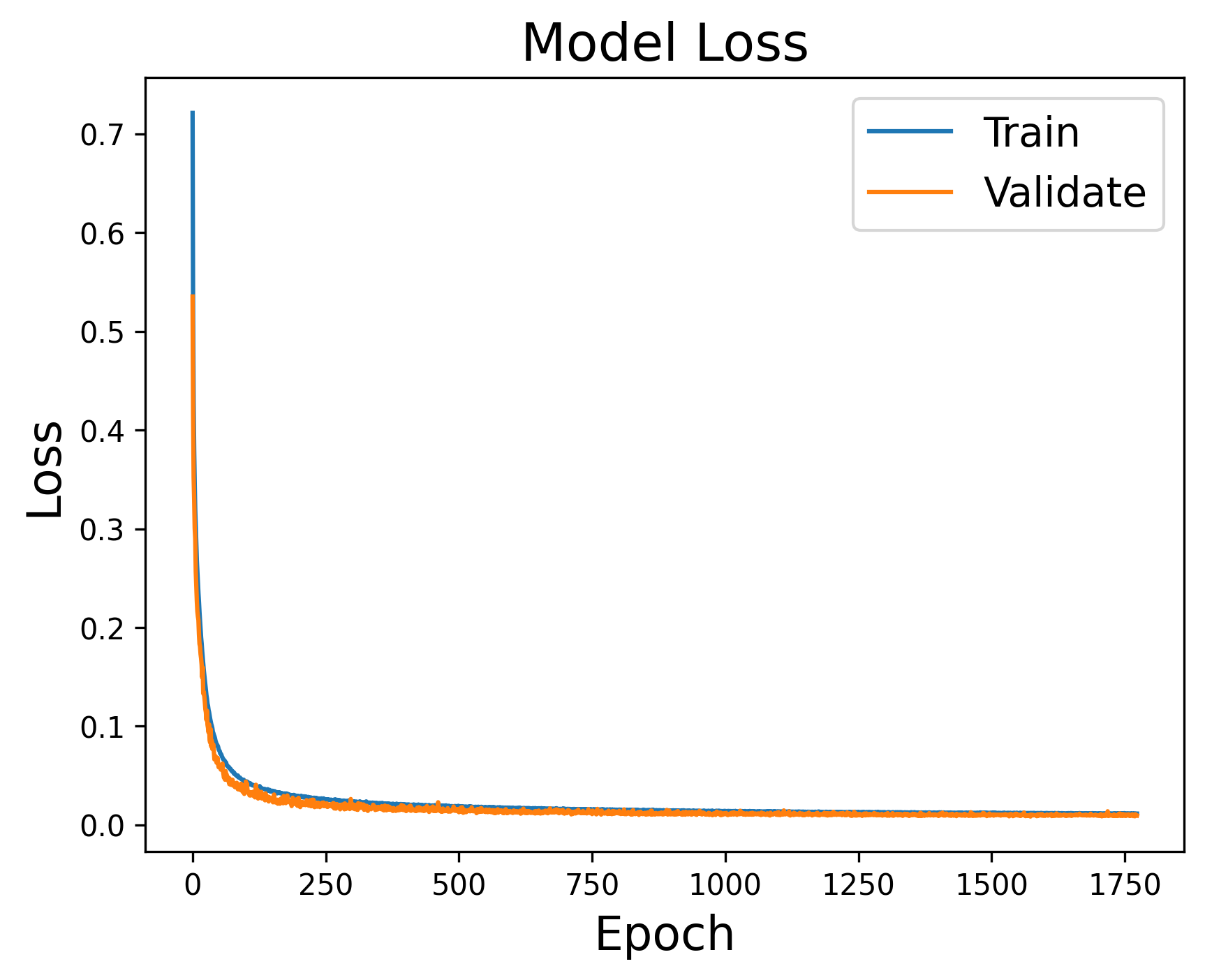} 
    \caption{The training and validation loss of the network on  star-shaped (fixed $\lambda$) obstacles.}
    \label{figAp8}
\end{figure}
\clearpage

\begin{minipage}{0.45\textwidth}
  \captionsetup{type=table}
  \captionof{table}{\label{tabAp9}Test set performance by coefficient for star-shaped (fixed $\lambda$) obstacles.}
  \centering
    \begin{tabular}{@{}lcc}
      \br
      Coefficient & $R^2$ score (\%) & RMSE    \\
      \mr
      $\alpha_0$  & 99.90             & 0.0019  \\
      $\alpha_1$  & 98.24             & 0.0230  \\
      $\alpha_2$  & 99.29             & 0.0146  \\
      $\alpha_3$  & 99.37             & 0.0136  \\
      $\alpha_4$  & 99.26             & 0.0149  \\
      $\alpha_5$  & 99.01             & 0.0172  \\
      $\beta_1$   & 97.13             & 0.0292  \\
      $\beta_2$   & 99.32             & 0.0143  \\
      $\beta_3$   & 99.38             & 0.0138  \\
      $\beta_4$   & 99.24             & 0.0150  \\
      $\beta_5$   & 98.91             & 0.0181  \\
      $x_0$       & 99.89             & 0.0039  \\
      $y_0$       & 99.87             & 0.0041  \\
      \br
    \end{tabular}
\end{minipage}
\hfill
\begin{minipage}{0.48\textwidth}
    \centering        
    \includegraphics[width=\textwidth]{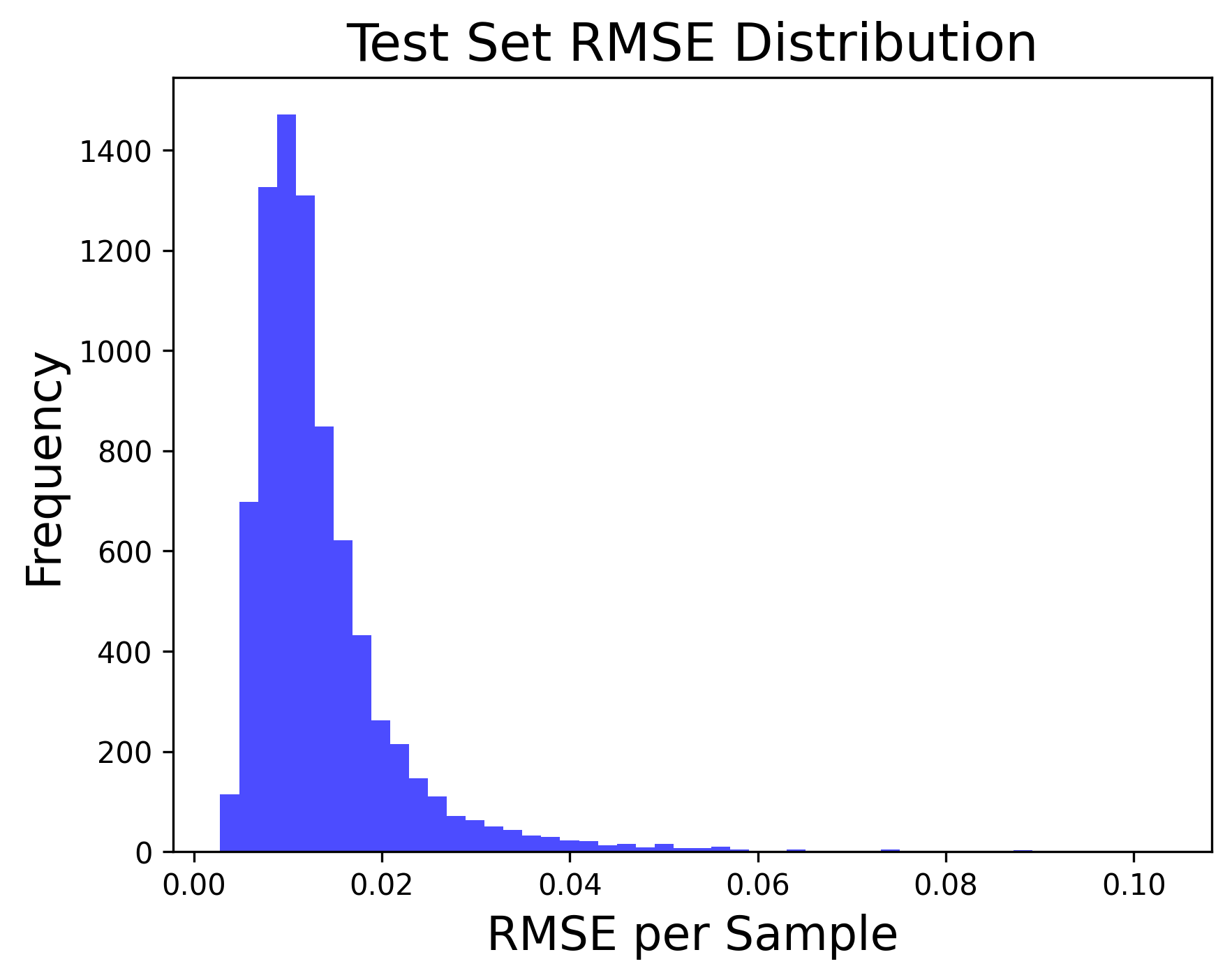} 
    \captionof{figure}{Distribution of the RMSE per sample in the test set of the star-shaped (fixed $\lambda$) obstacles.}
    \label{figAp9}
\end{minipage}%
\subsubsection{Variable impedance}
\label{A42}

\medskip
\begin{table}[h!]
\caption{Star-shaped obstacle regression: Initial approach with poor performance.}
\centering
\begin{tabular}{@{}lccc}
\br
Metric               & Train & Valid & Test   \\
\mr
$R^2$ score (\%)     & 98.71    & 90.05      & 90.05  \\
RMSE                 & 0.0531   & 0.0922     & 0.0915 \\
\br
\end{tabular}
\label{tabAp8}
\end{table} 
\medskip

\begin{table}[h!]
\caption{Circular CNN architecture for star-shaped obstacles.}
\begin{tabular}{@{}lp{0.7\textwidth}}
\br
Layer & Configuration \\
\mr
Input & $(T_0=128, C_0=8)$ tensor of real/imaginary electric and magnetic fields, from two incident waves \\
CircularConv1D$^{(1)}$ & $N_f^{(1)}=128$ filters, $K^{(1)}=5$, $S^{(1)}=1$, Swish activation; output: $(T_1=128, C_1=128)$ \\
CircularConv1D$^{(2)}$ & $N_f^{(2)}=128$ filters, $K^{(2)}=5$, $S^{(2)}=2$, Swish activation; output: $(T_2=64, C_2=128)$ (downsampling) \\
CircularConv1D$^{(3)}$ & $N_f^{(3)}=128$ filters, $K^{(3)}=15$, $S^{(3)}=1$, Swish activation; output: $(T_3=64, C_3=128)$ \\
CircularConv1D$^{(4)}$ & $N_f^{(4)}=128$ filters, $K^{(4)}=31$, $S^{(4)}=1$, Swish activation; output: $(T_4=64, C_4=128)$ \\
Angular Attention & Spatial mixing: CircularConv1D($K_{\text{mix}}=3$, Swish); LayerNormalization; Channel attention: SE-style (reduction $r=8$) \\
Pointwise Conv1D & $N_{\mathrm{b}}=64$ filters, $K=1$, Swish activation; output: $(T_L, N_{\mathrm{b}}=64)$ (bottleneck) \\
Flatten & $d_0 = T_L \times N_{\mathrm{b}} = 4096$ \\
Dense & $d_1=512$ nodes, Swish activation; LayerNormalization, Dropout ($p=0.3$), $L2$ regularization ($\lambda_{\mathrm{reg}}=10^{-4}$) \\
Dense & $d_2=256$ nodes, Swish activation; LayerNormalization, Dropout ($p=0.2$), $L2$ regularization ($\lambda_{\mathrm{reg}}=10^{-4}$) \\
Dense & $d_3=128$ nodes, Swish activation; LayerNormalization, Dropout ($p=0.1$), $L2$ regularization ($\lambda_{\mathrm{reg}}=10^{-4}$) \\
Output Dense & $d_{\mathrm{out}}=14$ units, no activation \\
\mr
Task Type & Regression \\
Optimizer & Adam, learning rate $\alpha=5\times10^{-5}$, gradient clipping $\gamma=1.0$ \\
Loss Function & Mean Squared Error (\(\mathcal{L}_{\mathrm{reg}}\)) \\
Metrics & Root Mean Squared Error (RMSE) \\
Training & Batch size $=128$; EarlyStopping: $\Delta=10^{-4}$, Patience$=200$; restore best weights; 42 sec/epoch on GPU \\
\br
\end{tabular}
\label{tabAp10}
\end{table}

\noindent
\begin{minipage}[t]{0.48\textwidth}
    \centering        
    \includegraphics[width=\textwidth]{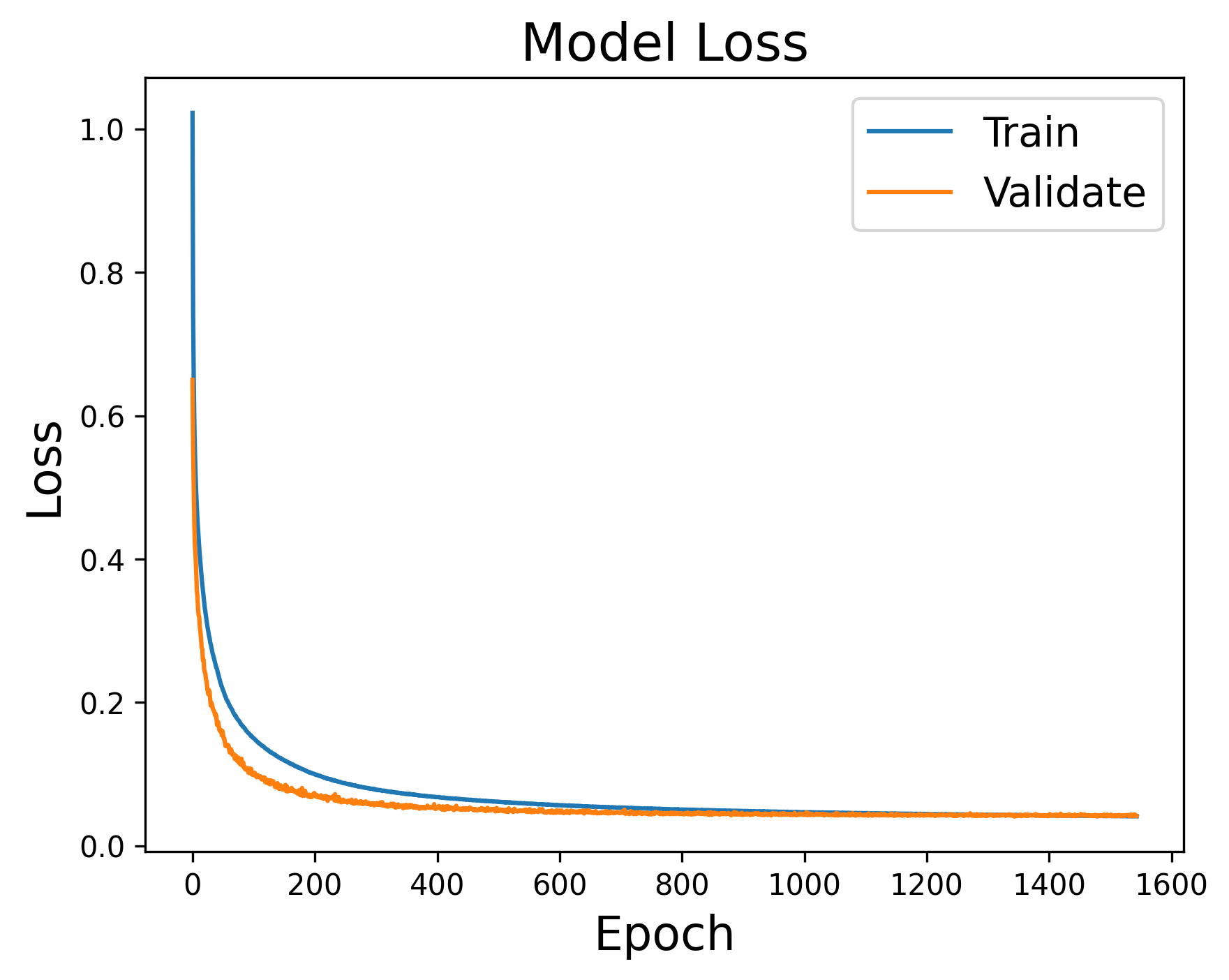} 
    \captionof{figure}{The training and validation loss of the network on star-shaped obstacles.}
    \label{figAp10}
\end{minipage}%
\hfill
\begin{minipage}[t]{0.48\textwidth}
    \centering        
    \includegraphics[width=\textwidth]{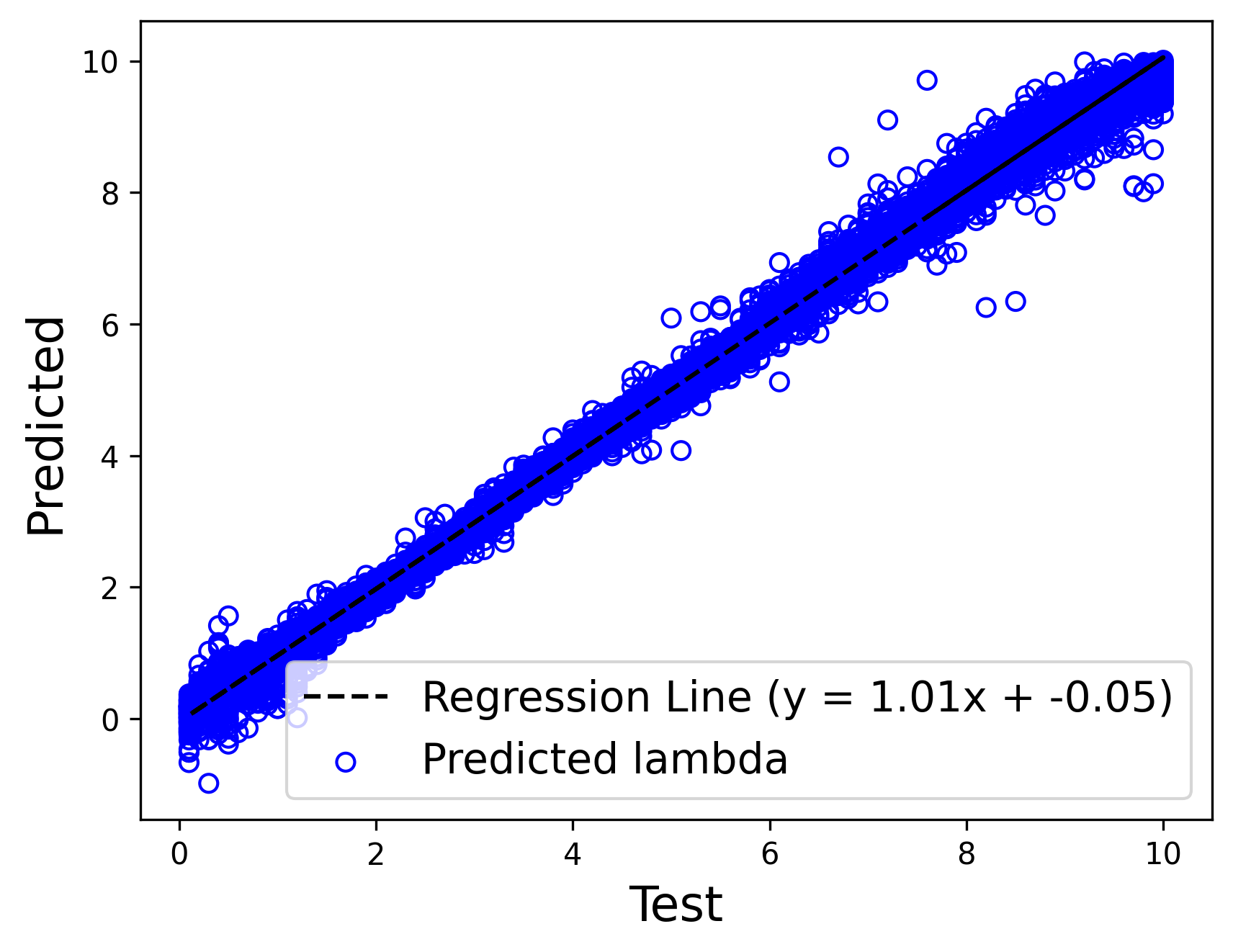} 
    \captionof{figure}{Regression results for the estimated impedance \(\lambda\) of star-shaped obstacles.}
    \label{figAp11}
\end{minipage}

\clearpage
\noindent
\begin{minipage}{0.48\textwidth}
  \captionsetup{type=table}
  \captionof{table}{\label{tabAp11}Test set performance by coefficient for star-shaped obstacles.}
  \centering
    \begin{tabular}{@{}lcc}
      \br
      Coefficient & $R^2$ score (\%) & RMSE    \\
      \mr
      $\alpha_0$  & 99.68             & 0.0034  \\
      $\alpha_1$  & 91.39             & 0.0506  \\
      $\alpha_2$  & 96.28             & 0.0337  \\
      $\alpha_3$  & 96.85             & 0.0309  \\
      $\alpha_4$  & 96.11             & 0.0341  \\
      $\alpha_5$  & 93.14             & 0.0452  \\
      $\beta_1$   & 89.55             & 0.0560  \\
      $\beta_2$   & 96.09             & 0.0343  \\
      $\beta_3$   & 96.81             & 0.0310  \\
      $\beta_4$   & 96.35             & 0.0332  \\
      $\beta_5$   & 93.31             & 0.0449  \\
      $x_0$       & 99.61             & 0.0072  \\
      $y_0$       & 99.57             & 0.0076  \\
      $\lambda$   & 99.56             & 0.1923  \\
      \br
    \end{tabular}
\end{minipage}%
\hfill
\begin{minipage}{0.48\textwidth}
    \centering        
    \includegraphics[width=\textwidth]{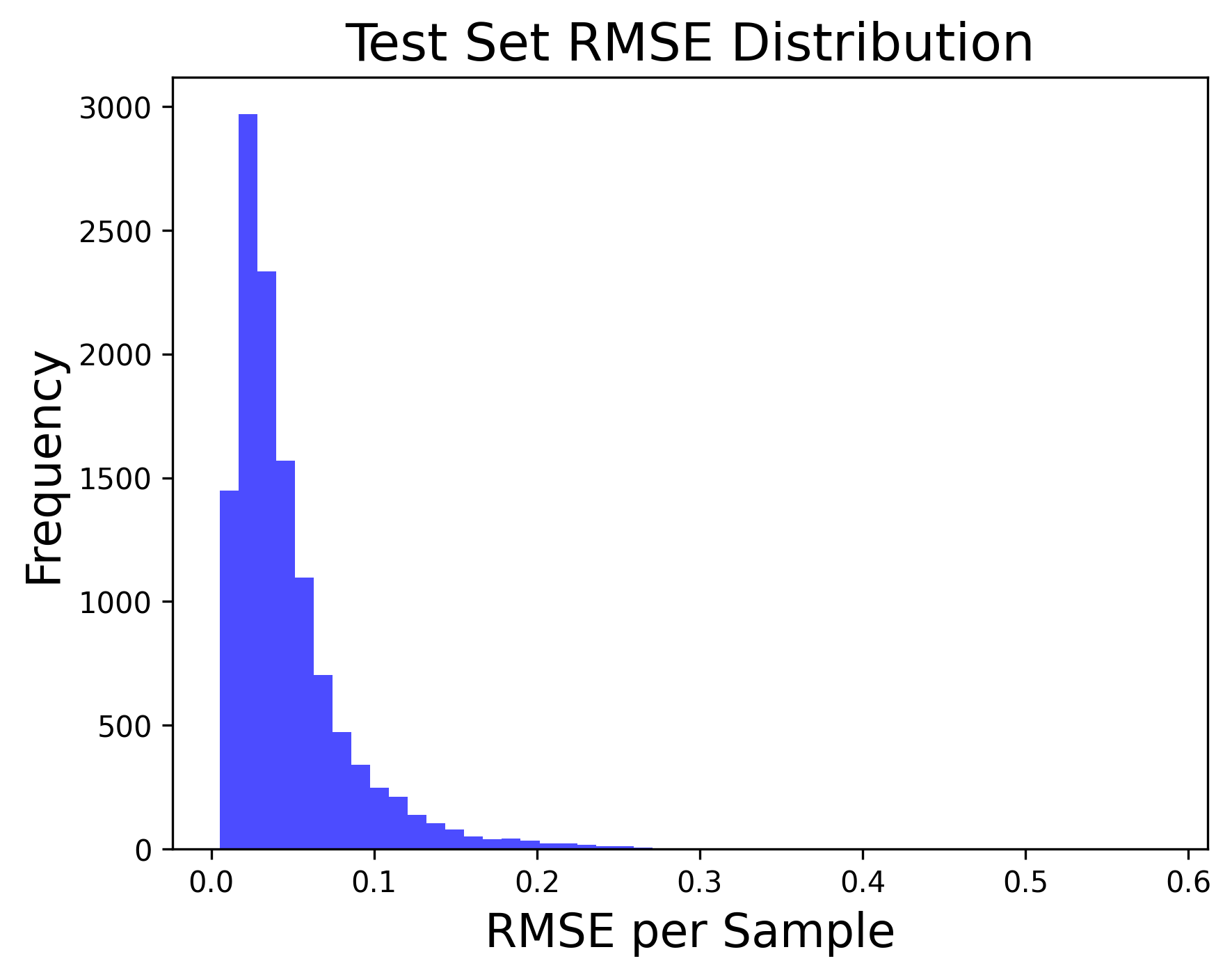} 
    \captionof{figure}{Distribution of the RMSE per sample in the test set of the star-shaped obstacles.}
    \label{figAp12}
\end{minipage}
\end{document}